\newif\iflong
\longtrue   % arXiv extended version
\iflong%
\documentclass{article}
\usepackage{graphicx} % Required for inserting images
\usepackage{subcaption}
\usepackage{amsmath,amsthm,amssymb}
\usepackage{dsfont}
\usepackage{xcolor}
\usepackage[margin=1in]{geometry}
\newcommand*{\dims}{d}
\newcommand{\calN}{\mathcal{N}}

\newcommand{\calV}{\mathcal{V}}
\newcommand{\calG}{\mathcal{G}}
\newcommand{\calE}{\mathcal{E}}
\newcommand{\calM}{\mathcal{M}}
\DeclareMathOperator*{\diag}{diag}
\DeclareMathOperator*{\diam}{diam}

\newcommand*{\N}{\ensuremath{\mathds{N}}}

\newcommand*{\R}{\ensuremath{\mathbb{R}}}

\newcommand{\domain}{\Omega}

\newcommand{\cons}{c}

\DeclareMathOperator*{\argmax}{argmax}
\DeclareMathOperator*{\argmin}{argmin}

\newcommand{\param}{\theta}

\newcommand{\E}{\mathbb{E}}

\newcommand{\ind}[1]{\hat\imath(#1)}

\newcommand{\M}{\mathcal{M}}

\newcommand{\eps}{\varepsilon}

\definecolor{sky}{rgb}{0.31, 0.69, 0.99}
\definecolor{apple}{rgb}{0.13, 0.6, 0.14}
\definecolor{grape}{rgb}{0.99, 0.35, 0.34}
\definecolor{hpurp}{rgb}{0.87, 0., 1.}
\definecolor{orong}{rgb}{1.0, 0.55, 0.0}
\definecolor{ponk}{rgb}{.8, 0.08, 0.58}
\definecolor{flamingopink}{rgb}{0.99, 0.46, 0.47}
\definecolor{teal}{rgb}{0.0, 0.5, 0.5}
\definecolor{brinkpink}{rgb}{0.98, 0.38, 0.5}
\definecolor{carnationpink}{rgb}{1.0, 0.65, 0.79}

\newcommand{\miparam}{\theta}

\newcommand{\FundingLogos}{%
  \raisebox{0pt}{\includegraphics[height=1.5cm]{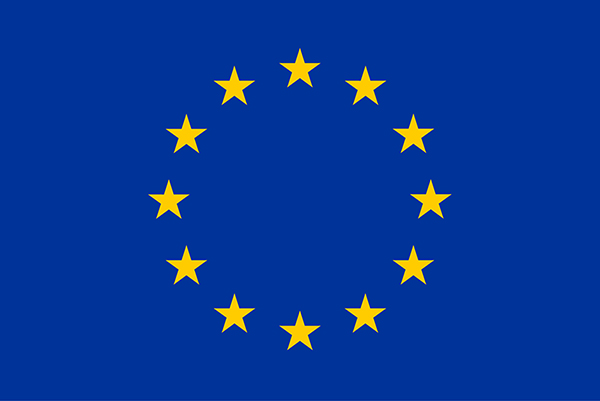}}%
  \hspace{1em}%
  \raisebox{0pt}{\includegraphics[height=1.5cm]{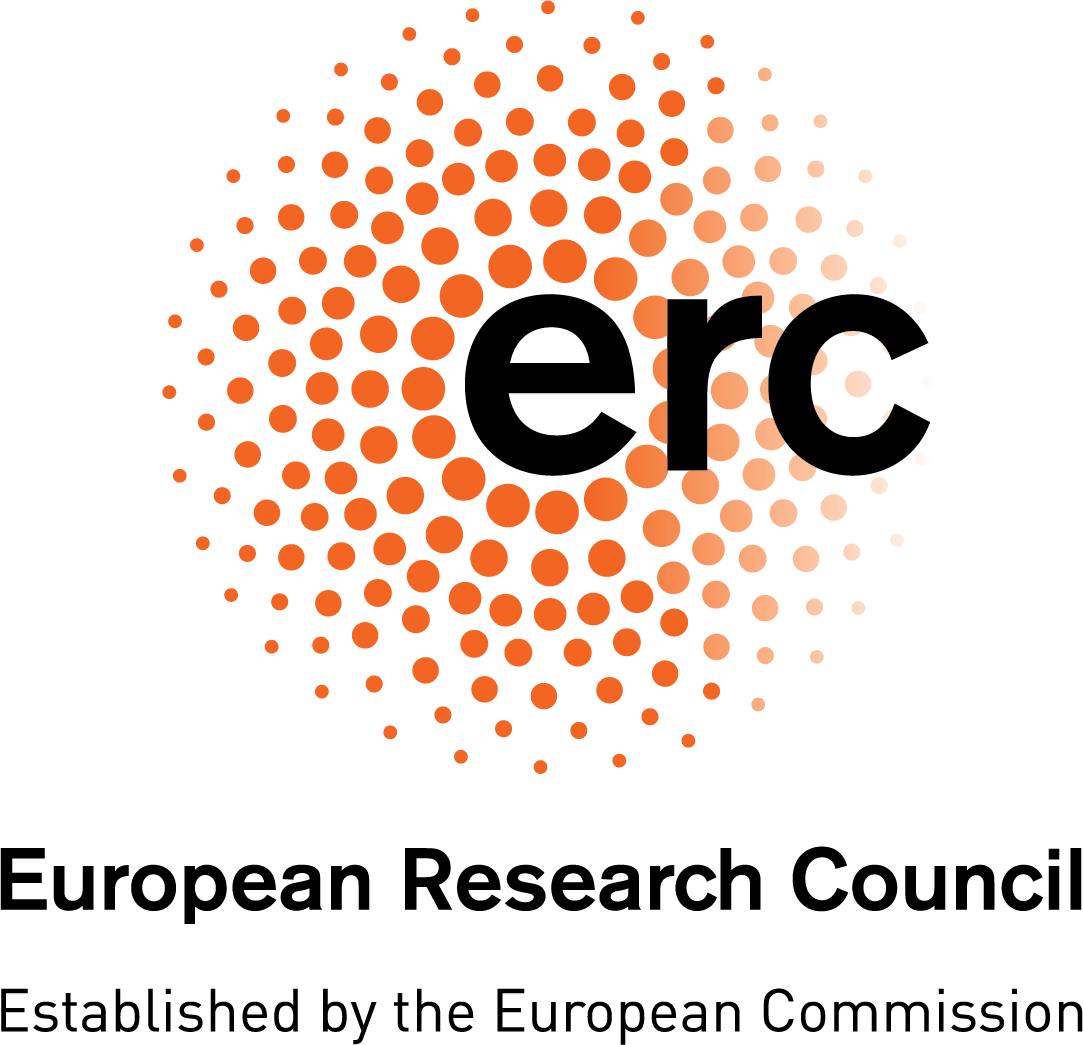}}%
}

\usepackage{tcolorbox}
\tcbset{colback=white,boxrule=0.5pt,arc=2pt,left=4pt,right=4pt,top=4pt,bottom=4pt}
\usepackage{todonotes}
\definecolor{cgreen}{rgb}{0.1, .5, 0.1}
\usepackage[
colorlinks=true,
linkcolor=blue,
filecolor=magenta,
citecolor=cgreen,
urlcolor=cyan,
]{hyperref}	
\usepackage{enumitem}
\usepackage{tikz}
\usepackage{booktabs}

\usepackage{algorithm}
\usepackage[noend]{algpseudocode}
\algblock{ParFor}{EndParFor}
\algnewcommand\algorithmicparfor{\textbf{for}}
\algnewcommand\algorithmicpardo{\textbf{do in parallel}}
\algnewcommand\algorithmicendparfor{\textbf{end\ parfor}}
\algrenewtext{ParFor}[1]{\algorithmicparfor\ #1\ \algorithmicpardo}
\makeatletter
\ifthenelse{\equal{\ALG@noend}{t}}%
  {\algtext*{EndParFor}}
  {}%
\makeatother

\usepackage[capitalize,nameinlink,noabbrev]{cleveref}

\setlist[enumerate,1]{label=(\roman*), ref=(\roman*)}

\newtheorem{theorem}{Theorem}[section]
\newtheorem{lemma}[theorem]{Lemma}
\newtheorem{assumption}[theorem]{Assumption}

\newtheorem{remark}[theorem]{Remark}

\newif\iflong
\longtrue   % arXiv extended version
\makeatletter
\newcommand{\longv}[1]{\iflong{#1}\else\@bsphack\@esphack\fi}
\newcommand{\shortv}[1]{\iflong\@bsphack\@esphack\else #1\fi}
\makeatother

\title{Consensus of Consensus: Decentralized Consensus-Based Optimization}
\author{S. Bellavia\thanks{Dipartimento di Ingegneria Industriale, Università degli studi di Firenze, Firenze, Italy. Member of the INdAM Research Group GNCS.}\and % 
M. Fornasier\thanks{School of Computation, Information and Technology, Technical University of Munich, Germany.} %
\textsuperscript{,}%
\thanks{Munich Center for Machine Learning, Munich, Germany.} %
\and %
G. Malaspina\footnotemark[1] \and %
T. Roith\footnotemark[2] \textsuperscript{,}\footnotemark[3]}
\date{November 2026}%
\fi%

\begin{document}

\maketitle%

\iflong%
\begin{abstract}
We propose and analyze a decentralized variant of Consensus-Based Optimization (CBO) for the global minimization of continuous, possibly non-convex functions accessible only through a stochastic, possibly biased, zeroth-order oracle with bounded mean-squared error. The agents are split among $N$ units of an undirected communication network. Each unit evolves its swarm of $M$ agents towards a local tracking variable, which follows the network average of the local consensus points through one round of neighbour-to-neighbour communication via a doubly stochastic mixing matrix, without any central coordinator. We show that, under an explicit condition linking the drift and diffusion parameters to $M$, all agents reach consensus almost surely, at an exponential mean-square rate governed by the spectrum of the mixing matrix. Using a finite-agent quantitative Laplace principle, we prove that, for large $\alpha$ and provided some agent visits a neighbourhood of the global minimizer with sufficient probability, the expected mean-squared distance of the swarm from the minimizer falls below any tolerance within an explicit number of iterations, at a rate that deteriorates for poorly connected networks. We verify this hitting assumption for a simplified dynamics. Experiments on benchmarks up to dimension $d=1000$, with stochastic objectives, and on a federated learning task demonstrate the method's effectiveness.
\end{abstract}
%
% \else%
% \input{abstract_SIAM}%
\fi%

% REQUIRED
%\begin{keywords}

% REQUIRED
%\begin{MSCcodes}
%\end{MSCcodes}
%\fi%

\section{Introduction}

Consensus-based optimization (CBO) is a versatile, zero-order metaheuristic for
global optimization: an ensemble of interacting agents explores the search
space and is driven, through a drift-diffusion dynamics, toward a weighted
average of the agents' positions, the \emph{consensus point}, computed
via a Gibbs-type weighting of the objective function. Since it never requires
gradient information, CBO is directly applicable to nonsmooth, nonconvex, and
even black-box objectives, while at the same time being amenable to rigorous
mean-field and finite-agent analysis, unlike most other metaheuristics of
swarm-intelligence type. The model was first introduced by Pinnau, Totzeck,
Tse and Martin~\cite{pinnau2017}, who proposed the underlying stochastic
differential equation and derived its mean-field limit as a nonlocal,
degenerate parabolic PDE. The first rigorous proofs that the method actually
converges to a global minimizer, under quantitative assumptions on the
objective and the initialization, were given by Carrillo, Choi, Totzeck and
Tse~\cite{carrillo2018} and, in fully unconditional 
form with respect to initialization and with agent analysis as well, by Fornasier, Klock and Riedl~\cite{fornasier2024} (see also the
independent proof by Ha, Jin and Kim~\cite{hajinkim2020}). Since these
foundational results, a large body of work has refined the model along
several directions: anisotropic, component-wise noise to remove the
dimension-dependence of the drift and to reduce computational cost via
mini-batching~\cite{carrillojinlizhu2021}, truncated
noise~\cite{fornasierrichtarikriedlsun2025}, jump-diffusion
noise~\cite{kalisesharmatretyakov2023}, memory
effects~\cite{riedl2023memory}, heterogeneous noise combined with random-batch
agent interactions~\cite{ko2022}, noisy/estimated objective
evaluations~\cite{bellavia2025}, and mechanisms such as a personal-best
memory term borrowed from particle swarm
optimization~\cite{totzeckwolfram2020}. In parallel, CBO has been extended to
constrained problems on the sphere and other
hypersurfaces~\cite{fornasier2020hypersurfaces} and shown to be competitive on
high-dimensional machine learning
tasks~\cite{carrillojinlizhu2021,fornasier2021spherejmlr}, has been connected
to consensus-based sampling~\cite{carrillohoffmannstuartvaes2022}, and has
recently been introduced to robotics for global trajectory and policy
optimization~\cite{sun2026robotics}; we refer to~\cite{totzeck2021trends} for
a broader survey of these developments.

Yet, in order to fully exploit the potential of CBO, in particular on the
large-scale, high-dimensional problems arising in modern machine learning, 
a genuinely parallel and distributed formulation of the algorithm would be
crucial. It is a recurring claim in the literature that CBO-type dynamics are
``highly parallelizable'' or ``inherently
scalable''~\cite{totzeck2021trends,fornasierhuangklemencmalaspina2026}, and
indeed the per-agent updates may typically be distributed across the cores of
a GPU or a multi-core machine~\cite{sun2026robotics}. However, this common
notion of parallelism is architecturally \emph{centralized}: every agent,
regardless of which processor evaluates its update, relies on the very same,
globally shared consensus point, so that a central synchronization step is
required at every iteration and no genuine network structure or local-only
communication is exploited. We  recently explored the parallelization of the centralized CBO in the companion paper \cite{bellavia2026parallelCBO}. A notable departure from a single, globally
shared consensus point is given by the kinetic theory-based optimization (KBO)
methods of~\cite{benfenati2022binary}, where the mean-field interaction of CBO
is replaced by binary interactions between randomly sampled pairs of
agents, described by a Boltzmann-type equation. There, however, any pair of
agents may interact and the update still involves a collective estimate
over the whole population, so that no network structure restricting
communication is present. To the best of our knowledge, no
distributed formulation of CBO with computational units communicating solely over a
network of local connections has appeared in the literature, and, with the
possible exception of the federated-learning variants of CBO proposed by
Garc\'ia Trillos and
collaborators~\cite{carrillogarciatrillos2024fedcbo,garciatrilloscb2o2024,garciatrillosdefending2025},
which nonetheless rely on a central server orchestrating the interaction
between clients, no convergence proof has been provided for any
parallel, decentralized implementation of the method.

In this paper we explore a parallel and distributed formulation of CBO with
fully decentralized computation. We assume the $N$ computational units to be
the nodes of a connected, undirected communication network $\calG=(\calV,\calE)$, each
holding a local swarm of agents and able to communicate only with its
network neighbors. Rather than relying on a single, globally shared consensus
point, each node computes its own \emph{local} consensus point from its own
agents and then updates a local estimate of the global consensus by
averaging it, through a doubly stochastic mixing matrix, with the estimates
held by its neighbors, a mechanism we call \emph{consensus of consensus},
in the spirit of classical linear gossip and distributed-averaging schemes
over networks~\cite{xiaoboyd2004,nedicozdaglar2009}. Since the network is
connected, information about the location of the global minimizer, although
initially known only locally and imprecisely at each node, is able to reach
every other node; along the dynamics this information is slowly diffused
across the network, and, as we prove, global consensus among all agents
and all local estimates is asymptotically formed, exactly as in the
centralized method, without any node ever having direct access to the whole
swarm.

A further consequence of this construction is that each node's local
consensus point is computed using only the agents, and hence only the
(possibly stochastic) evaluations of the objective function, held locally
at that node. This is particularly relevant for objectives of the type
arising in machine learning, where the objective function is itself an
empirical risk defined over possibly locally held data. The proposed scheme thus
allows for a fully decentralized optimization of such objectives, in which
nodes only ever exchange the low-dimensional consensus vectors $Z^i$,
and never share gradients, model updates, or any other disaggregated
information about their local data from which the data itself could
conceivably be reverse-engineered, as is instead possible for
gradient-sharing schemes~\cite{geiping2020invertinggradients,haim2022reconstructing}.
In this sense, the method is privacy-respecting essentially by design. This
sets it apart from existing CBO-inspired approaches to distributed learning,
such as FedCBO~\cite{carrillogarciatrillos2024fedcbo} and its bilevel
extensions CB$^2$O~\cite{garciatrilloscb2o2024} and its
attack-robust variant~\cite{garciatrillosdefending2025}, which use CBO
dynamics to reach (possibly clustered) consensus among the models of
different clients but still operate within the classical federated-learning
architecture, with a central server aggregating information from all
clients at every round\footnote{In clustered federated learning, many computational units (for example, smartphones) each hold private data and their own model, and the units are secretly divided into groups whose data come from the same distribution. Nobody knows which group anyone belongs to, or even how many groups there are. The goal is for every unit to end up with a model that works well for its own group. FedCBO treats each unit's model as an agent. At every round, each unit trains its model on its own data, downloads a sample of other units' models, and tests them on its own data. It then moves its model toward a weighted average of those models, where models with lower loss get exponentially more weight. Models from other groups perform poorly on the unit's data, so they receive almost no weight, while models from the same group dominate the average. As a result, units in the same group automatically converge to a shared model for that group, without any data being shared and without the cluster memberships ever being revealed.}. Our fully decentralized, server-free formulation
removes this single point of coordination, failure, and potential attack,
while retaining, as we show,  provable convergence guarantees; this
suggests it as a promising gradient-free building block for privacy-aware
decentralized training of machine learning and deep learning models over
peer-to-peer networks.

\paragraph{Organisation of the paper}
The rest of the paper is organised as follows. In
Section~\ref{sec:decentralized} we state our assumptions on the objective
function and on the stochastic oracle, and we introduce the decentralized
parallel CBO scheme together with the tracking mechanism that replaces the
global consensus point. Section~\ref{sec:conform} is devoted to consensus
formation. We show that, under an explicit condition on $\lambda$, $\sigma$
and $M$, all agents and tracking variables converge almost surely to a
common point, at a rate determined by the spectrum of the mixing matrix and by the parameters of the local dynamics.
In Section~\ref{sec:convergence} we derive a finite-agent quantitative
Laplace principle and use it to prove convergence in expectation of the swarm
towards the global minimizer. Section~\ref{sec:hitting} discusses the
hitting-probability assumption required by this result and verifies it for a
simplified decoupled dynamics. Finally, Section~\ref{sec:numerics} reports
numerical experiments illustrating the behaviour of the method.

\subsection{Notation}

Throughout the paper, $\N=\{0,1,2,\dots\}$ and, for $n\in\N$, the index range $i=1:n$ stands for $i=1,\dots,n$. 

Capital letters $X$, $Z$, $C$ denote the random vectors in $\mathbb{R}^d$ generated by the method (agents, tracking variables and consensus points, respectively), and the corresponding lowercase letters denote their components. As a general rule, superscripts are labels identifying an object (e.g. agent $\ell$), while subscripts refer to iteration indices and components. For instance, $x^\ell_{k,s}$ is the $s$-th component of the position $X^\ell_k$ of agent $\ell$ at iteration $k$. Calligraphic letters, such as $\mathcal{X}^s_k$ and $\mathcal{Z}^s_k$, denote the vectors collecting the $s$-th components across all agents or all computational units, so that, e.g., $(\mathcal{Z}^s_k)_i = z^i_{k,s}$. Capital letters are also used for matrices. We write $I$ for the identity matrix, $\mathbf{1}_N$ for the vector of all ones in $\R^N$, $\mathbf{e}_s$ for the $s$-th vector of the canonical basis of $\R^\dims$, and $\diag(\cdot)$ for a (block-)diagonal matrix with the given diagonal entries or blocks. The symbol $\Lambda$, with subscripts specifying the quantity being bounded, denotes a deterministic constant multiplying a geometric decay rate in an upper bound.

On $\R^\dims$ we denote by $\|\cdot\|$ the Euclidean norm. For vectors of any length, $\|\cdot\|_2$ and $\|\cdot\|_\infty$ denote the Euclidean and the maximum norm; for matrices, $\|\cdot\|_2$ is the induced spectral norm and $\rho(\cdot)$ the spectral radius. For a finite set $\mathcal{S}$, $|\mathcal{S}|$ denotes its cardinality.
All random variables are defined on a probability space $(\Omega,\mathcal{A},\mathbb{P})$. We denote by $\E$ the expectation and by $\E[\,\cdot\,|\,\mathcal{S}]$ the expectation conditioned on an event $\mathcal{S}$, and ``i.i.d.'' stands for independent and identically distributed.

 %Subscripts on constants such as $D_{\mathcal{X},s}$ and $\Lambda_{\mathcal{Z},s}$ specify the quantity being bounded.
\section{Decentralized parallel CBO}\label{sec:decentralized}
We consider a global optimization problem of the form 
$$
x^* = \argmin_{x\in\R^\dims} f(x)
$$
where $f:\R^\dims\longrightarrow\R$ is a continuous, generally nonconvex function.
Moreover, we assume that the function $f$ cannot be evaluated exactly. That is, given a point $x\in\R^\dims$, we assume that only a stochastic estimator $ \widehat f(x,\omega) \simeq f(x)$ can be computed, where 
$\omega\in\Omega$ for a suitable probability space $(\Omega, \mathcal{A}, \mathbb{P})$. Here we assume that the global minimizer $x^*$ is unique. 
Given a scalar $r>0$ we denote
\begin{equation}\label{def:B_r}
    \begin{aligned}
        &\mathcal B_r:=\{x\in\R^d\enspace|\enspace \|x-x^*\|\leq r \},\quad f_r = \max_{x \in \mathcal B_r} f(x).
    \end{aligned}
\end{equation}
We make the following assumptions on the objective function $f$ and the stochastic estimator $\widehat f$.
\begin{assumption}\label{a:f} 
\begin{enumerate}[label=(\roman*)]
\item $f:\R^d\longrightarrow\R$ is continuous and has a unique global minimizer $x^* = \argmin_{x\in\R^d} f(x)$.
\item There exist constants $f_{\infty}, R_0, \beta, \nu>0$ such that
\begin{align*}
(f(x)-f(x^*))^{\nu}&\geq \beta\|x-x^*\|\enspace \text{for every}\ x\in \mathcal B_{R_0},\\
f(x)&\geq f(x^*)+f_\infty\enspace \text{for every}\ x\in \R^\dims\setminus \mathcal B_{R_0}.
\end{align*}

\end{enumerate}

\end{assumption}
\begin{assumption}\label{a:f_hat}
There exists a constant $\sigma_f\geq 0$ such that 
\begin{equation}\label{noise_bound}
    \E_{\omega\sim \mathbb{P}}\left[|f(x)-\widehat f(x,\omega)|^2\right]\leq \sigma_f^2\enspace\text{for every}\ x\in\R^\dims.
\end{equation}
\end{assumption}
Notice that we do not assume that $\E_{\omega}[\widehat f(x,\omega)] = f(x)$, hence the oracle may not be unbiased.
To ease the notation we will also assume, without loss of generality, that $f(x^*)=0$.

\subsection{Discrete-time finite-particle regime}
Assume we have $N$ computational units on a communication network $\calG=(\cal V, \cal E)$, where $\calV=\{1,\dots,N\}$ denotes the set of vertices and $\calE$ the set of edges. Given $\calG$ we define a \emph{mixing matrix} $W\in\R^{N\times N}$ associated with the network. We make the following assumption on $\calG$ and $W$.
\begin{assumption}\label{a:W}
\begin{enumerate}[label=(\roman*)]
    \item The network $\calG = (\cal V,\cal E)$ is undirected, connected, and without multiple edges.
    \item The mixing matrix $W\in\R^{N\times N}$ is  symmetric, doubly stochastic, and such that if  $(i,j)\notin\cal E$ then $w_{ij}=0$, while $w_{ij} > 0$ whenever $(i,j) \in\calE $.
\end{enumerate}
\end{assumption}

For every node $i$ we denote by $\calN_i$ the neighborhood of $i$ in $\cal G$, namely 
$$\calN_i = \{j\in\calV\ |\ (i,j)\in\calE\}.$$ We assume that node $i$ can communicate directly only with nodes $j\in\calN_i$. 
In practice, it is natural to assume that each node also uses its own tracking variable, that is, $(i,i)\in\mathcal{E}$ and $w_{ii}>0$. However, the analysis below does not require it, so in Assumption \ref{a:W} we do not require the network to have self-loops at every node.

We assume each node $i$ holds a set of $M$ agents.\longv{\footnote{With slightly more effort and notational burden we may consider also a model where the number $M_i$ of agents may differ for different units $i$.}} Without loss of generality we can assume node $i$ holds agents $(X^\ell)_{\ell\in \calM_i}\subset\R^\dims$ with $ \calM_i = \{(i-1)M+1,(i-1)M+2, \dots,iM\}.$ \longv{Implicitly, we are ordering the agents within the units' indexes $\M_1,\dots,\M_N$.} Moreover, we assume that every node $i$ holds an additional \emph{tracking variable} $Z^i\in\R^\dims.$

Given an initial distribution $\rho_0$ with finite fourth moment, the agents evolve according to the following equations:
\begin{equation}\label{eq:CBOpar}
\begin{cases}
X^\ell_0 \overset{\text{i.i.d.}}{\sim}  \rho_0 & \text{ for } \ell\in\calM_i,\ i=1,\dots,N,\\
Z^i_0 = \widehat C^i_0 & \text{ for } i=1,\dots,N,\\
X^\ell_{k+1} = X^\ell_k - \lambda(X^\ell_k-Z^i_k)+\sigma \sum_{s=1}^\dims (x^\ell_{k,s}- \widehat\cons^i_{k,s}) \xi^\ell_{k,s}\mathbf{e}_s  & \text{ for } \ell\in\calM_i,\ i=1,\dots,N,\\
Z^i_{k+1} = \sum_{j\in\calN_i}w_{ij}Z^j_{k} +  \widehat C^i_{k+1}- \widehat C^i_{k} &  \text{ for } i=1,\dots,N,
\end{cases}
\end{equation}
where $\lambda\in(0,1)$ and $\sigma>0$, $(\xi^\ell_{k,s})$ are i.i.d. standard Gaussian random variables, and $ \widehat C^i_k$ is the consensus point for the swarm at node $i$, namely for $\alpha > 0$
\begin{equation}\label{eq:cons_i}
 \widehat C^i_k = \frac{\sum_{\ell\in\M_i} X^\ell_k \exp(-\alpha  \widehat f(X^\ell_k, \omega^\ell_k))}{\sum_{\ell\in\M_i}\exp(-\alpha  \widehat f(X^\ell_k,\omega^\ell_k))} = \sum_{\ell\in\M_i} X^\ell_k  \widehat v^{i,\ell}_k, 
 \quad  \widehat v^{i,\ell}_k = \frac{\exp(-\alpha  \widehat f(X^\ell_k,\omega^\ell_k))}{\sum_{p\in\M_i}\exp(-\alpha  \widehat f(X^p_k,\omega^p_k))}.
\end{equation}
The oracle samples $\omega^\ell_k$ are drawn independently at every evaluation. Moreover, the initial positions $(X^\ell_0)$, the noise $(\xi^\ell_{k,s})$ and the oracle samples $(\omega^\ell_k)$ are mutually independent.

\longv{We emphasize that $\sigma$ here denotes the strength of the variance in the consensus dynamics \eqref{eq:CBOpar}, and it is not to be confused with $\sigma_f$ describing the noise strength in the objective function in \eqref{noise_bound}.}

\subsubsection{Consensus formation}\label{sec:conform}
In this subsection we prove that, under suitable assumptions on the parameters of the method, the agents achieve consensus.
For any $s=1,\dots,d$, $i=1,\dots,N$ and $k\in\N$ we define
\begin{equation}\label{eq:defYU}
\begin{aligned}
    &\mathcal{X}_k^s = \left(\begin{array}{c}
         x^1_{k,s}\\
         \vdots\\
         x^{NM}_{k,s}
    \end{array}\right)\in\R^{NM}, \qquad
    \mathcal{X}_k^{s,i} = \left(x^\ell_{k,s}\right)_{\ell\in\calM_i}\in\R^{M},\\
    &\mathcal{Z}^s_k = \left(\begin{array}{c}
         z^1_{k,s}\\
         \vdots\\
         z^{N}_{k,s}
    \end{array}\right)\in\R^{N},
    \qquad
    \widehat{\mathcal{C}}^s_k = \left(\begin{array}{c}
         \widehat\cons^1_{k,s}\\
         \vdots\\
         \widehat\cons^N_{k,s}
    \end{array}\right)\in\R^{N},\quad \mathcal{H}^s_k = \mathcal{Z}^s_k-\widehat{\mathcal{C}}^s_k\in\R^N.
\end{aligned}
\end{equation}
\longv{That is, $\mathcal{X}^s_k$ collects the $s$-th component of all agents, while $\mathcal{X}^{s,i}_k$ is the block of $\mathcal{X}^s_k$ corresponding to the agents held by node $i$. Analogously, $\mathcal{Z}^s_k$ and $\widehat{\mathcal{C}}^s_k$ collect the $s$-th component of the tracking variables and the consensus points of all computational units, respectively. The vectors $\mathcal{H}^s_k$ measure the discrepancy between the tracking variables $Z^i_k$ and the consensus points $\widehat C^i_k$.}

We denote by $J$  the matrix  
\begin{equation}\label{eq:defJh}
    J = \frac{1}{N}\mathbf{1}_N\mathbf{1}_N^\top\in\R^{N\times N},
\end{equation}
and will use the following standard result about the mixing matrix $W$, which is a direct consequence of the Perron--Frobenius Theorem (see \cite[Chapter 8]{matan}).
\begin{lemma}\label{lemma:W}
Let $W$ satisfy Assumption \ref{a:W} and let us denote by $\gamma_1\geq\gamma_2\geq\dots\geq\gamma_N$ the eigenvalues of $W$. Then
\begin{enumerate}[label=(\roman*)]
    \item $W\mathbf{1}_N = \mathbf{1}_N$ and $\mathbf{1}_N^\top W = \mathbf{1}_N^\top$; in particular $JW=WJ=J$;
    \item $\gamma_1=1$ and $-1\leq\gamma_N\leq\gamma_2<1$.
\end{enumerate}
\end{lemma}

Let $\overline{Z}_k$ and $\overline{z}_{k,s}$ be the average of the tracking variables and its $s$-th component, respectively. That is, 
$$\overline{Z}_k=\frac1N\sum_{i=1}^N Z^i_k = \left(\begin{array}{c}\overline{z}_{k,1}\\ \vdots\\ \overline{z}_{k,d}\end{array}\right)\in\R^d.$$
Notice that $J\mathcal{Z}^s_k=\overline{z}_{k,s}\mathbf{1}_N$ and $(I-J)\mathcal{Z}^s_k=\mathcal{Z}^s_k-\overline{z}_{k,s}\mathbf{1}_N$. Moreover, since $Z^i_0=\widehat C^i_0$ and $\mathbf{1}_N^\top W=\mathbf{1}_N^\top$, \eqref{eq:CBOpar} implies
\begin{equation}\label{def:averages}
    \overline{Z}_k=\frac1N\sum_{i=1}^N\widehat C^i_k\qquad\text{for every }k\in\N,
\end{equation}
i.e., the tracking variables preserve the average of the local consensus points.
Given any $n\in\N$ and any vector $y\in\R^n$ we define the \emph{diameter} of $y$ as
\begin{equation}\label{eq:diam}
    \diam(y) = \max_{i=1:n}y_i - \min_{i=1:n}y_i = \max_{i,j=1:n} |y_i-y_j|,
\end{equation}
and the following two matrices
\begin{equation}\label{eq:defA}
     P =\left(\begin{array}{cc}I-J&0\\ 0&I-J\end{array}\right)\in\R^{2N\times 2N},\qquad A_\lambda=\left(\begin{array}{cc}W&\lambda I\\ W-I&I\end{array}\right)\in\R^{2N\times 2N}.
\end{equation}
\longv{The matrix $(I-J)$ is the orthogonal projection onto the subspace orthogonal to the \emph{consensus subspace} $\mathrm{span}\{\mathbf{1}_N\}.$ Analogously, $P$ is the orthogonal projection onto the subspace orthogonal to $\mathrm{span}\{(\mathbf{1}_N;0),\ (0;\mathbf{1}_N) \}.$
The matrix $A_\lambda$ governs the evolution of the disagreement among the tracking variables, $(I-J)\mathcal{Z}^s_k$, and of the tracking errors $\mathcal{H}^s_k$ (see Lemma \ref{lemma:consensus_structure}).}

The following lemma shows that the powers of $A_\lambda$ decay geometrically on the range of $P$.

\begin{lemma}\label{lemma:A}
Let $W$ satisfy Assumption \ref{a:W}, $P$ and $A_\lambda$ be given in \eqref{eq:defA}  and $\lambda\in(0,1)$. There exist $\kappa \in(0,1)$ and $\Lambda_A>0$,
depending on $\lambda$ and $W$, such that
\begin{equation}\label{eq:Apower}
    \|A_\lambda^k P v\|_2\leq \Lambda_A\kappa ^{\,k}\| P v\|_2\qquad\text{for every } v\in\R^{2N},\ k\in\N.
\end{equation}
\end{lemma}

\begin{proof}
Since $WJ=JW=J$, we have $A_\lambda P = P A_\lambda$, and since $ P ^2= P $,
we get $A_\lambda^k P v=(A_\lambda P )^k P v$ for every $k\in\N$ and $v\in\R^{2N}$.
Let us denote by $\rho(A_\lambda P )$ the spectral radius of $A_\lambda P $.
We prove that  $\rho(A_\lambda P )<1$.
Let us denote by $Q = (q_1,\dots,q_N)\in\R^{N\times N}$ an orthonormal matrix of eigenvectors of $W$ with $q_1 = N^{-1/2}\mathbf{1}_N$,  and let us define the block-diagonal matrix $\widehat Q = \diag(Q,Q).$
Then, we have 
\begin{equation*}
Q^\top WQ = \diag(\gamma_1,\dots,\gamma_N) =: \Gamma, \quad Q^\top(I-J)Q = I - Q^\top q_1q_1^\top Q = \diag(0,1,\dots,1) =: \Sigma,
\end{equation*}
therefore
\begin{equation}\begin{aligned}
\widehat Q^\top A_\lambda P  \widehat{Q} &= \left(\begin{matrix}
Q^\top W(I-J)Q & \lambda Q^\top(I-J)Q\\
Q^\top(W-I)(I-J)Q & Q^\top(I-J)Q  
\end{matrix}\right) =
\left(\begin{matrix}
\Gamma\Sigma & \lambda \Sigma\\
(\Gamma-I)\Sigma & \Sigma  
\end{matrix}\right).
\end{aligned}\end{equation}
In particular, there exists a permutation matrix $\Pi$ such that $\Pi\widehat Q^\top A_\lambda P  \widehat{Q}\Pi^\top = \diag(T_1,\dots, T_N)$ with 
$$T_1 = \left(\begin{matrix}
    0& 0\\
    0 & 0  
\end{matrix}\right), \quad T_j = \left(\begin{matrix}\gamma_j&\lambda\\ \gamma_j-1&1 \end{matrix}\right)\ \text{for }j=2,\dots,N.$$
The eigenvalues of $A_\lambda P $ are given by the eigenvalues of the blocks $T_1,\dots,T_N$. 
For every $j=2,\dots,N$, the eigenvalues of $T_j$ are the solutions of
\begin{equation*}
    t^2-(1+\gamma_j)t+\gamma_j+\lambda(1-\gamma_j) = 0.
\end{equation*}
Solving the equation one can see that the spectral radius of $T_j$ is given by
\begin{equation}
\rho(T_j) = \begin{cases}
    \left( \gamma_j+\lambda(1-\gamma_j)\right)^{1/2} & \text{if } (1-\gamma_j)(1-4\lambda-\gamma_j)<0,\\
    \frac12\left( 1+\gamma_j+((1-\gamma_j)(1-4\lambda-\gamma_j))^{1/2}\right) & \text{otherwise}.
\end{cases}     
\end{equation}
Since $-1\leq\gamma_j<1$ and $\lambda\in(0,1)$, we have that $\rho(T_j)<1$
and therefore $\rho(A_\lambda P )<1$. Furthermore, due to Gelfand's formula \cite[Corollary 5.6.14]{matan}
\begin{equation}\label{power_spectral}
\rho(A_\lambda P ) = \lim_{k\to+\infty}\|(A_\lambda P )^k\|_2^{1/k},
\end{equation}
for any choice of $\kappa \in(\rho(A_\lambda P ),1)$ there exists $k_0\in\N$ such that $\|(A_\lambda P )^k\|_2\leq\kappa ^{\,k}$ for every $k\geq k_0$. Defining
\begin{equation*}
    \Lambda_A=\max\Big\{1,\ \max_{0\leq k<k_0}\kappa ^{\,-k}\|(A_\lambda P )^k\|_2\Big\},
\end{equation*}
we get $\|(A_\lambda P )^k\|_2\leq \Lambda_A\kappa ^{\,k}$ for every $k\in\N$, and therefore
\begin{equation*}
    \|A_\lambda^k P v\|_2=\|(A_\lambda P )^k P v\|_2
    \leq \Lambda_A\kappa ^{\,k}\| P v\|_2,
\end{equation*}
which is the thesis.
\end{proof}
The following lemma collects some properties of the dynamics \eqref{eq:CBOpar} that we use in the analysis\longv{: the evolution of the diameter of each local swarm, the linear recursion satisfied by the tracking variables, and a bound on the diameter of the whole swarm in terms of these quantities.}\shortv{.}

\begin{lemma}\label{lemma:consensus_structure}
Let Assumption \ref{a:W} hold. For every $s=1,\dots,d$ and $k\in\N$ the following hold.
\begin{enumerate}[label=(\roman*)]
    \item For every $i=1,\dots,N$,
    \begin{equation}\label{eq:local_diam}
        \diam(\mathcal{X}^{s,i}_{k+1})\leq\Big(1-\lambda+2\sigma\max_{\ell\in\calM_i}|\xi^\ell_{k,s}|\Big)\diam(\mathcal{X}^{s,i}_k).
    \end{equation}
    \item $\mathcal{H}^s_0=0$, $J\mathcal{H}^s_k = 0$, and
    \begin{equation}\label{eq:tracking_rec}
        \left(\begin{array}{c}(I-J)\mathcal{Z}^s_{k+1}\\ \mathcal{H}^s_{k+1}\end{array}\right)
        =A_\lambda\left(\begin{array}{c}(I-J)\mathcal{Z}^s_{k}\\ \mathcal{H}^s_{k}\end{array}\right)
        +\left(\begin{array}{c}(I-J)\zeta^s_k\\ 0\end{array}\right),
    \end{equation}
    with $A_\lambda$ given in \eqref{eq:defA} and 
    \begin{equation}\label{eq:zeta}
\zeta^s_k=\widehat{\mathcal{C}}^s_{k+1}-(1-\lambda)\widehat{\mathcal{C}}^s_k-\lambda \mathcal{Z}^s_k\in\R^N.       
    \end{equation} Moreover, for every $i=1,\dots,N$
    \begin{equation}\label{eq:zeta_bound}
        |\zeta_{k,i}^s|\leq\Big(1-\lambda+\sigma\max_{\ell\in\calM_i}|\xi^\ell_{k,s}|\Big)\diam(\mathcal{X}^{s,i}_k).
    \end{equation}
    
    \item The following inequality holds:
    \begin{equation}\label{eq:global_diam}
        \diam(\mathcal{X}^s_k)\leq2\max_{i=1:N}\diam(\mathcal{X}^{s,i}_k)+2\|\mathcal{H}^s_k\|_2+2\|(I-J)\mathcal{Z}^s_k\|_2.
    \end{equation}
\end{enumerate}
\end{lemma}

\begin{proof}
For every $\ell,m\in\calM_i$, by \eqref{eq:CBOpar}, we have
\begin{equation*}
    x^\ell_{k+1,s}-x^m_{k+1,s}=(1-\lambda)(x^\ell_{k,s}-x^m_{k,s})
    +\sigma\left((x^\ell_{k,s}-\widehat\cons^i_{k,s})\xi^\ell_{k,s}-(x^m_{k,s}-\widehat\cons^i_{k,s})\xi^m_{k,s}\right).
\end{equation*}
Since $\widehat\cons^i_{k,s}$ is a convex combination of the entries of $\mathcal{X}^{s,i}_k$, we have
$|x^\ell_{k,s}-\widehat\cons^i_{k,s}|\leq\diam(\mathcal{X}^{s,i}_k)$ for every $\ell\in\calM_i$. Therefore, taking the absolute value and the maximum over $\ell,m\in\calM_i$ we get \eqref{eq:local_diam}.

Let us now prove $(ii)$.  %
Since $Z^i_0=\widehat C^i_0$ for every $i$, we have $\mathcal{H}^s_0=0$. Moreover, \longv{\eqref{eq:CBOpar} yields $\mathcal{Z}^s_{k+1}=\widehat{\mathcal{C}}^s_{k+1}-\widehat{\mathcal{C}}^s_{k}+W\mathcal{Z}^s_k$, and therefore}
\begin{equation}\label{auxh}
    \mathcal{H}^s_{k+1}=\mathcal{Z}^s_{k+1}-\widehat{\mathcal{C}}^s_{k+1}=W\mathcal{Z}^s_k-\widehat{\mathcal{C}}^s_k=\mathcal{H}^s_k+(W-I)\mathcal{Z}^s_k.
\end{equation}
Multiplying both sides by $J$ and using the fact that $J(W-I)=0$, we have that $J\mathcal{H}^s_{k+1} = J\mathcal{H}^s_k$ for every $k$ which, together with the fact that
$\mathcal{H}^s_0=0$ ensures $J\mathcal{H}^s_k = 0$ for every $k.$
Since $(W-I)\mathbf{1}_N=0$, equation \eqref{auxh} also implies
\begin{equation}\label{auxh2}
    \mathcal{H}^s_{k+1}=\mathcal{H}^s_k+(W-I)(I-J)\mathcal{Z}^s_k.
\end{equation}
Equality \eqref{eq:tracking_rec} then follows directly from the equation above, \eqref{auxh} and the definition of $\mathcal{H}^s_k,\ \zeta^s_{k}$ and $A_\lambda.$

Using \eqref{eq:CBOpar} and the fact that $\widehat\cons^i_{k+1,s}=\sum_{\ell\in\calM_i}\widehat v^{i,\ell}_{k+1}x^\ell_{k+1,s}$
with $\widehat v^{i,\ell}_{k+1}>0$ and $\sum_{\ell\in\calM_i}\widehat v^{i,\ell}_{k+1}=1$, we have
\begin{equation*}
    \zeta^s_{k,i}=\widehat\cons^i_{k+1,s}-(1-\lambda)\widehat\cons^i_{k,s}-\lambda z^i_{k,s}
    =(1-\lambda)\sum_{\ell\in\calM_i}\widehat v^{i,\ell}_{k+1}(x^\ell_{k,s}-\widehat\cons^i_{k,s})
    +\sigma\sum_{\ell\in\calM_i}\widehat v^{i,\ell}_{k+1}(x^\ell_{k,s}-\widehat\cons^i_{k,s})\xi^\ell_{k,s}.
\end{equation*}
Taking the absolute value, using again the fact that $|x^\ell_{k,s}-\widehat\cons^i_{k,s}|\leq\diam(\mathcal{X}^{s,i}_k)$ and that  $\sum_{\ell\in\calM_i}\widehat v^{i,\ell}_{k+1}=1$, we get \eqref{eq:zeta_bound}, which concludes the proof of $(ii).$

We now prove $(iii)$. %
For any $\ell\in\calM_i$ and $m\in\calM_j$, using the fact that $(J\mathcal{Z}^s_k)_i=(J\mathcal{Z}^s_k)_j$, we have 
\begin{equation*}\begin{aligned}
    |x^\ell_{k,s}-x^m_{k,s}|&\leq|x^\ell_{k,s}-\widehat\cons^i_{k,s}|+|\widehat\cons^i_{k,s}-z^i_{k,s}|
    +|z^j_{k,s}-\widehat\cons^j_{k,s}|+|\widehat\cons^j_{k,s}-x^m_{k,s}|+|z^i_{k,s}-z^j_{k,s}|\\
    &\leq \diam(\mathcal{X}_k^{s,i})+\|\mathcal{H}^s_k\|_{\infty}+\diam(\mathcal{X}_k^{s,j})+\|\mathcal{H}^s_k\|_{\infty}+|((I-J)\mathcal{Z}^s_k)_i-((I-J)\mathcal{Z}^s_k)_j|\\
    &\leq \diam(\mathcal{X}_k^{s,i})+\diam(\mathcal{X}_k^{s,j})+2\|\mathcal{H}^s_k\|_{\infty}+2\|(I-J)\mathcal{Z}^s_k\|_{\infty}.
\end{aligned}\end{equation*}
Taking the maximum over all $\ell, m$ and using the fact that $\|\cdot\|_{\infty}\leq\|\cdot\|_2$, we complete the proof.
\end{proof}
Since the random variables $\xi^\ell_{k,s}$
are independent and normally distributed, for every $s=1,\dots,d$, $i=1,\dots,N$ and $k\in\N$ we have \cite{gaussianbound}:
\begin{equation}\label{eq:gauss_max}
    \mathbb{E}\left[\max_{\ell\in\calM_i}|\xi^\ell_{k,s}|\right]\leq 2\big(\log(\sqrt2M)\big)^{1/2},\qquad
    \mathbb{E}\left[\max_{\ell\in\calM_i}|\xi^\ell_{k,s}|^2\right]\leq 4\log(\sqrt2M).
\end{equation}
We define the following constant, which depends on the dynamics parameters and on the number of per-swarm agents $M$
\begin{equation}\label{eq:vartheta}
    \vartheta = 1-\lambda+4\sigma\big(\log(\sqrt2M)\big)^{1/2}.
\end{equation}
The constant $\vartheta$ bounds the contraction factor of the diameter of each local swarm over one iteration (see Lemma \ref{lemma:local_moments} below). \longv{The term $1-\lambda$ accounts for the drift towards the tracking variable, while the second term accounts for the noise, through the expected maximum of $M$ independent Gaussian variables in \eqref{eq:gauss_max}. The condition $\vartheta<1$ will be our main assumption for consensus formation.}
Moreover, for every $s=1,\dots,d$, we define the following quantities, which depend on the initial distribution of the particles
\begin{equation}\label{eq:D0E0}
    D_{\mathcal{X},s} = \sum_{i=1}^N\mathbb{E}\left[\diam(\mathcal{X}^{s,i}_0)^2\right],\qquad
    D_{\mathcal{C},s} = \mathbb{E}\left[\|(I-J)\widehat{\mathcal{C}}^s_0\|_2^2\right],\qquad
    D_{\mathcal{X}} = \sum_{s=1}^d D_{\mathcal{X},s}.
\end{equation}
\longv{The quantities $D_{\mathcal{X},s}$ and $D_{\mathcal{C},s}$ measure the spread of the particles within each computational unit and the disagreement among the consensus points of different units, respectively, at $k=0.$}

The following lemma shows that the local diameters and the  terms $\zeta^s_k$ in \eqref{eq:tracking_rec} decay geometrically in mean square, with rate $\vartheta$.
\begin{lemma}\label{lemma:local_moments}
For every $s=1,\dots,d$, $i=1,\dots,N$ and $k\in\N$ we have
\begin{enumerate}[label=(\roman*)]
    \item $\mathbb{E}\left[\diam(\mathcal{X}^{s,i}_k)^2\right]\leq\vartheta^{2k}\,\mathbb{E}\left[\diam(\mathcal{X}^{s,i}_0)^2\right]$;
    \item $\mathbb{E}\left[\|\zeta^s_k\|_2^2\right]\leq\vartheta^{2(k+1)}D_{\mathcal{X},s}$.
\end{enumerate}
\end{lemma}
\begin{proof}
Applying \eqref{eq:local_diam} recursively and taking the expectation\longv{, using the fact that the random variables $\xi^\ell_{t,s}$ are independent,} we get
\begin{equation}\label{eq:local_prod}\begin{aligned}
    \E[\diam(\mathcal{X}^{s,i}_k)^2]&\leq \E\left[\diam(\mathcal{X}^{s,i}_0)^2\prod_{t=0}^{k-1}(1-\lambda+2\sigma\max_{\ell\in\calM_i}|\xi^\ell_{t,s}|)^2\right] \\
    &= \E[\diam(\mathcal{X}^{s,i}_0)^2]\prod_{t=0}^{k-1}\E[(1-\lambda+2\sigma\max_{\ell\in\calM_i}|\xi^\ell_{t,s}|)^2].
\end{aligned}\end{equation}
By \eqref{eq:gauss_max}, we have
\begin{equation*}
    \mathbb{E}\left[(1-\lambda+2\sigma\max_{\ell\in\calM_i}|\xi^\ell_{t,s}|)^2\right]
    \leq(1-\lambda)^2+8\sigma(1-\lambda)\big(\log(\sqrt2M)\big)^{1/2}+16\sigma^2\log(\sqrt2M)=\vartheta^2 .
\end{equation*}
This, together with \eqref{eq:local_prod}, gives $(i).$
By \eqref{eq:zeta_bound} we have 
$$|\zeta^s_{k,i}|\leq(1-\lambda+\sigma\max_{\ell\in\calM_i}|\xi^\ell_{k,s}|)\diam(\mathcal{X}^{s,i}_k).$$ Proceeding as above we get $$\mathbb{E}\left[|\zeta^s_{k,i}|^2\right]\leq\vartheta^{2(k+1)}\,\mathbb{E}\left[\diam(\mathcal{X}^{s,i}_0)^2\right],$$ and summing over $i=1,\dots,N$ we conclude the proof.
\end{proof}
The following lemma shows that the tracking errors and the diameter of the whole swarm decay geometrically, at a rate $\varrho$ given by the slower of the contraction of the local swarms ($\vartheta$, Lemma \ref{lemma:local_moments}) and the contraction of the communication network ($\kappa$, Lemma \ref{lemma:A}).
\begin{lemma}\label{lemma:global_consensus}
Let Assumption \ref{a:W} hold, let $\kappa $ and $\Lambda_A$ be given by Lemma \ref{lemma:A}, with
$\kappa \neq\vartheta$, and define $\varrho=\max\{\kappa ,\vartheta\}$. Then for every $s=1,\dots,d$ and $k\in\N$ we have
\begin{enumerate}[label=(\roman*)]
    \item $\mathbb{E}\left[\|(I-J)\mathcal{Z}^s_k\|_2^2+\|\mathcal{H}^s_k\|_2^2\right]\leq\Lambda_{\mathcal{Z},s}^2\,\varrho^{2k}$, with
    $\Lambda_{\mathcal{Z},s}=\Lambda_A\Big(D_{\mathcal{C},s}^{1/2}+\dfrac{\vartheta\,D_{\mathcal{X},s}^{1/2}}{|\kappa -\vartheta|}\Big)$;
    \item $\mathbb{E}\left[\diam(\mathcal{X}^s_k)^2\right]\leq\Lambda_{\mathcal{X},s}^2\,\varrho^{2k}$, with $\Lambda_{\mathcal{X},s}=2D_{\mathcal{X},s}^{1/2}+2\sqrt2\,\Lambda_{\mathcal{Z},s}$;
    \item for every $\bar\varrho>\varrho$ there exists a $\mathbb{P}$-almost surely finite random variable $\Theta_s$, such that for every $k\in\N$
    \begin{equation*}
        \diam(\mathcal{X}^s_k)+\|(I-J)\mathcal{Z}^s_k\|_2+\|\mathcal{H}^s_k\|_2\leq \Theta_s\,\bar\varrho^{\,k}.
    \end{equation*}
\end{enumerate}
\end{lemma}

\begin{proof}
Recursively applying 
\eqref{eq:tracking_rec}, using the definition of $P$ given in \eqref{eq:defA}, the fact that $(I-J)^2=I-J$ and $\mathcal{H}^s_0 = 0$, we get
\begin{equation*}\begin{aligned}
    \left(\begin{array}{c}(I-J)\mathcal{Z}^s_{k}\\ \mathcal{H}^s_{k}\end{array}\right)&=A_\lambda^k\left(\begin{array}{c}(I-J)\mathcal{Z}^s_{0}\\ 0 \end{array}\right)+\sum_{t=0}^{k-1}A_\lambda^{k-1-t}\left(\begin{array}{c}(I-J)\zeta^s_t\\ 0\end{array}\right) \\&= A_\lambda^k P \left(\begin{array}{c}(I-J)\mathcal{Z}^s_{0}\\ 0\end{array}\right)+\sum_{t=0}^{k-1}A_\lambda^{k-1-t} P \left(\begin{array}{c}\zeta^s_t\\ 0\end{array}\right),
\end{aligned}\end{equation*}
and together with Lemma \ref{lemma:A} we have
\begin{equation*}
    \left\|\left(\begin{array}{c}(I-J)\mathcal{Z}^s_{k}\\ \mathcal{H}^s_{k}\end{array}\right)\right\|_2\leq \Lambda_A\kappa ^{\,k}\|(I-J)\widehat{\mathcal{C}}^s_0\|_2+\Lambda_A\sum_{t=0}^{k-1}\kappa ^{\,k-1-t}\|\zeta^s_t\|_2.
\end{equation*}
Using the fact that for any random vectors $x,y$ it holds $\E[\|x+y\|^2]^{1/2}\leq \E[\|x\|^2]^{1/2}+\E[\|y\|^2]^{1/2}$, the definition of the constants $D_{\mathcal{C},s}$ and $D_{\mathcal{X},s}$ given in \eqref{eq:D0E0}
and Lemma \ref{lemma:local_moments} we then have
\begin{equation*}
    \mathbb{E}\left[\left\|\left(\begin{array}{c}(I-J)\mathcal{Z}^s_{k}\\ \mathcal{H}^s_{k}\end{array}\right)\right\|_2^2\right]^{1/2}
    \leq \Lambda_A\kappa ^{\,k}D_{\mathcal{C},s}^{1/2}+\Lambda_AD_{\mathcal{X},s}^{1/2}\sum_{t=0}^{k-1}\kappa ^{\,k-1-t}\vartheta^{\,t+1}
    =\Lambda_A\kappa ^{\,k}D_{\mathcal{C},s}^{1/2}+\Lambda_AD_{\mathcal{X},s}^{1/2}\,\vartheta\,\frac{\kappa ^{\,k}-\vartheta^{k}}{\kappa -\vartheta}.
\end{equation*}
Part $(i)$ of the statement follows by squaring both sides of the inequality above, by the definition of $\varrho$ and $\Lambda_{\mathcal{Z},s}$.

We now prove $(ii)$. By \eqref{eq:global_diam}, using $\max_i\diam( \mathcal{X}^{s,i}_k)\leq\big(\sum_i\diam( \mathcal{X}^{s,i}_k)^2\big)^{1/2}$ and
$\|\mathcal{H}^s_k\|_2+\|(I-J)\mathcal{Z}^s_k\|_2\leq\sqrt2\left\|\left((I-J)\mathcal{Z}^s_{k}; \mathcal{H}^s_{k}\right)\right\|_2$, we have
\begin{equation*}
\diam( \mathcal{X}^s_k)\leq2\Big(\sum_{i=1}^N\diam( \mathcal{X}^{s,i}_k)^2\Big)^{1/2}+2\sqrt2\left\|\left(\begin{array}{c}(I-J)\mathcal{Z}^s_{k}\\ \mathcal{H}^s_{k}\end{array}\right)\right\|_2 .
\end{equation*}
Taking the expected value, using Lemma \ref{lemma:local_moments} and $(i)$, we get $(ii)$.

By the first and second statement we have that for every $k\in\N$
\begin{equation*}
    \mathbb{E}[(\diam( \mathcal{X}^s_k)+\|(I-J)\mathcal{Z}^s_k\|_2+\|\mathcal{H}^s_k\|_2)^2]
    \leq3(\Lambda_{\mathcal{X},s}^2+\Lambda_{\mathcal{Z},s}^2)\varrho^{2k}.
\end{equation*}
 Therefore
\begin{equation*}
    \mathbb{E}\left[\sum_{k=0}^{\infty}\bar\varrho^{\,-2k}(\diam( \mathcal{X}^s_k)+\|(I-J)\mathcal{Z}^s_k\|_2+\|\mathcal{H}^s_k\|_2)^2\right]\leq 3(\Lambda_{\mathcal{X},s}^2+\Lambda_{\mathcal{Z},s}^2)\sum_{k=0}^\infty\Big(\frac{\varrho}{\bar\varrho}\Big)^{2k}<+\infty,
\end{equation*}
which implies that $\sum_{k}\bar\varrho^{\,-2k}(\diam( \mathcal{X}^s_k)+\|(I-J)\mathcal{Z}^s_k\|_2+\|\mathcal{H}^s_k\|_2)^2<+\infty$ almost surely. In particular the terms of the series
go to zero almost surely, and part $(iii)$ follows.
\end{proof}

We can now state our first main result: if $\vartheta<1$, the agents asymptotically achieve consensus almost surely.
\begin{theorem}\label{thm:consensus}
Let Assumption \ref{a:W} hold and let $\lambda\in(0,1)$ and $\sigma>0$ be such that $\vartheta<1$, with $\vartheta$ given in \eqref{eq:vartheta}. Then
\begin{equation*}
    \lim_{k\to+\infty}\max_{\ell,p=1:NM}\|X^\ell_k-X^p_k\|_\infty=0\quad\text{almost surely}.
\end{equation*}
Moreover, there exists $X_\infty\in\R^d$ such that, almost surely,
\begin{equation*}
    \lim_{k\to+\infty}X^\ell_k=X_\infty\ \text{ for every } \ell=1,\dots,NM,\qquad
    \lim_{k\to+\infty}Z^i_k=X_\infty\ \text{ for every } i=1,\dots,N.
\end{equation*}
\end{theorem}

\begin{proof}
By definition of $\diam(\cdot)$ and part $(iii)$ of Lemma \ref{lemma:global_consensus}, we have that for every $\bar\varrho>\varrho = \max\{\kappa, \vartheta\}$ the following inequality holds almost surely
\begin{equation*}
    \max_{\ell,p=1:NM}\|X^\ell_k-X^p_k\|_\infty=\max_{s=1:d}\max_{\ell,p=1:NM}|x^\ell_{k,s}-x^p_{k,s}|=\max_{s=1:d}\diam( \mathcal{X}^s_k)\leq \bar\varrho^k \max_{s=1:\dims}\Theta_s.
\end{equation*}
Since $\kappa,\vartheta<1$, in particular we can take $\bar\varrho\in(\varrho,1)$
and the first part of the statement follows by taking the limit as $k$ goes to infinity.

Let us now fix $s\in\{1,\dots,d\}$. 
 From the definition of 
 $\overline{Z}_{k}$ and  \eqref{def:averages} it follows 
\begin{equation*}
    \overline{z}_{k+1,s}-\overline{z}_{k,s}=\frac1N\mathbf{1}_N^\top\big(\widehat{\mathcal{C}}^s_{k+1}-\widehat{\mathcal{C}}^s_k\big) = \frac1N\mathbf{1}_N^\top\big( \lambda \mathcal{H}^s_k+\zeta^s_k\big),
\end{equation*}
where $\zeta^s_k$ is defined in \eqref{eq:zeta}.
The Cauchy--Schwarz inequality and  
and $\mathbf{1}_N^\top\mathcal{H}^s_k=0$ yield
\begin{equation*}
    |\overline{z}_{k+1,s}-\overline{z}_{k,s}|\leq\frac1N\|\mathbf{1}_N\|_2\,\|\zeta^s_k\|_2=\frac{1}{\sqrt N}\|\zeta^s_k\|_2 .
\end{equation*}
Taking the expected value, by Jensen's inequality and Lemma \ref{lemma:local_moments},
we have
\begin{equation*}\begin{aligned}
    \mathbb{E}\left[\sum_{k=0}^{\infty}|\overline{z}_{k+1,s}-\overline{z}_{k,s}|\right]
    &\leq\frac{1}{\sqrt N}\sum_{k=0}^{\infty}\mathbb{E}\left[\|\zeta^s_k\|_2\right]
    \leq\frac{1}{\sqrt N}\sum_{k=0}^{\infty}\mathbb{E}\left[\|\zeta^s_k\|_2^2\right]^{1/2}\\
    &\leq\frac{D_{\mathcal{X},s}^{1/2}}{\sqrt N}\sum_{k=0}^{\infty}\vartheta^{k+1}
    =\frac{\vartheta\,D_{\mathcal{X},s}^{1/2}}{\sqrt N(1-\vartheta)}<+\infty.
\end{aligned}\end{equation*}
Therefore $\sum_{k=0}^{\infty}|\overline{z}_{k+1,s}-\overline{z}_{k,s}|<+\infty$ almost surely which implies that  
$\sum_{t=0}^{\infty}(\overline{z}_{t+1,s}-\overline{z}_{t,s})\in\R$ almost surely. This, together with the telescopic relation
\begin{equation*}
    \overline{z}_{k,s}=\overline{z}_{0,s}+\sum_{t=0}^{k-1}(\overline{z}_{t+1,s}-\overline{z}_{t,s})\qquad\text{for every }k\in\N,
\end{equation*}
implies that $(\overline{z}_{k,s})_k$ converges almost surely to some $x_{\infty,s}\in\R$.

We now have to prove that $z^i_{k,s}$ and $x^\ell_{k,s}$ converge to $x_{\infty,s}.$
Let us consider $i\in\{1,\dots,N\}$ and $\ell\in\calM_i$. Using the fact that $z^i_{k,s}-\overline{z}_{k,s}=((I-J)\mathcal{Z}^s_k)_i$ and
$\widehat\cons^i_{k,s}-z^i_{k,s}=-(\mathcal{H}^s_k)_i$,
we have
\begin{align*}
    |z^i_{k,s}-\overline{z}_{k,s}|&\leq\|(I-J)\mathcal{Z}^s_k\|_2,\\
    |x^\ell_{k,s}-\overline{z}_{k,s}|&\leq|x^\ell_{k,s}-\widehat\cons^i_{k,s}|+|\widehat\cons^i_{k,s}-z^i_{k,s}|+|z^i_{k,s}-\overline{z}_{k,s}|
    \leq\diam( \mathcal{X}^s_k)+\|\mathcal{H}^s_k\|_2+\|(I-J)\mathcal{Z}^s_k\|_2.
\end{align*}
 By part $(iii)$ of Lemma \ref{lemma:global_consensus} the right-hand side
converges to zero almost surely. Therefore $z^i_{k,s}\to x_{\infty,s}$ and $x^\ell_{k,s}\to x_{\infty,s}$ almost surely, for every
$i$, every $\ell$ and every $s$, and the thesis follows with $X_\infty=(x_{\infty,1},\dots,x_{\infty,d})$.
\end{proof}

\subsection{Convergence in expectation}\label{sec:convergence}
Let us define 
$$V_k = \frac{1}{NM}\sum_{i=1}^N\sum_{\ell\in\mathcal{M}_i}\|X^\ell_k-x^*\|^2,$$
which measures the mean squared distance of the particles from the global minimizer $x^*$.
We will need the following recursive inequality for $\E[V_k].$
\begin{lemma}\label{lemma:Vk}
For every $k\in\N$, the following inequality holds, with $ D_{\mathcal{X}} $ defined in \eqref{eq:D0E0}
\begin{equation}\label{eq:Vrec}
    \E[V_{k+1}]\leq(1-\lambda)\E[V_k] + \frac\lambda N\sum_{i=1}^N\E\big[\|Z^i_k-x^*\|^2\big] + \frac{\sigma^2}{N}D_{\mathcal{X}}\vartheta^{2k}.
\end{equation}
\end{lemma}

\begin{proof}
In the following, given any $\ell\in\{1,\dots,NM\}$ we denote by $\ind{\ell}$ the index in $\{1,\dots,N\}$ such that $\ell\in\calM_{\ind{\ell}}$.
For every $k\in\N$ and every $\ell=1,\dots,NM$, let us define
\begin{equation*}
    A^\ell_k = (1-\lambda)(X^\ell_k-x^*)+\lambda(Z^{\ind{\ell}}_k-x^*),\qquad
    B^\ell_k = \sigma\sum_{s=1}^\dims(x^\ell_{k,s}-\widehat\cons^{\ind{\ell}}_{k,s})\xi^\ell_{k,s}\mathbf{e}_s.
\end{equation*}
By \eqref{eq:CBOpar} we have $X^\ell_{k+1}-x^*=A^\ell_k+B^\ell_k$, and therefore
\begin{equation}\label{aux_xnext}
    \E\big[\|X^\ell_{k+1}-x^*\|^2\big] = \E\big[\|A^\ell_k\|^2\big]+\E\big[\|B^\ell_k\|^2\big]+2\E\big[(A^\ell_k)^\top B^\ell_k\big].
\end{equation}
Let us now bound the three terms on the right-hand side. Since $\|\cdot\|^2$ is convex,
\begin{equation}\label{aux_A}
    \E\big[\|A^\ell_k\|^2\big] \leq (1-\lambda)\E\big[\|X^\ell_k-x^*\|^2\big]+\lambda\E\big[\|Z^{\ind{\ell}}_k-x^*\|^2\big].
\end{equation}
The random variables $(\xi^\ell_{k,s})$ are independent of $X^\ell_k$, $Z^{\ind{\ell}}_k$ and $\widehat C^{\ind{\ell}}_k$,
and $\E[(\xi^\ell_{k,s})^2]=1$. Moreover, since $\widehat\cons^{\ind{\ell}}_{k,s}$ is a convex combination of the entries of $ \mathcal{X}^{s,\ind{\ell}}_k$, we have
$|x^\ell_{k,s}-\widehat\cons^{\ind{\ell}}_{k,s}|\leq\diam( \mathcal{X}^{s,\ind{\ell}}_k)$. Therefore, using Lemma \ref{lemma:local_moments},
\begin{equation}\label{aux_B}\begin{aligned}
    \E\big[\|B^\ell_k\|^2\big] & =  \sigma^2\sum_{s=1}^\dims\E\big[(x^\ell_{k,s}-\widehat\cons^{\ind{\ell}}_{k,s})^2(\xi^\ell_{k,s})^2\big]
    = \sigma^2\sum_{s=1}^\dims\E\big[(x^\ell_{k,s}-\widehat\cons^{\ind{\ell}}_{k,s})^2\big]\\
    &\leq \sigma^2 \sum_{s=1}^\dims \E\big[\diam( \mathcal{X}^{s,\ind{\ell}}_k)^2\big]\leq \sigma^2 \vartheta^{2k}\sum_{s=1}^\dims \E\big[\diam( \mathcal{X}^{s,\ind{\ell}}_0)^2\big].
\end{aligned}\end{equation}
Finally, since $\xi^\ell_{k,s}$ is independent of $X^\ell_k$, $Z^{\ind{\ell}}_k$ and $\widehat C^{\ind{\ell}}_k$, and  $\E[\xi^\ell_{k,s}]=0$, we have $\E\big[(A^\ell_k)^\top B^\ell_k\big]=0$.
This, together with \eqref{aux_xnext}, \eqref{aux_A}, \eqref{aux_B} and the fact that $|\calM_i|=M$ for every $i$, yields
\begin{equation*}\begin{aligned}
    \E[V_{k+1}] &\leq \frac{1}{NM}\sum_{\ell=1}^{NM}\left((1-\lambda)\E\big[\|X^\ell_k-x^*\|^2\big]+\lambda\E\big[\|Z^{\ind{\ell}}_k-x^*\|^2\big]
    +\sigma^2 \vartheta^{2k}\sum_{s=1}^\dims \E\big[\diam( \mathcal{X}^{s,\ind{\ell}}_0)^2\big]\right)\\
    & = (1-\lambda)\E[V_k] + \frac\lambda N \sum_{i=1}^N\E\big[\|Z^{i}_k-x^*\|^2\big]
    + \frac{\sigma^2}{N}\vartheta^{2k}\sum_{s=1}^\dims\sum_{i=1}^N \E\big[\diam( \mathcal{X}^{s,i}_0)^2\big],
\end{aligned}\end{equation*}
and the thesis follows from the definition of $D_{\mathcal{X},s}$ and $D_{\mathcal{X}}$.
\end{proof}

Let us denote by $C^i_k$ the consensus point of the swarm at node $i$ computed using the true function $f$, and with $C_k$ the consensus point of the whole swarm, again computed with $f$, namely
\begin{equation}\label{eq:globalcons}
    C^i_k = \frac{\sum_{\ell\in\M_i} X^\ell_k \exp(-\alpha  f(X^\ell_k))}{\sum_{\ell\in\M_i}\exp(-\alpha f(X^\ell_k))},\quad C_k = \frac{\sum_{\ell=1}^{NM}X^\ell_k\exp(-\alpha f(X^\ell_k))}{\sum_{\ell=1}^{NM}\exp(-\alpha f(X^\ell_k))}.
\end{equation}
Notice that $C^i_k$ and $C_k$ are defined only for the purpose of this theoretical analysis and they are never computed by the algorithm.

Given $r>0$ and $\mathcal{B}_r$ defined in \eqref{def:B_r}, we define
\begin{equation}
    \begin{aligned}
        & \mathcal{I}^i_{k,r} = \{\ell\in\calM_i\ |\ X^\ell_k\in\mathcal B_r\},\quad  \mathcal{I}_{k,r} = \bigcup_{i=1}^N\mathcal{I}^i_{k,r}.
    \end{aligned}
\end{equation}

We have the following finite agents version of the quantitative Laplace principle \cite{fornasier2024}.
\begin{lemma}[Quantitative Laplace Principle]\label{lemma:laplaceprinciple}
Let Assumption \ref{a:f} hold and $C_k$ be given in   \eqref{eq:globalcons}. Consider $r\in(0,R_0]$ and $q>0$ such that $q+f_r\leq f_{\infty}$.
For every $k\in\N$ we have
\begin{enumerate}[label=(\roman*)]
    \item if  $|{\cal I}_{k,r}|>0$ then
$$\|C_k-x^*\|\leq \frac{1}{\beta}(q+f_r)^\nu+\frac{e^{-\alpha q}}{ |{\cal I}_{k,r}|}NM V_k^{1/2}$$ 
    \item if  $|{\cal I}_{k,r}|=0$ it still holds that
    $$\|C_k-x^*\|\leq \sqrt{NM}V_k^{1/2}.$$
\end{enumerate}

\end{lemma}

\begin{proof}
To prove $(i)$ we follow the proof of \cite[Proposition 4.5]{fornasier2024}. For $\ell=1,\dots,NM$ let $w_\ell=e^{-\alpha f(X^\ell_k)}$.
By definition of $C_k$ in \eqref{eq:globalcons} and the triangle inequality, for any $\tilde r\geq r$ we have
\begin{equation}\label{diffsol}\begin{aligned}
\|C_k-x^*\| & = \left\|\frac{1}{\sum_{\ell=1}^{NM} w_\ell}\sum_{\ell=1}^{NM} w_\ell(X^\ell_k-x^*)\right\|
\leq \frac{1}{\sum_{\ell=1}^{NM} w_\ell}\sum_{\ell=1}^{NM} w_\ell\|X^\ell_k-x^*\| \\
&=\frac{1}{\sum_{\ell=1}^{NM} w_\ell}\left( \sum_{\ell\in\mathcal{I}_{k,\tilde r}} w_\ell\|X^\ell_k-x^*\|
+ \sum_{\ell\notin\mathcal{I}_{k,\tilde r}} w_\ell\|X^\ell_k-x^*\|\right)\\
&\leq \tilde r + \frac{\max_{\ell\notin\mathcal{I}_{k,\tilde r}} w_\ell}{\sum_{\ell=1}^{NM} w_\ell}
\sum_{\ell\notin\mathcal{I}_{k,\tilde r}}\|X^\ell_k-x^*\| \\
&\leq \tilde r + \frac{\exp\big(-\alpha\inf_{x\notin\mathcal B_{\tilde r}}f(x)\big)}{\sum_{\ell=1}^{NM} w_{\ell}}
\sum_{\ell=1}^{NM}\|X^\ell_k-x^*\|.
\end{aligned}\end{equation}
We take $\tilde r=\frac{1}{\beta}(q+f_r)^{\nu}$. Notice that $\tilde r\geq r$, since by Assumption \ref{a:f}
\begin{equation*}
    \tilde r = \frac{1}{\beta}(q+f_r)^{\nu} \geq \frac{1}{\beta}f_r^\nu=\max_{x\in \mathcal B_r}\frac{1}{\beta}f(x)^\nu
    \geq \max_{x\in \mathcal B_r}\|x-x^*\| = r.
\end{equation*}
Moreover, by Assumption \ref{a:f} and the definition of $\tilde r$, for every $x\notin\mathcal B_{\tilde r}$ either
$x\notin\mathcal B_{R_0}$ and $f(x)\geq f_\infty$, or $x\in\mathcal B_{R_0}$ and $f(x)\geq(\beta\|x-x^*\|)^{1/\nu}\geq(\beta\tilde r)^{1/\nu}$.
Since $(\beta\tilde r)^{1/\nu}=q+f_r\leq f_\infty$, we then have
\begin{equation}\label{inff}
    \inf_{x\notin\mathcal B_{\tilde r}}f(x)\geq\min\left\{f_\infty,(\beta\tilde r)^{1/\nu}\right\}= q+f_r .
\end{equation}
Since $f(x)\leq f_r$ for every $x\in\mathcal B_r$ we have
\begin{equation*}
    \sum_{\ell=1}^{NM} w_\ell\geq \sum_{\ell\in\mathcal{I}_{k,r}} e^{-\alpha f(X^\ell_k)}\geq |\mathcal{I}_{k,r}|\, e^{-\alpha f_r}.
\end{equation*}
Using this inequality and \eqref{inff} in \eqref{diffsol} we get
\begin{equation*}
    \|C_k-x^*\|\leq\frac1\beta(q+f_r)^\nu+\frac{e^{-\alpha(q+f_r)}}{|\mathcal{I}_{k,r}|\,e^{-\alpha f_r}}\sum_{\ell=1}^{NM}\|X^\ell_k-x^*\|
    =\frac1\beta(q+f_r)^\nu+\frac{e^{-\alpha q}}{|\mathcal{I}_{k,r}|}\sum_{\ell=1}^{NM}\|X^\ell_k-x^*\|,
\end{equation*}
which, together with the definition of $V_k$ yields part $(i)$ of the thesis.

Let us now consider part $(ii).$ Since $C_k$ is a convex combination of the agents $(X^\ell_k)_{\ell=1}^{NM}$, we have
\begin{equation}
    \begin{aligned}
        \|C_k-x^*\| \leq \max_{\ell=1:NM}\|X^\ell_k-x^*\|\leq \left(\sum_{\ell=1}^{NM}\|X^\ell_k-x^*\|^2\right)^{1/2},
    \end{aligned}
\end{equation}
and the thesis follows from the definition of $V_k.$
\end{proof}

Let us introduce the following error vector $\mathcal E_k$.
\begin{equation} \label{def:error_vec}
\mathcal{E}_{k}:=\left(\left|f(X^\ell_k)-\widehat f(X^\ell_k, \omega^\ell_k)\right|\right)_{\ell=1:NM}\in\R^{NM}.
\end{equation}
The next lemma quantifies the effect of the noise in the function evaluations on the local consensus points.
\begin{lemma}\label{lemmadiffxalpha}
Let Assumption \ref{a:f_hat} hold, and let $\widehat C^i_k$ and $C^i_k$ be given in \eqref{eq:cons_i} and \eqref{eq:globalcons}, respectively. Then, for every $k\in\N$
$$\frac1N\sum_{i=1}^N\E\big[\|\widehat C^i_k-C^i_k\|^2\big]\leq\frac{\alpha^2M\sigma_f^2D_{\mathcal{X}}}{N}\,\vartheta^{2k}.$$
\end{lemma}

\begin{proof}

Let $\mathcal F_k$ be the $\sigma$-algebra generated by the initial positions of the agents, $(\xi^\ell_{t,s})_{t<k}$ and
$(\omega^\ell_t)_{t<k}$. By Assumption \ref{a:f_hat} we have $\E[\mathcal E_{k,\ell}^2\ |\ \mathcal F_k]\leq \sigma_f^2$ for every $\ell$, and
therefore
\begin{equation*}
    \E\Big[\max_{\ell\in\calM_i}\mathcal E_{k,\ell}^2\,\Big|\,\mathcal F_k\Big]\leq\sum_{\ell\in\calM_i}\E\big[\mathcal E_{k,\ell}^2\ |\ \mathcal F_k\big]\leq M\sigma_f^2.
\end{equation*}

Moreover, applying \cite[Lemma 3.3]{bellavia2025} to the swarm $(X^\ell_k)_{\ell\in\calM_i}$ held by node $i$, we have that for every $s=1,\dots,d$ and every $k\in\N$
$$|\widehat\cons^i_{k,s}-\cons^i_{k,s}|\leq \alpha \max_{\ell\in\calM_i}|\mathcal{E}_{k,\ell}|\,\diam( \mathcal{X}^{s,i}_k).$$

Using the law of total expectation, these inequalities, and Lemma \ref{lemma:local_moments}, we have
\begin{equation*}\begin{aligned}
    \E\big[\|\widehat C^i_k- C^i_k\|^2\big]
    &\leq\alpha^2\sum_{s=1}^\dims\E\Big[\diam( \mathcal{X}^{s,i}_k)^2\,\E\Big[\max_{\ell\in\calM_i}\mathcal E_{k,\ell}^2\,\Big|\,\mathcal F_k\Big]\Big]
    \leq\alpha^2M\sigma_f^2\sum_{s=1}^\dims\E\big[\diam( \mathcal{X}^{s,i}_k)^2\big]\\
    &\leq\alpha^2M\sigma_f^2\,\vartheta^{2k}\sum_{s=1}^\dims\E\big[\diam( \mathcal{X}^{s,i}_0)^2\big].
\end{aligned}\end{equation*}
Averaging over $i=1,\dots,N$ and using the definition of $D_{\mathcal{X}}$ we get the thesis.
\end{proof}
Combining Lemma \ref{lemma:global_consensus} and Lemma \ref{lemmadiffxalpha}, we can now bound the distance of the tracking variables from $x^*$ in terms of the distance of the exact global consensus point $C_k$ from $x^*$, up to a term that decays geometrically.
\begin{lemma}\label{lemma:Zbound}
Let the assumptions of Lemmas \ref{lemma:global_consensus} and \ref{lemmadiffxalpha} hold. Then for every $k\in\N$
\begin{equation*}
    \frac1N\sum_{i=1}^N\E\big[\|Z^i_k-x^*\|^2\big]\leq4\,\E\big[\| C_k-x^*\|^2\big]+\Lambda_{Z}\,\varrho^{2k},
    \qquad \Lambda_{Z}=4\left(\frac{\alpha^2M\sigma_f^2D_{\mathcal{X}}}{N}+\sum_{s=1}^\dims\Big(\frac{\Lambda_{\mathcal{Z},s}^2}{N}+\Lambda_{\mathcal{X},s}^2\Big)\right).
\end{equation*}
\end{lemma}

\begin{proof}

By \eqref{def:averages} we have 
\begin{equation*}
    \begin{aligned}
        \|Z^i_k-x^*\|^2 &\leq 4\| Z^i_k-\overline{Z}_k\|^2 +4\left\|\frac{1}{N}\sum_{j=1}^N \widehat{ C}_k^j - \frac{1}{N}\sum_{j=1}^N  C_k^j  \right\|^2 + 4\left\|\frac{1}{N}\sum_{j=1}^N  C_k^j -  C_k\right\|^2 +4\| C_k-x^*\|^2 \\
        &\leq 4\| Z^i_k-\overline{Z}_k\|^2 + \frac{4}{N}\sum_{j=1}^N\|\widehat  C_k^j- C^j_k\|^2 + \frac{4}{N}\sum_{j=1}^N\| C^j_k- C_k\|^2 +4\| C_k-x^*\|^2.
    \end{aligned}
\end{equation*}
This implies 
\begin{equation}\label{aux_Z}
\begin{aligned}
    \frac1N\sum_{i=1}^N\E[\|Z^i_k-x^*\|^2]&\leq   \frac4N\left (\sum_{i=1}^N\E[\|Z^i_k-\overline Z_k\|^2]+\sum_{j=1}^N\E[\|\widehat C^j_k- C^j_k\|^2]+\sum_{j=1}^N\E[\| C^j_k- C_k\|^2]\right )\\&\quad+4\E[\| C_k-x^*\|^2].
\end{aligned}
\end{equation}
We now bound the first three terms on the RHS.
By definition of $\mathcal{Z}^s_k$ and Lemma \ref{lemma:global_consensus} we have
\begin{equation}\label{aux_Zbar}
    \sum_{i=1}^N\E[\|Z^i_k-\overline Z_k\|^2] = \sum_{i=1}^N\sum_{s=1}^\dims\E[((I-J)\mathcal{Z}^s_k)_i^2] =  \sum_{s=1}^\dims\E[\|(I-J)\mathcal{Z}^s_k\|^2]\leq \varrho^{2k}\sum_{s=1}^d\Lambda_{\mathcal{Z},s}^2.
\end{equation}
By Lemma \ref{lemmadiffxalpha}, since $\vartheta\leq\varrho,$
 \begin{equation}\label{aux_c}
 \sum_{j=1}^N\E[\|\widehat C^j_k- C^j_k\|^2]\leq \alpha^2M\sigma_f^2D_{\mathcal{X}}\vartheta^{2k}\leq\alpha^2M\sigma_f^2D_{\mathcal{X}}\varrho^{2k}.
\end{equation}
Since $ C^j_k$ and $ C_k$ are both convex
combinations of $(X^\ell_k)_{\ell=1}^{NM}$, by Lemma \ref{lemma:global_consensus}, we have
\begin{equation}\label{aux_c2}
 \sum_{j=1}^N\E[\| C^j_k- C_k\|^2]\leq N\sum_{s=1}^\dims \E[\diam( \mathcal{X}^s_k)^2]\leq N\varrho^{2k}\sum_{s=1}^\dims\Lambda_{\mathcal{X},s}^2.
\end{equation}
Using \eqref{aux_Zbar}, \eqref{aux_c} and \eqref{aux_c2} in \eqref{aux_Z},
we get the thesis.
\end{proof}

We can now prove our main result.

\begin{theorem}\label{thm:convergence}
Let Assumptions \ref{a:f}--\ref{a:W} hold. Moreover, let $\lambda\in(0,1)$ and $\sigma>0$ be such that $\vartheta<1$. 
Given $\varepsilon>0$ and $\tau\in(0,1)$ such that $1-\lambda(1-\tau)\neq\varrho^2$, with $\varrho$ defined in Lemma \ref{lemma:global_consensus}, let us define
\begin{equation}\label{eq:thm_constants}
\begin{aligned}
    &\mu=\max\{1-\lambda(1-\tau),\varrho^2\},\qquad
    \Lambda_{V}=\E[V_0]+\frac{\left(\frac{\sigma^2D_{\mathcal{X}}}{N}+\lambda \Lambda_{Z}\right)}{|1-\lambda(1-\tau)-\varrho^2|},\qquad
    k^*=\max\Big\{0,\Big\lceil\log_{1/\mu}\Big(\frac{\Lambda_{V}}{\varepsilon}\Big)\Big\rceil\Big\},
\end{aligned}\end{equation}
with $\Lambda_{Z}$ given in Lemma \ref{lemma:Zbound}, and
\begin{equation}\label{eq:qr_eps}
    q_\varepsilon=\frac12\min\left\{f_\infty,\Big(\frac\beta2\Big(\frac{\tau\varepsilon}{4}\Big)^{1/2}\Big)^{1/\nu}\right\},\qquad
    r_\varepsilon=\max\{s\in[0,R_0]\ |\ f(x)\leq q_\varepsilon\ \ \forall x\in\mathcal B_s\}.
\end{equation}
Assume that for every $k\leq k^*$ we have
\begin{equation}\label{eq:deltak}
    \mathbb P\big(|\mathcal{I}_{k,r_\varepsilon}|=0\big)\leq\delta_k\qquad\text{with}\qquad
    \delta_k\leq\frac{\tau^2}{64\,N^2M^2}\,\frac{\E[V_k]^2}{\E[V_k^2]} .
\end{equation}
If
\begin{equation}\label{eq:alphacond}
    \alpha\geq\frac{1}{q_\varepsilon}\log\Big(\frac{16NM}{\sqrt\tau}\Big),
\end{equation}
then there exists $\bar k\leq k^*$ such that $\E[V_{\bar k}]\leq\varepsilon$. Moreover, $\E[V_k]\leq \Lambda_{V}\,\mu^k$ for every $k=0,\dots,\bar k$.
\end{theorem}

\begin{remark}
    The rate is given by $\mu\geq\varrho^2 = \max\{\vartheta^2, \kappa ^2\}$, where $\kappa $ is larger than the spectral radius of $A_\lambda P $, which in turn depends on the eigenvalues $\gamma_2,\dots,\gamma_N$ of the mixing matrix $W$.  In particular we can observe that if $\gamma_2$ approaches 1 then $\kappa $ (and therefore the rate $\mu$) also approaches 1. Since the spectral gap $1-\gamma_2$ measures the connectivity of the underlying communication network $\mathcal{G}$, this reflects the fact that the rate guaranteed by the theorem is worse when the network is poorly connected.
    \end{remark}

\begin{proof}
Let $k_0$ be the smallest $k\in\N$ such that
\begin{equation}\label{eq:kbar}
    \E[V_k]\leq\varepsilon\qquad\text{or}\qquad\E\big[\| C_k-x^*\|^2\big]>\frac\tau4\E[V_k].
\end{equation}
For every $k<k_0$ we have $\E[\| C_k-x^*\|^2]\leq\frac\tau4\E[V_k]$, and therefore,
by Lemma \ref{lemma:Vk}, Lemma \ref{lemma:Zbound} and $\vartheta\leq\varrho$,
\begin{equation*}\begin{aligned}
    \E[V_{k+1}]&\leq(1-\lambda)\E[V_k]+4\lambda\E\big[\| C_k-x^*\|^2\big]+\lambda \Lambda_{Z}\varrho^{2k}+\frac{\sigma^2D_{\mathcal{X}}}{N}\varrho^{2k}
    \\&\leq(1-\lambda(1-\tau))\E[V_k]+\left(\frac{\sigma^2D_{\mathcal{X}}}{N}+\lambda \Lambda_{Z}\right)\varrho^{2k}.
\end{aligned}\end{equation*}
Applying this inequality recursively we get, for every $k\leq k_0$
\begin{equation}\label{eq:Vdecay}\begin{aligned}
    \E[V_k]&\leq(1-\lambda(1-\tau))^k\E[V_0]+\left(\frac{\sigma^2D_{\mathcal{X}}}{N}+\lambda \Lambda_{Z}\right)\sum_{t=0}^{k-1}(1-\lambda(1-\tau))^{k-1-t}\varrho^{2t}\\&=(1-\lambda(1-\tau))^k\E[V_0]+\left(\frac{\sigma^2D_{\mathcal{X}}}{N}+\lambda \Lambda_{Z}\right)\frac{(1-\lambda(1-\tau))^k-\varrho^{2k}}{(1-\lambda(1-\tau))-\varrho^2}\leq \Lambda_{V}\,\mu^k.
\end{aligned}\end{equation}
If $k^*\leq k_0$,  by \eqref{eq:Vdecay} and the definition of $k^*$ it holds $\E[V_{k^*}]\leq \Lambda_{V}\mu^{k^*}\leq\varepsilon$, therefore we have the thesis and $\bar k= k^*.$\\
If $k^*> k_0$ and $\E[V_{ k_0}]\leq \varepsilon$, we also have the thesis and $\bar k= k_0.$ \\
We prove that the remaining case, namely $\E[V_{ k_0}]>\varepsilon$ and $\E\big[\| C_{ k_0}-x^*\|^2\big]>\frac\tau4\E[V_{ k_0}]$ cannot happen, provided that $\alpha$ fulfills \eqref{eq:alphacond} and that \eqref{eq:deltak} holds. 

By the law of total expectation we have
\begin{equation}\label{eq:totalexp}\begin{aligned}
    \E[\| C_{k_0}-x^*\|^2] &= \E[\| C_{k_0}-x^*\|^2 \ |\  |\mathcal{I}_{k_0,r_\varepsilon}|>0]\mathbb P(|\mathcal{I}_{k_0,r_\varepsilon}|>0)
    + \E[\| C_{k_0}-x^*\|^2 \ | \ |\mathcal{I}_{k_0,r_\varepsilon}|=0]\mathbb P(|\mathcal{I}_{k_0,r_\varepsilon}|=0).
\end{aligned}\end{equation}
We now bound the two terms on the RHS.

Let us consider the first term. Applying part $(i)$ of Lemma \ref{lemma:laplaceprinciple}, with $q=q_\varepsilon$ and $r=r_\varepsilon$, and using 
\begin{equation*}
\frac1\beta(q_\varepsilon+f_{r_\varepsilon})^\nu\leq\frac1\beta(2q_\varepsilon)^\nu\leq\frac12\Big(\frac{\tau\varepsilon}4\Big)^{1/2}<\frac{\sqrt\tau}{4}\E[V_{k_0}]^{1/2},
\end{equation*}
we have that, if $|\mathcal{I}_{k_0,r_\varepsilon}|>0$, then
\begin{equation*}
\begin{aligned}
    \| C_{k_0}-x^*\|^2&\leq\Big(\frac{\sqrt\tau}{4}\E[V_{k_0}]^{1/2}+NMe^{-\alpha q_\varepsilon}V_{k_0}^{1/2}\Big)^2
    \\&=\frac{\tau}{16}\E[V_{k_0}]+\frac{\sqrt\tau}{2}NMe^{-\alpha q_\varepsilon}\E[V_{k_0}]^{1/2}V_{k_0}^{1/2}+N^2M^2e^{-2\alpha q_\varepsilon}V_{k_0}.
\end{aligned}
\end{equation*}
Taking the conditional expectation this yields
\begin{equation*}\begin{aligned}
    \E[\| C_{k_0}-x^*\|^2 \ |\ |\mathcal{I}_{k_0,r_\varepsilon}|>0]
    &\leq\frac{\tau}{16}\E[V_{k_0}]+\frac{\sqrt\tau}{2}NMe^{-\alpha q_\varepsilon}\E[V_{k_0}]^{1/2}\E[V_{k_0}\ |\ |\mathcal{I}_{k_0,r_\varepsilon}|>0]^{1/2}\\
    &\quad+N^2M^2e^{-2\alpha q_\varepsilon}\E[V_{k_0}\ |\ |\mathcal{I}_{k_0,r_\varepsilon}|>0],
\end{aligned}\end{equation*}
and therefore, by the assumption on $\alpha$, and the fact that for any non-negative random variable $y$ and any event $\mathcal{S}$ it holds 
\begin{equation}\label{aux_ineqprob}
    \E[y|\mathcal{S}]\mathbb{P}(\mathcal{S})\leq \E[y],
\end{equation}
we have 
\begin{equation}\label{aux_case1}
    \E[\| C_{k_0}-x^*\|^2 \ |\ |\mathcal{I}_{k_0,r_\varepsilon}|>0]\,\mathbb P(|\mathcal{I}_{k_0,r_\varepsilon}|>0)
    \leq\Big(\frac{\tau}{16}+\frac{\sqrt\tau}{2}NMe^{-\alpha q_\varepsilon}+N^2M^2e^{-2\alpha q_\varepsilon}\Big)\E[V_{k_0}]
    \leq\frac\tau8\E[V_{k_0}].
\end{equation}

By part $(ii)$ of Lemma \ref{lemma:laplaceprinciple}, Jensen's inequality, the assumption on $\delta_k$, and \eqref{aux_ineqprob} we have
\begin{equation}\label{aux_case2}\begin{aligned}
    \E[\| C_{k_0}-x^*\|^2 \ |\ |\mathcal{I}_{k_0,r_\varepsilon}|=0]\,\mathbb P(|\mathcal{I}_{k_0,r_\varepsilon}|=0)
    &\leq NM\,\E[V_{k_0} \ |\   |\mathcal{I}_{k_0,r_\varepsilon}|=0]\,\mathbb P(|\mathcal{I}_{k_0,r_\varepsilon}|=0)\\
    &\leq NM\,\E[V_{k_0}^2 \ |\   |\mathcal{I}_{k_0,r_\varepsilon}|=0]^{1/2}\,\mathbb P(|\mathcal{I}_{k_0,r_\varepsilon}|=0)\\
    &= NM\,\Big(\E[V_{k_0}^2 \ |\   |\mathcal{I}_{k_0,r_\varepsilon}|=0]\,\mathbb P(|\mathcal{I}_{k_0,r_\varepsilon}|=0)\Big)^{1/2}\,\mathbb P(|\mathcal{I}_{k_0,r_\varepsilon}|=0)^{1/2}\\
    &\leq NM\,\E[V_{k_0}^2]^{1/2}\,\delta_{k_0}^{1/2}\leq\frac\tau8\E[V_{k_0}].
\end{aligned}\end{equation}
Using \eqref{aux_case1} and \eqref{aux_case2} in \eqref{eq:totalexp} we get $\E[\| C_{k_0}-x^*\|^2]\leq\frac\tau4\E[V_{k_0}]$,
which gives a contradiction, and the proof is complete.
\end{proof}

\subsubsection{On assumption \eqref{eq:deltak}}\label{sec:hitting}
Assumption \eqref{eq:deltak} is fundamental in order to apply the quantitative Laplace principle (Lemma \ref{lemma:laplaceprinciple}) to ensure $ C_k\approx x^*$. Yet,
assumption \eqref{eq:deltak} can appear quite strong, as it requires that there is a high probability that a particle hits a neighborhood of the global minimizer in finite time. When this happens, it is no longer surprising that, by the consensus drift, the entire population will also converge to the same neighborhood. While Theorem \ref{thm:convergence} does not guarantee such a hitting time with high probability, it shows anyway that the dynamics is contractive, namely that $\E[V_k]\leq \Lambda_{V}\,\mu^k$ holds for every $k=0,\dots,\bar k$.
In what follows, we would like to provide a justification of Assumption \eqref{eq:deltak}, at least for a simplified model.
As Assumption \eqref{eq:deltak} is needed in order to apply the quantitative Laplace principle, the latter can no longer be used to ensure $ C_k\approx x^*$. Hence, we argue by geometrical considerations that if $x^*$ is in the convex hull of the particles, the function $f$ is sufficiently smooth, and $\alpha$ is not too large, then $ C_k\approx x^*$ is a legitimate assumption, even without invoking the Laplace principle or requiring that any particle actually lies in the vicinity of $x^*$. 
Accordingly, in the following, we will consider the simplified model for which $\widehat C^i_{k}= C_k:=x^*$. With this (admittedly strong) assumption, the model reduces to a totally decoupled  dynamics of multiple independent particles, no matter what their graph connectivity $W$ is. Now, we  prove that \eqref{eq:deltak} holds for this simplified model.

More formally, let us consider the decoupled dynamics 
\begin{equation}\label{simple_dyn}
    X^\ell_{k+1} = X^\ell_k - \lambda(X^\ell_k-x^*)+\sigma \sum_{s=1}^\dims (x^\ell_{k,s}- x^*_s) \xi^\ell_{k,s}\mathbf{e}_s,
\end{equation}
where ${(X^\ell_0)}_{\ell=1}^{NM}$ are i.i.d. random variables with law $\rho_0$, and for every $s=1:d,\ \ell=1:NM,\ k\in\N$ the random variables $\xi^\ell_{k,s}$ are independent and normally distributed. 

\begin{lemma}\label{lemma:prob_LHS}
Let us assume $(X^\ell_k)$ satisfies \eqref{simple_dyn} with $\lambda\in(0,1),\ \sigma>0$ such that $\vartheta<1$. Then the following hold. 
\begin{enumerate}[label=(\roman*)]
    \item For every $s=1:d,\ \ell=1:NM,\ k\in\N$
    $$(x^\ell_{k,s}-x^*_s) = (x^\ell_{0,s}-x^*_s)\prod_{t=0}^{k-1}(1-\lambda+\sigma\xi^\ell_{t,s}).$$
    \item For every $k\in\N$ the random variables $(X^\ell_k)_{\ell=1}^{NM}$ are i.i.d. Moreover, for every $s=1:d,\ k\in\N$ the random variables $(x^\ell_{k,s}-x^*_s)_{\ell=1}^{NM}$ are i.i.d.
    \item For every $k\in\N$ and every $r>0$
            $$\mathbb{P}(|\mathcal{I}_{k,r}|=0)=\left(1-\mathbb{P}(\|X^1_k-x^*\|\leq r)\right)^{NM}.$$
    \item For every $k\in\N$ and every $r>0$
    \begin{equation*}
        \begin{aligned}
            \mathbb{P}(\|X^1_k-x^*\|\leq r)&\geq \mathbb{P}(\|X^1_0-x^*\|\leq r)\max\left\{\eta^k,\ 1-((1-\lambda)^2+\sigma^2)^k\right\}^d\\
            &\geq \mathbb{P}(\|X^1_0-x^*\|\leq r)\min\left\{\frac12,\ \eta^T\right\}^d=:\underline{p},
        \end{aligned}
    \end{equation*}
    with $\eta= 1-\frac{\sqrt{2}\sigma}{\sqrt{\pi}\lambda}e^{-\frac{\lambda^2}{2\sigma^2}},$ and $T = \frac{\log(2)}{|\log((1-\lambda)^2+\sigma^2)|}$.
\end{enumerate}
\end{lemma}
\begin{proof}
    By \eqref{simple_dyn} we have that for every $s=1:d,\ \ell=1:NM,\ k\in\N$
    $$(x^\ell_{k,s}-x^*_s) = (1-\lambda+\sigma\xi^\ell_{k-1,s})(x^\ell_{k-1,s}-x^*_s).$$
    Recursively applying this identity, we get $(i)$.
    Part $(ii)$ of the statement follows directly from $(i)$ and the fact that the random variables $\xi^\ell_{t,s}$ are i.i.d. and independent of $(X^\ell_0)_{\ell=1}^{NM},$ for every $\ell,\ t$ and $s.$
    By definition of $\mathcal{I}_{k,r}$, using the fact that $(X^\ell_k)_{\ell=1}^{NM}$ are i.i.d. we have
    \begin{equation*}
        \begin{aligned}
            \mathbb{P}(|\mathcal{I}_{k,r}|=0) &= \mathbb{P}(\|X^\ell_k-x^*\|> r\ \forall\ \ell=1:NM) = \prod_{\ell=1}^{NM}\mathbb{P}(\|X^\ell_k-x^*\|> r)\\
            &=\mathbb{P}(\|X^1_k-x^*\|> r)^{NM} = \left(1-\mathbb{P}(\|X^1_k-x^*\|\leq r)\right)^{NM},
        \end{aligned}
    \end{equation*}
    which is $(iii).$
    We now prove part $(iv)$ of the statement. Let us denote by $\mathcal{S}_k$ the event
    $$\mathcal{S}_k = \{\|X^1_0-x^*\|\leq r\}\cap \left\{ \prod_{t=0}^{k-1}\left|1-\lambda+\sigma\xi^1_{t,s}\right|\leq 1\ \forall\ s=1:d\right\}.$$
    If $\mathcal{S}_k$ is true, by part $(i)$ and the definition of $\mathcal{S}_k$, we have
    \begin{equation*}
        \begin{aligned}\|X^1_k-x^*\|^2&=\sum_{s=1}^d|x^1_{k,s}-x^*_s|^2= \sum_{s=1}^d|x^1_{0,s}-x^*_s|^2\prod_{t=0}^{k-1}\left|1-\lambda+\sigma\xi^1_{t,s}\right|^2\\
    &\leq\sum_{s=1}^d|x^1_{0,s}-x^*_s|^2\leq r^2.\end{aligned}
    \end{equation*}
    This, together with the fact that $\xi^1_{t,s}$ are i.i.d., implies 
    \begin{equation}\label{pq}\mathbb{P}(\|X^1_k-x^*\|\leq r)\geq\mathbb{P}(\mathcal{S}_k) = \mathbb{P}(\|X^1_0-x^*\|\leq r)q_k^d,\end{equation}
    with $q_k = \mathbb{P}\left(\prod_{t=0}^{k-1}\left|1-\lambda+\sigma\xi^1_{t,s}\right|\leq 1\right). $ We now bound $q_k.$
    On the one hand since the variables $\xi$ are independent and normally distributed, we have 
      \begin{equation}\label{q_1}
        \begin{aligned}
        q_k &\geq \mathbb{P}\left(\left|1-\lambda+\sigma\xi^1_{t,s}\right|\leq 1\ \forall t=0:k-1\right) = \mathbb{P}\left(\left|1-\lambda+\sigma\xi^1_{0,s}\right|\leq 1\right)^k\\
        & = \mathbb{P}\left(-\frac{2-\lambda}{\sigma}\leq\xi^1_{0,s}\leq\frac{\lambda}{\sigma}\right)^k\geq \mathbb{P}\left(\left|\xi^1_{0,s}\right|\leq \frac{\lambda}{\sigma}\right)^k \\
        &= \left(1- \mathbb{P}\left(\left|\xi^1_{0,s}\right|> \frac{\lambda}{\sigma}\right)\right)^k \geq \left(1-\frac{\sqrt{2}\sigma}{\sqrt{\pi}\lambda}e^{-\frac{\lambda^2}{2\sigma^2}}\right)^k=\eta^k.
        \end{aligned}
    \end{equation}
    On the other hand, by the Markov inequality and again the fact that the variables $\xi$ are i.i.d., we have 
    \begin{equation}\label{q_2}
        \begin{aligned}
            1-q_k &= \mathbb{P}\left(\prod_{t=0}^{k-1}\left|1-\lambda+\sigma\xi^1_{t,s}\right|^2> 1\right)  \leq  \E\left[\prod_{t=0}^{k-1}\left(1-\lambda+\sigma\xi^1_{t,s}\right)^2\right]\\
            &=\E\left[\left(1-\lambda+\sigma\xi^1_{0,s}\right)^2\right]^k = ((1-\lambda)^2+\sigma^2)^k.
        \end{aligned}
    \end{equation}
    Putting together \eqref{q_1} and \eqref{q_2} we have 
    \begin{equation}\label{q_k}
        q_k\geq\max\left\{\eta^k,\ 1-((1-\lambda)^2+\sigma^2)^k\right\}
    \end{equation}
     which, together with \eqref{pq} gives the first inequality in part $(iv)$. To conclude the proof, it is enough to notice that if $k\geq T$ then by \eqref{q_k}
     $q_k\geq 1-((1-\lambda)^2+\sigma^2)^k\geq \frac12,$
     while if $k<T$ then 
     $q_k\geq\eta^T$.   
\end{proof}

\begin{lemma}\label{lemma:prob_RHS}
Let us assume $(X^\ell_k)$ satisfies \eqref{simple_dyn}. Then for every $k\in\N$ 
    $$\frac{\E[V_k]^2}{\E[V_k^2]}\geq \left(\frac{((1-\lambda)^2+\sigma^2)^2}{(1-\lambda)^4+6(1-\lambda)^2\sigma^2+3\sigma^4}\right)^k\frac{\E[\|X^1_0-x^*\|^2]^2}{\E[\|X^1_0-x^*\|^4]}.$$
    \end{lemma}
    \begin{proof}
    By definition of $V_k$ and part $(ii)$ of Lemma \ref{lemma:prob_LHS}, we have
        \begin{equation*}
                \E[V_k]= \frac{1}{NM}\sum_{\ell=1}^{NM}\E\left[\|X^\ell_k-x^*\|^2\right] = \E[\|X^1_k-x^*\|^2],
        \end{equation*}
    and 
    \begin{equation*}
            \begin{aligned}
                \E[V_k^2] &= \frac{1}{(NM)^2}\E\left[\left(\sum_{\ell=1}^{NM}\|X^\ell_k-x^*\|^2\right)^2\right] \\
                &= \frac{1}{(NM)^2}\E\left[\sum_{\ell=1}^{NM}\|X^\ell_k-x^*\|^4 +
                \sum_{\ell=1}^{NM}\sum_{m\neq\ell}\|X^\ell_k-x^*\|^2\|X^m_k-x^*\|^2\right]\\ & = \frac{1}{NM}\E[\|X^1_k-x^*\|^4] + \frac{NM-1}{NM} \E[\|X^1_k-x^*\|^2]^2.
            \end{aligned}
        \end{equation*}
    This, using the fact that $\E[\|X^1_k-x^*\|^4]\geq\E[\|X^1_k-x^*\|^2]^2$, implies
    \begin{equation}\label{ratio_V}
        \begin{aligned}
            \frac{\E[V_k]^2}{\E[V_k^2]} = \frac{NM}{\frac{\E[\|X^1_k-x^*\|^4]}{\E[\|X^1_k-x^*\|^2]^2}-1+NM}
            \geq \frac{\E[\|X^1_k-x^*\|^2]^2}{\E[\|X^1_k-x^*\|^4]}. 
        \end{aligned}
    \end{equation}
    We now bound the right hand side. 
    Let us first notice that, since the variables $\xi$ are independent and normally distributed, for every $k\in\N$ and every $s=1,\dots,d$ we have 
    \begin{equation*}
            \E\left[\prod_{t=0}^{k-1}(1-\lambda+\sigma\xi^1_{t,s})^4\right] = \left((1-\lambda)^4+6(1-\lambda)^2\sigma^2+3\sigma^4\right)^k,
    \end{equation*}
    and for every $k\in\N$ and every $s,z=1,\dots,d$ with $s\neq z$
    \begin{equation*}
        \begin{aligned}
            &\E\left[\prod_{t=0}^{k-1}(1-\lambda+\sigma\xi^1_{t,s})^2\prod_{t=0}^{k-1}(1-\lambda+\sigma\xi^1_{t,z})^2\right] =     \E\left[\prod_{t=0}^{k-1}(1-\lambda+\sigma\xi^1_{t,s})^2\right]^2\\ &=\left((1-\lambda)^2+\sigma^2\right)^{2k}\leq \left((1-\lambda)^4+6(1-\lambda)^2\sigma^2+3\sigma^4\right)^k.
        \end{aligned}
    \end{equation*}
   These, together with part $(i)$ of Lemma \ref{lemma:prob_LHS}, imply 
   \begin{equation}\label{ratio_E4}
       \begin{aligned}
           &\E[\|X^1_k-x^*\|^4]= \E\left[\left(\sum_{s=1}^d(x^1_{k,s}-x^*_s)^2\right)^2\right] \\
           &= \E\left[\sum_{s=1}^d(x^1_{k,s}-x^*_s)^4+ \sum_{s=1}^d\sum_{z\neq s}(x^1_{k,s}-x^*_s)^2(x^1_{k,z}-x^*_z)^2\right] \\
           &= \E\left[\sum_{s=1}^d(x^1_{0,s}-x^*_s)^4\prod_{t=0}^{k-1}(1-\lambda+\sigma\xi^1_{t,s})^4\right]\\
           &+\E\left[\sum_{s=1}^d\sum_{z\neq s}(x^1_{0,s}-x^*_s)^2(x^1_{0,z}-x^*_z)^2\prod_{t=0}^{k-1}(1-\lambda+\sigma\xi^1_{t,s})^2\prod_{t=0}^{k-1}(1-\lambda+\sigma\xi^1_{t,z})^2\right]\\
           &\leq \left((1-\lambda)^4+6(1-\lambda)^2\sigma^2+3\sigma^4\right)^k\E[\|X^1_0-x^*\|^4].
       \end{aligned}
   \end{equation}
   Analogously 
      \begin{equation}\label{ratio_E2}
       \begin{aligned}
           &\E[\|X^1_k-x^*\|^2]^2 =  \left((1-\lambda)^2+\sigma^2\right)^{2k}\E[\|X^1_0-x^*\|^2]^2.
       \end{aligned}
   \end{equation}
   Using \eqref{ratio_E4} and \eqref{ratio_E2} in \eqref{ratio_V} we get the thesis.
    \end{proof}

By Lemma \ref{lemma:prob_LHS}, parts $(iii)$ and $(iv)$ we have that for every $k\leq k^*$
$\mathbb{P}(|\mathcal{I}_{k,r_\eps}|=0)\leq (1-\underline{p})^{NM}$,
while by Lemma \ref{lemma:prob_RHS}, the right hand side of \eqref{eq:deltak} satisfies
$$\frac{\tau^2}{64N^2M^2}\frac{\E[V_k]^2}{\E[V_k^2]}\geq \frac{\tau^2}{64N^2M^2}\left(\frac{((1-\lambda)^2+\sigma^2)^2}{(1-\lambda)^4+6(1-\lambda)^2\sigma^2+3\sigma^4}\right)^{k^*}\frac{\E[\|X^1_0-x^*\|^2]^2}{\E[\|X^1_0-x^*\|^4]}=: \frac{\Delta}{N^2M^2}$$
By putting together these two inequalities and enforcing $(1- \underline{p})^{NM} \leq \frac{\Delta}{N^2M^2}$ we have that if 
\begin{equation}
    \label{eq:Nparticles}
NM\geq \frac{2}{|\log(1-\underline{p})|}\log\left(\frac{16}{\Delta\,|\log(1-\underline{p})|^2}\right)
\end{equation}
then inequality \eqref{eq:deltak} holds for every $k\leq k^*$.\\
Notice that the number of particles above depends on: the initial distribution of the particles (through the second and fourth moment of $\|X^1_0-x^*\|$ in $\Delta$, and the term $\mathbb{P}(\|X^1_0-x^*\|\leq r)$ in $\underline{p}$); the parameters $\lambda$ and $\sigma$; the dimension $d$, through the term $\min\left\{\frac12,\ \eta^T\right\}^d$ in $\underline{p}$.
\begin{remark} $(i)$ The estimate \eqref{eq:Nparticles} above implies the need for an exponentially large number $N M$ of particles with respect to the dimensionality $d$. This 
should not be surprising, as our assumptions on $f$ are minimal and the class of optimization problems they comprise can include NP-hard ones (intrinsically affected by the curse of dimensionality). Yet, there are classes of functions for which CBO methods do not have exponential complexity \cite{bonandin2026exploiting}.\\
$(ii)$    An adaptation of the above argument to the original model \eqref{eq:CBOpar} allowing for $\widehat C^i_{k} \neq x^*$ is technically very hard. It would require the use of the quantitative Laplace principle (Lemma \ref{lemma:laplaceprinciple}), whose applicability  we would like to justify in the first place. Hence, there is a circular argument that we cannot break. We conclude that a full convergence argument for the model \eqref{eq:CBOpar} is currently beyond the scope of the current techniques we have at our disposal. We should mention that the continuous-in-time approach as in \cite{fornasier2024} could allow for a more complete proof of convergence, but it would also require the analysis of a PDE with stochastic coefficients -- as we are considering the stochastic objective $\widehat f$ to define the consensus points -  describing the mean-field law of the particle system. This analysis is completely open to this date and is beyond the scope of the present paper.
\end{remark}

\section{Numerical examples}\label{sec:numerics}

In this section we test the performance of the proposed decentralized CBO algorithm on different benchmark problems. The code for the experiments is available at \href{https://github.com/TimRoith/DecentralizedCBO}{github.com/TimRoith/DecentralizedCBO} and is based on the \texttt{CBXPy} package from \cite{Bailo_CBX_Python_and_2024}.

\subsection{Notes on the implementation and heuristics}
\begin{algorithm}
\caption{Decentralized Consensus-based optimization}\label{alg:CBO}
\textbf{Input:}  Initial ensemble $X_{0}$, hyperparameters $\sigma\geq 0, \Delta t>0$, mixing matrix $W$
\begin{algorithmic}[1]
\State Initialize $\widehat{C}_{0}^i = \textbf{ComputeConsensus}(X_{0}^{\mathcal{M}_i},\alpha_0^i)$, 
$Z_{0}^i =\widehat{C}_{0}^i$ for $i=1:N$
\While{\textbf{Not Terminate()}}
    \ParFor{$i=1:N$}
    \State $\tilde{X}_{k+1}^{\mathcal{M}_i} =X_{k}^{\mathcal{M}_i} - \Delta t\, (X_{k}^{\mathcal{M}_i} -  Z_k^i) + \sigma \,\textbf{Noise}(X_{k}^{\mathcal{M}_i} -  \widehat{C}^i_{k}, \Delta t)$
    \State \textbf{PostStep}()
    \State $\widehat{C}_{k+1}^i = \textbf{ComputeConsensus}(X_{k+1}^{\mathcal{M}_i}, \alpha_{k}^i)$
    \EndParFor
    \For{$i=1:N$}
    \State $Z_{k+1}^i = \sum_{j\in\mathcal{N}_i} W_{ij} Z_{k}^j + (\widehat{C}_{k+1}^i  - \widehat{C}_{k}^i)$
    \EndFor
\EndWhile
\end{algorithmic}
\end{algorithm}%
We briefly describe the numerical implementation and specify, which heuristics are employed in the experiments. The pseudo-code for the proposed algorithm in \eqref{eq:CBOpar} is presented in \cref{alg:CBO}. Therein the function \textbf{Terminate()} denotes a generic termination criterion, which we detail in the following sections. The function \textbf{ComputeConsensus} computes the consensus point as specified in \eqref{eq:cons_i} in a numerically stable way. In our case this is realized via the LogSumExp function, as in \cite[Alg.~2]{bungert2025mirrorcbo}. 
Beyond that the role of the function \textbf{PostStep} and the sequence of parameters $(\alpha^i_k)$ will be explained in \cref{sec:heuristics}. 
We slightly modify the notation in \eqref{eq:CBOpar} to have the usual numerical interpretation as a time-discretized scheme. For this we set $\lambda = \Delta t$ and multiply an additional factor $\sqrt{\Delta t}$ to the noise part. Concretely, the function $\textbf{Noise}$ in \cref{alg:CBO} will be the anisotropic noise model from \cite{carrillojinlizhu2021}, meaning that $\textbf{Noise}(z, \Delta t) = \sqrt{\Delta t}\, z \odot \xi$, where $\xi\sim\mathcal{N}(0,I)$. We refer to \cite{fiedler2026consensus,bonandin2026exploiting} for a recent discussion on the numerical behavior of this model and use it throughout our experiments.

We highlight that in \cref{alg:CBO}, the first loop in lines 3--6 can indeed be parallelized, since there is no interaction across nodes. The synchronization happens in lines 7--8, where each node communicates with its neighbors. 
The total number of communications per step is $\sum_{i=1}^N\left|\mathcal{N}_i\setminus\{i\} \right|$, where each communication exchanges the respective tracking variable of dimension $d$.

\begin{remark}
The communication could also be realized in an asynchronous fashion, which would remove the forced synchronization step. In this case each node sends the tracking variable to its neighbors, once it computed it. However, both the analysis and the numerical study of such a scheme go beyond the scope of the present work. Regarding the communication cost per step, we want to remark that this is not purely determined by the number of communications, but in particular by the practical realization of node communication and bandwidth per node. In our numerical examples we do not consider these details and instead only focus on the behavior of the scheme in \eqref{eq:CBOpar} in a synthetic setting.
\end{remark}
\begin{remark}
Within this section we will often employ $N$ nodes running CBO in parallel without interaction as a baseline comparison. 
This can be recovered by choosing $W=I$ and thus $\mathcal{N}_i = \{i\}$ in the above formulation. The initialization ensures $Z_0^i = \widehat{C}_0^i$ and by induction, if $Z_k^i = \widehat{C}_k^i$, then $W=I$ yields
\begin{align*}
Z_{k+1}^{i} = Z_k^i + (\widehat{C}_{k+1}^i  - \widehat{C}_{k}^i) = \widehat{C}_{k}^i + (\widehat{C}_{k+1}^i  - \widehat{C}_{k}^i) = 
\widehat{C}_{k+1}^i.
\end{align*}
\end{remark}

\subsubsection{Heuristics}\label{sec:heuristics}
\paragraph{Heuristic 1: Scheduling $\alpha$}

The choice of $\alpha$ can greatly impact the performance of the CBO algorithm. In order not to obfuscate the results obtained in this work, we adopt a simplistic strategy used in other works like \cite{bungert2025mirrorcbo,beddrich2026constrained}: At the start of the dynamics the particles should explore more, while towards the end they should be confined more strongly towards the current best particle. This intuition leads to the following simple scheduling $\alpha_k^i = \min\{\iota^k\cdot \alpha_0, \alpha_{\max}\}$, where throughout this work, if not stated otherwise, we choose $\iota=1.05, \alpha_0=10^{-4},\alpha_{\max}=10^9$.

\paragraph{Heuristic 2: Particle Clipping} As noted in \cite{beddrich2026constrained} it is favorable for CBO algorithms to enhance the exploration by choosing $\sigma$ larger than theoretic results allow. In the following experiments we will also employ this heuristic, which however risks that the ensemble diverges, if particles are moved too far by high noise. To this end, we employ another heuristic, that restricts the particles to a problem dependent, rectangular domain $\domain_f=[-a_f,a_f]^d$. Namely, in \textbf{PostStep}, for each node $i=1:N$ we then perform the following modification:
\begin{align*}
X^\ell_{k+1,s} \gets \min\{\max\{X^\ell_{k+1,s}, -a_f\}, a_f\},\qquad\text{ for all } \ell\in\mathcal{M}_i, s=1:d .
\end{align*}
This can be seen as infusing the problem with prior information about the function $f$, which however in our case will be only used as a safeguard, with very mild assumptions about the function, see \cref{tab:problem-thresholds}.

\paragraph{Heuristic 3: Resampling} 
In practical scenarios, CBO algorithms can suffer from an early collapse, which means that all particles are effectively at the same position, although they did not reach the global minimizer yet. In our multi-node setting this does not immediately yield to stationarity, but instead to interesting dynamics, we describe in \cref{sec:conform}. However, as in the standard CBO case, this situation is unfavorable from the optimization point of view, because once the particles per node collapse, the information from objective function evaluations is disregarded and convergence towards the global minimizer is essentially lost. In order to circumvent such situations, we adopt the collapse prevention strategy from \cite{carrillojinlizhu2021}: If at iteration $k$ a certain collapse criterion is met at node $i$, the routine in \textbf{PostStep} will add independent Gaussian noise of scale $\sigma^{\text{indep}}\geq 0$, to expand the ensemble again:
\begin{align*}
X^{\ell}_{k+1} \gets X^{\ell}_{k+1} + \sigma^{\text{indep}} \sqrt{\Delta t} \xi^\ell,\quad \xi^\ell\sim\mathcal{N}(0,I), \quad\text{for all}\quad \ell \in \mathcal{M}_i.
\end{align*}
We consider two different collapse criteria, the first being based on how much the ensemble per node changed during one iteration. Namely, this criterion is met at node $i$, if $\sum_{\ell\in\mathcal{M}_i} \|X^\ell_{k+1} - X^\ell_{k}\|^2 < 
\delta^{\text{thresh}}$ with some given parameter $\delta^{\text{thresh}}$. The other criterion we consider, checks whether the best found loss value decreased during the last $K\in\mathbb{N}$ iterations. Thus for $k>K$ this criterion is met if
\begin{align*}
\min_{t=(k-K):k ,\ell\in\mathcal{M}_i} f(X^\ell_t) \geq \min_{t=1:(k-K-1),\ell\in\mathcal{M}_i} f(X^\ell_t).
\end{align*}
We specify the concrete choice of these criteria below.

\subsubsection{Choice of mixing matrix}\label{sec:choice}

We present the different mixing matrices $W\in\R^{N\times N}$ employed in the following. Note, that this choice also determines the communication topology, through the sparsity pattern of $W$.

\paragraph{Full} Allowing all nodes to communicate with each other with the same strength, we obtain the following doubly stochastic mixing matrix

\vspace{.5em}\begin{minipage}{.6\textwidth}%
\begin{align*}
W^{\text{Full}} = 
\frac{1}{N}
\begin{pmatrix}
1 & 1 & \cdots & 1\\
1 & 1 & \cdots & 1\\
\vdots & \vdots & \ddots & \vdots\\
1 & 1 & \cdots & 1\\
\end{pmatrix}.
\end{align*}\end{minipage}%
\begin{minipage}{.35\textwidth}\centering%
\includegraphics[width=.4\textwidth]{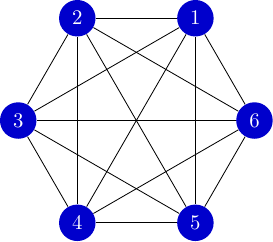}\\
Full topology for $N=6$.
\end{minipage}\vspace{.5em}

\paragraph{Ring topology} Figuratively placing the available nodes in a ring, we only allow them to communicate with their direct neighbours. We assume a uniform communication strength controlled by a parameter $\miparam\in[0,1]$ and thus the communication is modelled by the following mixing matrix:

\vspace{.5em}\begin{minipage}{.6\textwidth}%
\begin{align}\label{eq:matgo}
W_\miparam^{\text{ring}} = 
\begin{pmatrix}
\miparam               & \frac{1 - \miparam}{2} & 0                     & \ldots &0& \frac{1 - \miparam}{2}\\
\frac{1 - \miparam}{2} & \miparam               & \frac{1 - \miparam}{2} & 0 &\ldots & 0\\
&&\ddots & \ddots & \ddots \\
 \frac{1 - \miparam}{2} & 0 & \ldots & 0& \frac{1 - \miparam}{2} & \miparam
\end{pmatrix},
\end{align}
\end{minipage}%
\begin{minipage}{.35\textwidth}\centering%
\includegraphics[width=.4\textwidth]{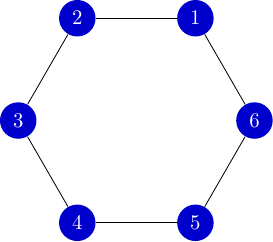}\\
Ring topology for $N=6$.
\end{minipage}\vspace{.5em}

\noindent\def\miiparam{\mu}%
which is indeed symmetric and doubly stochastic. In particular, we observe that $W_{1}^{\text{ring}} = I$, i.e., $\miparam=1$ recovers the case where $N$ nodes perform independent CBO runs. We also consider the case where the node can connect to neighbors that are more than one ``hop'' away, for example the following $2$-hop ring connects each node to $4$ neighbours, defining $\miiparam = \tfrac{1-\param}{4}$:

\vspace{.5em}\begin{minipage}{.6\textwidth}%
\begin{align*}
W_\miparam^{2\text{hop}} = 
\begin{pmatrix}
\miparam               & \miiparam & \miiparam & 0                     & \ldots &0& \miiparam& \miiparam\\
\miiparam & \miparam               & \miiparam & \miiparam & 0 &\ldots & 0 & \miiparam\\
&&&\ddots & \ddots & \ddots \\
 \miiparam & \miiparam& 0 & \ldots & 0& \miiparam &\miiparam & \miparam
\end{pmatrix}.
\end{align*}
\end{minipage}%
\begin{minipage}{.35\textwidth}\centering%
\includegraphics[width=.4\textwidth]{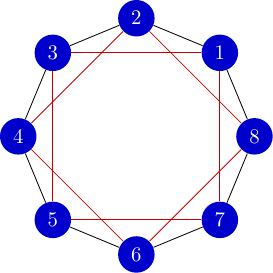}\\
Ring topology for $N=8$.
\end{minipage}\vspace{.5em}
\paragraph{Star topology}  We can assign a single node to be central in the sense, that all other nodes only communicate with this centre node. In the graph, this can be seen as a star topology, which is modelled as follows:

\vspace{.5em}\begin{minipage}{.61\textwidth}%
\begin{align*}
W^{\text{star}} = 
\begin{pmatrix}
\frac{1}{N}  & \frac{1}{N} & \frac{1}{N} & \frac{1}{N} &\cdots & \frac{1}{N}\\
\frac{1}{N} & 1 - \frac{1}{N} & 0 & 0 & \cdots & 0\\
\frac{1}{N} & 0               & 1 - \frac{1}{N}& 0 & \cdots & 0\\
& \ddots &  & \ddots\\
\frac{1}{N} & 0 & \cdots & & 0 & 1 - \frac{1}{N}
\end{pmatrix}.%
\end{align*}
\end{minipage}%
\begin{minipage}{.34\textwidth}\centering%
\includegraphics[width=.4\textwidth]{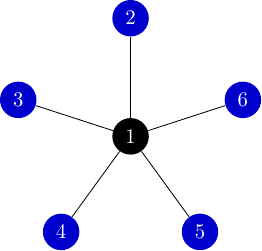}\\
Star topology for $N=6$.
\end{minipage}\vspace{.5em}
\subsection{Consensus formation}\label{sec:conformnum}
\begin{figure}
\centering
\includegraphics[width=\linewidth]{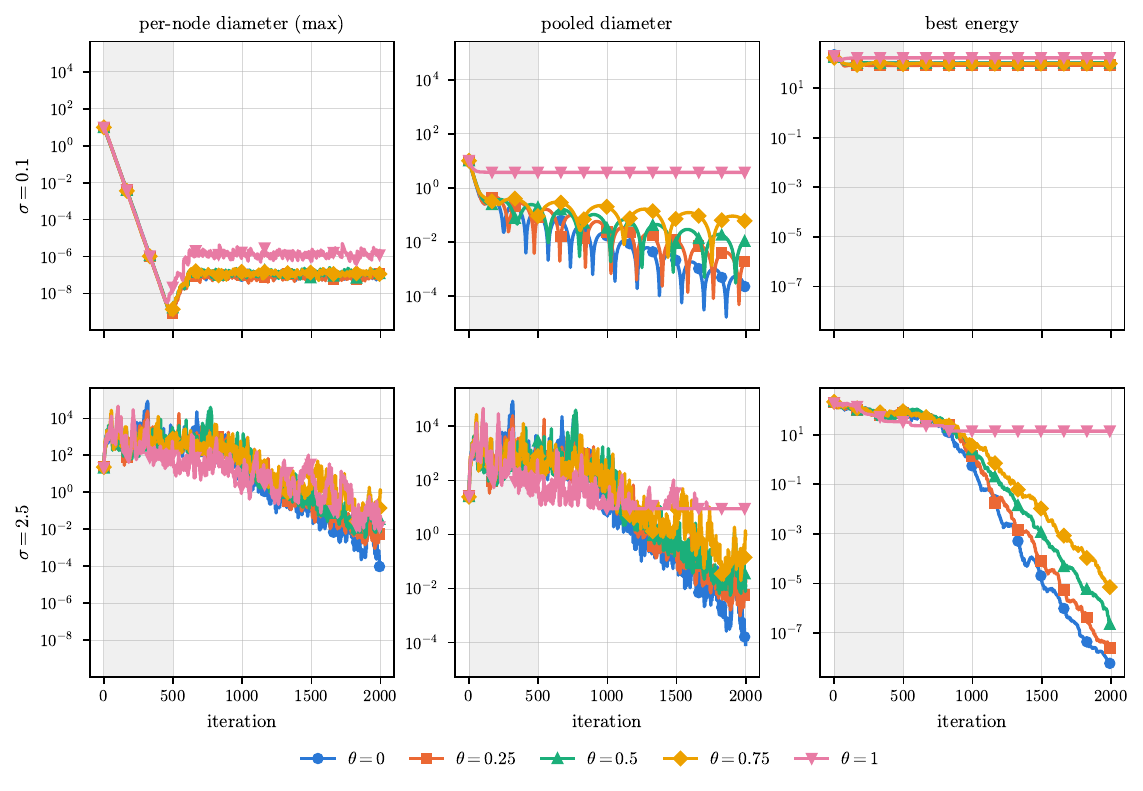}
\caption{Consensus formation described in \cref{sec:conformnum} underlining the theory in \cref{sec:conform}. We choose $N=50$ nodes with $M=20$ particles per node in the ring topology $W^{\text{ring}}_\miparam$ for different values of $\miparam\in[0,1]$.}
\label{fig:conform}
\end{figure}

As the first numerical study, we study the consensus formation in \eqref{eq:CBOpar}, which was theoretically analyzed in \cref{sec:conform}. We focus on the ring mixing matrix $W^{\text{ring}}_\miparam$, since the parameter $\miparam$ allows us to vary the strength of the interaction, with $\miparam=1$ corresponding to the case of no interaction. For the test problem we choose the Rastrigin function, see \cite{hansen2009real} for the definition, in dimension $d=17$ and let \cref{alg:CBO} run for $2000$ iterations with a step size of $\Delta t=0.05$. We display the results in \cref{fig:conform} for different values of $\miparam$ and noise parameters $\sigma$. The third column displays the evolution of the best current energy across all nodes, i.e.,
\begin{align*}
e_k = \min_{\ell=1:NM} f(X_k^\ell).
\end{align*}
The first column displays the diameter within each node studied in \cref{lemma:consensus_structure}. Namely we consider the diameter of each node separately and take the maximum value, across nodes $i=1:N$ and dimensions $s=1:d$, i.e., the tracked quantity is given as
\begin{align*}
d_k = \max_{s=1:d}\max_{i=1:N} \operatorname{diam}(\mathcal{X}^{s,i}_k).
\end{align*}
The second column then considers the diameter studied in \cref{lemma:global_consensus}, which now takes the diameter across nodes and we thus refer to it as \emph{pooled}, since particles from all nodes are pooled together. We then take the maximum value across the dimensions $s=1:d$ and thus track the following quantity
\begin{align*}
d_k^{\text{pooled}} = \max_{s= 1:d} \operatorname{diam}(\mathcal{X}^{s}_k).
\end{align*}
For the runs shown in \cref{fig:conform} we deliberately exclude heuristic 2 and 3, in order to study cleanly the effect of consensus formation. 
\paragraph{Small noise case} Focusing on the first row, we observe that with a small noise parameter $\sigma$, for all parameters $\miparam$ there is a quick collapse of the per-node diameter, which has a negative impact on the optimization performance. For the pooled diameter, we first observe that in the case $\miparam=1$ there is visible no-convergence behavior. This is expected since in this case the different nodes do not interact with each other and there is no additional force pulling different nodes together. One could expect convergence here, if each separate node would converge to the unique minimizer of $f$, which however is not the case. This is in line with \cref{lemma:global_consensus} since for $\miparam=1$ the corresponding graph $\mathcal{G}$ is not connected thus violating Assumption \ref{a:W}. For the case $\miparam<1$ we curiously observe oscillations, which are discussed in \cref{sec:ballet} and we refer to as the \emph{consensus ballet}. This is not the preferred situation from an optimization point of view, but is interesting form collective behavior perspective. While not being monotonic, the sequence $d_k^{\text{pooled}}$ still appears to tend to zero confirming \cref{lemma:global_consensus}.

\paragraph{Intermediate noise case} In the second row, we observe the numerical benefit of the algorithm when choosing appropriate hyperparameters. Due to the higher noise, the per-node collapse is slowed down and allows for the necessary exploration needed in CBO. The cases where the interaction is active, i.e., $\miparam<1$ greatly outperform the no-interaction case, which again fails to converge towards a good approximation in this setting. For $\miparam<1$ we also observe that $d_k^{\text{pooled}}$ tends to zero confirming \cref{lemma:global_consensus}. For $\miparam=1$ the stagnation of $d_k^{\text{pooled}}$ can be explained in the same way as in the small noise case. Interestingly, even though the higher exploration noise obfuscates the clear oscillations in the diameter, one can still observe them in the evolution of the energy. This hints at the fact, the consensus ballet described below occurs even for non-clustered nodes.
\iflong\subsection{Oscillatory behavior in the consensus formation}\label{sec:ballet}
\subsubsection{A connection to the graph Laplacian}
In our numerical study of the consensus formation in \cref{fig:conform} we observe oscillations appearing in the pooled diameter, especially in the case with small noise parameters $\sigma$. With low noise, the ensemble within each node quickly collapses and thus for illustrative purposes, we assume in the following, that each ensemble is already fully collapsed, i.e.,
\begin{align*}
P^i_k:= X^\ell_k = X^{\tilde{\ell}}_k = \widehat{C}^i_k,\qquad \text{ for all } \ell,\tilde{\ell}\in\mathcal{M}_i\text{ and all } i\in\{1,\ldots,N\}.
\end{align*}
Since in this case the noise term vanishes and the system is deterministic, \cref{alg:CBO}  reduces to a linear update with two coupled states $P_k :=(P^1_k,\ldots, P^N_k), Z_k:=(Z_k^1,\ldots, Z_k^N)$:
\begin{align}\label{eq:simpleODE}
\begin{aligned}
P_{k+1} &= P_k - \Delta t (P_k - Z_k),\\%
Z_{k+1} &= W Z_k + P_{k+1} - P_k = W Z_k -\Delta t (P_k - Z_k),
\end{aligned}
\end{align}
which can also be written more compactly in the state $U_k=(P_k,Z_k)$
\begin{align*}
U_{k+1} = B U_k,\qquad \text{ with }
B = 
\begin{pmatrix}
(1-\Delta t) I & \Delta t I\\
-\Delta t I    & W + \Delta t I
\end{pmatrix}.
\end{align*}
Let again $\gamma_1\geq\ldots\geq\gamma_N$ denote the eigenvalues of $W$, then the eigenvalues of the matrix $B$ are given as
\begin{align*}
\gamma^B_i = \frac{1+\gamma_i}{2} \pm \frac{1}{2}\sqrt{(1-\gamma_i)(1-\gamma_i - 4\Delta t)}.
\end{align*}
Thus we observe that the matrix $B$ has complex eigenvalues, if one $i=1:N$ satisfies 
\begin{align*}
\gamma_i < 1\quad\text{ and }\quad \frac{1-\gamma_i}{4} < \Delta t,
\end{align*}
which is also the case, in which we expect the solution of the simplified system in \eqref{eq:simpleODE} to exhibit oscillatory behavior. For the ring matrix in \eqref{eq:matgo} we have $\gamma_1=1$ and the second largest eigenvalue is given as
\begin{align*}
\gamma_2 = (1-\theta)\cos(\tfrac{2\pi}{N}) + \theta,
\end{align*}
thus taking $\Delta t$ less or equal than $\tfrac{1-\gamma_2}{4}$ would ensure that all eigenvalues of $B$ are real.

Alternatively, we can also try to interpret \eqref{eq:simpleODE} as an explicit discretization of an ODE, for which we write
\begin{align*}
U_{k+1} = U_k - \Delta t\, C U_k, \qquad\text{ with } \quad
C = 
\begin{pmatrix}
I & -I\\
I & \frac{1}{\Delta t} (I-W) -  I
\end{pmatrix}.
\end{align*}
If we want to consider the limit $\Delta t\to 0$ we need to ensure that the lower right entry of $C$ is well-defined, which is for example the case if $W= I - \Delta t L$ for some matrix $L$. In the following we argue, that this matrix $L$ can be seen as a graph Laplacian on the connectivity graph of the nodes. Since $W$ is a doubly stochastic matrix its weighted degree matrix (with self-loops being counted once) is given as $D=I$. Thus $L=\frac{1}{\Delta t}(I-W)$ is indeed the weighted graph Laplacian of the graph given by $\frac{1}{\Delta t} W$. Consider now the special case of the ring topology $W^{\text{ring}}_\theta$, and choose $\theta = 1 -\Delta t$, then we have that
\begin{align*}
L &= \frac{1}{\Delta t}(I - W^{\text{ring}}_\theta) \\&= 
\frac{1}{\Delta t}
\begin{pmatrix}
\Delta t               & \frac{-\Delta t}{2} & 0                     & \ldots &0& \frac{-\Delta t}{2}\\
\frac{-\Delta t}{2} & \Delta t               & \frac{-\Delta t}{2} & 0 &\ldots & 0\\
&&\ddots & \ddots & \ddots \\
\frac{-\Delta t}{2} & 0 & \ldots & 0& \frac{-\Delta t}{2} & \Delta t
\end{pmatrix}
=
I - 
\begin{pmatrix}
0 & \frac{1}{2} & 0                     & \ldots &0& \frac{1}{2}\\
\frac{1}{2} & 0               & \frac{1}{2} & 0 &\ldots & 0\\
&&\ddots & \ddots & \ddots \\
\frac{1}{2} & 0 & \ldots & 0& \frac{1}{2} & 0
\end{pmatrix}
,
\end{align*}
i.e., $L$ is the random-walk graph Laplacian of the undirected graph connecting each node to its neighbors. Sending ${\Delta t}\to 0$ in \eqref{eq:simpleODE} with $W^{\text{ring}}_{1- \Delta t}$, we obtain the ODE
\begin{align}\label{eq:ODE}
\dot{U} = 
\begin{pmatrix}
-I & I\\
-I & I-L
\end{pmatrix} U.
\end{align}
The stationary states of this system are given by $U=(P,Z)$ such that
\begin{align*}
P=Z,\quad\text{ and }\quad LZ = 0,
\end{align*}
i.e., states where the consensus points and tracking variables coincide and solve the graph Laplace equation in each dimension $s=1:d$. Since the graph associated to $L$ is connected, this can only be true, if 
\begin{align*}
Z^1 = Z^2 = \ldots = Z^{N},
\end{align*}
i.e., all nodes share the same location in $\R^d$. While the numerical experiments indeed confirm a convergence to a single point, the dynamics exhibit oscillations, which we illustrate in the next section.
\subsubsection{The consensus ballet}
\begin{figure}
\centering%
\begin{subfigure}{\textwidth}%
\centering%
\includegraphics[width=.5\textwidth]{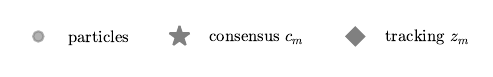}
\end{subfigure}\\
\begin{subfigure}{.42\textwidth}%
\includegraphics[width=\textwidth,trim={.5cm .9cm 0cm .5cm},clip]{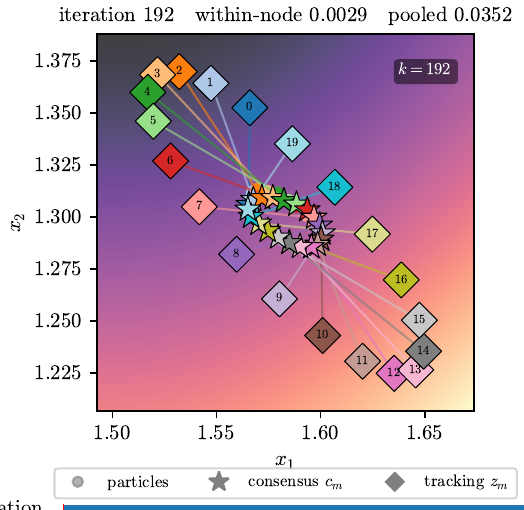}
\end{subfigure}\hfill%
\begin{subfigure}{.42\textwidth}%
\includegraphics[width=\textwidth,trim={.5cm .9cm 0cm .5cm},clip]{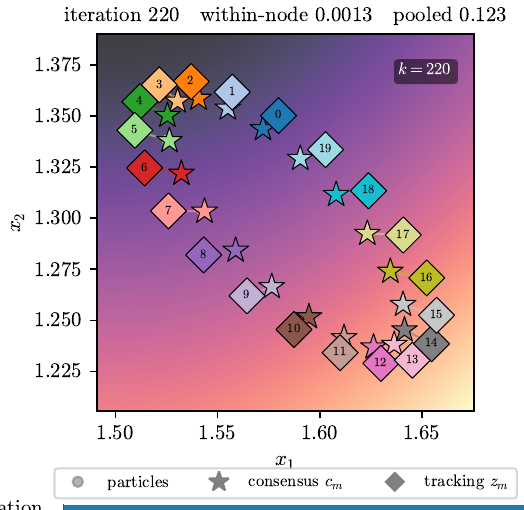}
\end{subfigure}\\
\begin{subfigure}{.42\textwidth}%
\includegraphics[width=\textwidth,trim={.5cm .9cm 0cm .5cm},clip]{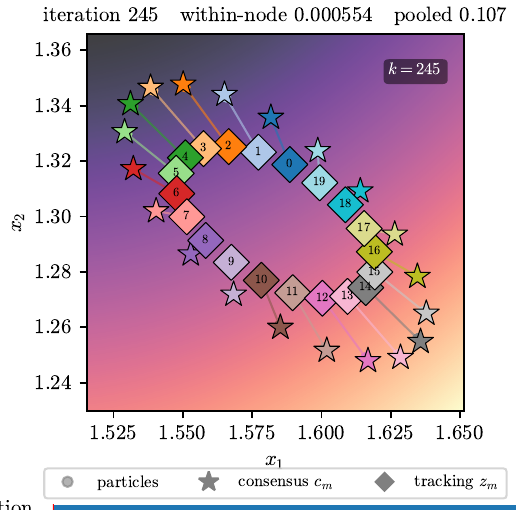}
\end{subfigure}\hfill%
\begin{subfigure}{.42\textwidth}%
\includegraphics[width=\textwidth,trim={.5cm .9cm 0cm .5cm},clip]{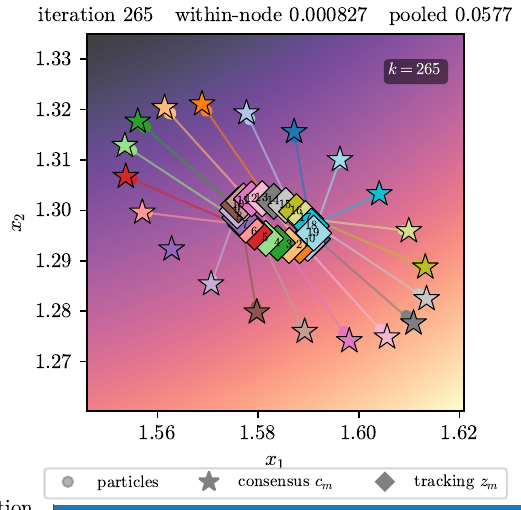}
\end{subfigure}\\
\begin{subfigure}[t]{.42\textwidth}%
\includegraphics[width=\textwidth,trim={.5cm .9cm 0cm .5cm},clip]{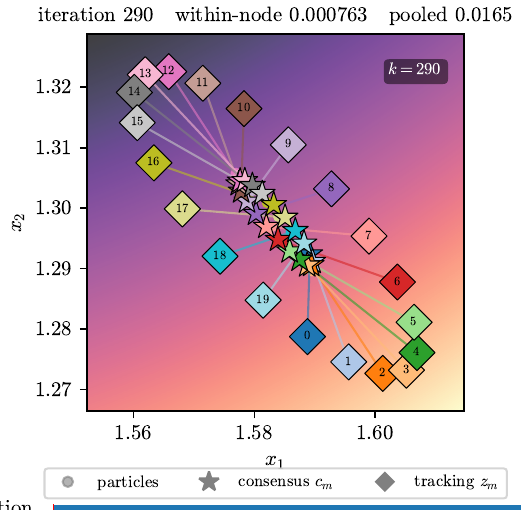}
\end{subfigure}\hfill%
\begin{subfigure}[t]{.42\textwidth}%
\includegraphics[width=\textwidth]{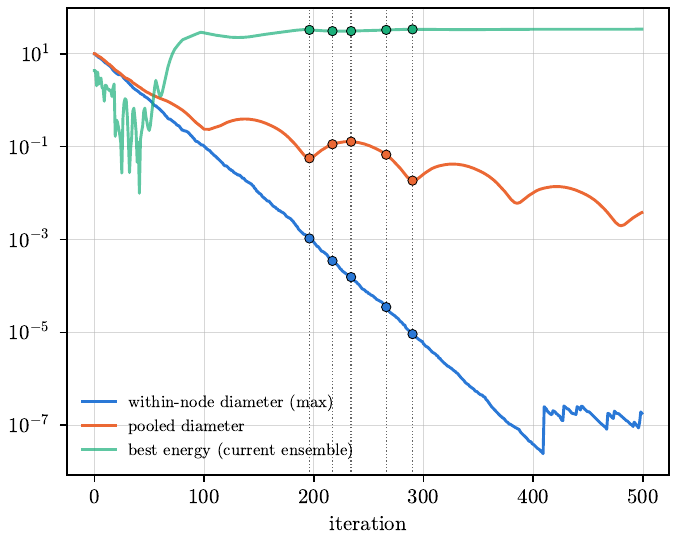}
\end{subfigure}
\caption{The consensus ballet. The underlying test function is the Rastrigin function, we choose $N=20$ nodes, $\Delta t=0.05,\sigma=0.01$ and employ $M=20$ particles per node.}\label{fig:ballet}
\end{figure}
In \cref{fig:ballet} we plot the evolution of the system with a small noise level $\sigma=0.01$ in dimension $d=2$ with the ring topology and $\theta=0.5$. As expected we can observe that per node the particles exhibit strong clustering towards their respective consensus point and are visually not distinguishable in the plot. The consensus points and tracking variables $Z^i_k$ arrange in two approximately concentric rings, with the outer one consisting of the tracking variables. We display the number of each node that is assigned to it in the graph ordering in \cref{fig:ballet} and observe that they are indeed ordered in the ring. We then observe that consensus points and tracking variables are pulled towards each other, but there is an overshooting effect, yielding each pair of consensus point and tracking variable to pass each other. This again yields two approximately concentric rings, but now the outer ring is the consensus points. Now the pairs are again pulled towards each other, finally resulting in a shrunken and flipped version of the initial state. The previously described process is then repeated cyclically while shrinking the diameter in each step. Due to rhythmic and coordinated behavior we refer to this as the consensus ballet.

We interpret this behavior using the ODE in \eqref{eq:ODE} and assume w.l.o.g. that $d=1$ since the dynamics are decoupled in the dimension. Denote by $v_1,\ldots,v_N\in\R^N$ the orthonormal basis of eigenvectors of the graph Laplacian $L$, which exists since $L$ is symmetric and by $\gamma_1^L\geq \ldots\geq \gamma_N^L$ the corresponding real eigenvalues. We then consider the linear spaces $V_i = \operatorname{span}(\{(v_i, 0), (0, v_i)\})\subset\R^{2N}$ and observe that $V_1+\ldots+V_N=\R^{2N}$. The specific structure of the system matrix in \eqref{eq:ODE} leaves each $V_i$ invariant under the dynamics, namely take $u = (\mu_1 v_i, \mu_2 v_i)^\top\in V_i$ with $\mu_1,\mu_2\in\R$, then
\begin{align*}
\begin{pmatrix}
-I & I\\
-I & I-L
\end{pmatrix}
u =
\begin{pmatrix}
(\mu_2-\mu_1) v_i\\
(\mu_2(1-\gamma^L_i) - \mu_1) v_i
\end{pmatrix}\in V_i.
\end{align*}
Thus when decomposing the solution $U$ of \eqref{eq:ODE} as follows
\begin{align*}
U(t) = \sum_{i=1}^N 
\binom{\mu_{1,i}(t) v_i}{\mu_{2,i}(t) v_i}
\qquad\text{with}\quad
\mu_{1,i}(t) = \langle U(t), (v_i,0)^\top\rangle,\, 
\mu_{2,i}(t) = \langle U(t), (0,v_i)^\top\rangle,
\end{align*}
plugging it in the ODE yields
\begin{align*}
\sum_{i=1}^N
\binom{\dot{\mu}_{1,i}(t) v_i}{ \dot{\mu}_{2,i}(t) v_i} = 
\sum_{i=1}^N 
\binom{(\mu_{2,i}(t) - \mu_{1,i}(t))v_i}{(\mu_{2,i}(t)(1-\gamma^L_i) - \mu_{1,i}(t)) v_i}.
\end{align*}
Thus it suffices to solve for each $i=1:N$ the system
\begin{align*}
\dot{\mu}_{1,i}(t) &= \mu_{2,i}(t) - \mu_{1,i}(t),\\
\dot{\mu}_{2,i}(t) &= 
\mu_{2,i}(t)(1-\gamma^L_i) - \mu_{1,i}(t)
\end{align*}
which expressed solely in $\mu_{1,i}$ is the ODE for the damped harmonic oscillator
\begin{align*}
\ddot{\mu}_{1,i}(t) + 
\gamma_i^L(\dot{\mu}_{1,i}(t)+
{\mu}_{1,i}(t))=0
\end{align*}
with characteristic polynomial
\begin{align*}
a^2 + \gamma_i^L a + \gamma_i^L = 0.
\end{align*}
Since $L$ is also positive semi-definite, we have $\gamma_i^L\geq 0$ for all $i=1:N$. As noted before, the null space of $L$ is one dimensional and thus $\gamma^L_N=0$ and $\gamma^L_i>0$ for $i=1:N-1$. Thus the equation for the $N$-th coefficient is solved for any $C_{1,N},C_{2,N}\in\R$ by
\begin{align*}
\mu_{1,N}(t) = C_{1,N} + C_{2,N} t.
\end{align*}
We furthermore have that $0<\gamma_i^L<4$ for all $i=1:N-1$ and thus the roots of the characteristic polynomial are complex valued in this case. Thus the coefficient ODE is solved for any $C_{1,i},C_{2,i}\in\R$ by
\begin{align*}
\mu_{1,i}(t) = 
e^{\nu_i t}\left(
C_{1,i} \cos(\omega_{i} t) + 
C_{2,i} \sin(\omega_{i} t)
\right),\qquad \nu_i = -\frac{\gamma_i^L}{2},\quad \omega_i = \frac{\sqrt{4\gamma_i^L-(\gamma_i^L)^2}}{2}.
\end{align*}
The ballet described in \cref{fig:ballet} is then the result of superimposing the above oscillating solitons across $i=1:N$ run independently in each dimension. The possible dynamics are much more diverse than the ring structure displayed in the figure and are determined by the initial state and number of nodes. For illustrative purposes we provide as a companion an interactive tool at \href{https://timroith.github.io/DecentralizedCBO/companion.html}{timroith.github.io/DecentralizedCBO/companion.html}, which visualizes the solutions to the ODE for different configurations.\else\fi%
\subsection{Hyperparameter study}

In this section we examine how the proposed algorithm is influenced by its hyperparameters $\Delta t$ and $\sigma$. We observe in \cref{fig:conform} that the choice of $\sigma$ has a significant influence on the numerical performance and thus we want to study it in more detail here. As the test problem, we choose the Rastrigin function in dimension $d=100$ and perform $5,000$ iterations for each scheme in the following. For the mixing matrix, we choose in each experiment the ring matrix $W^{\text{ring}}_\miparam$ with $\miparam$ being specified for each test.

\subsubsection{Interplay between noise strengths and number of particles}
\begin{figure}
\centering%
\begin{subfigure}{.33\textwidth}%
\includegraphics[width=\linewidth]{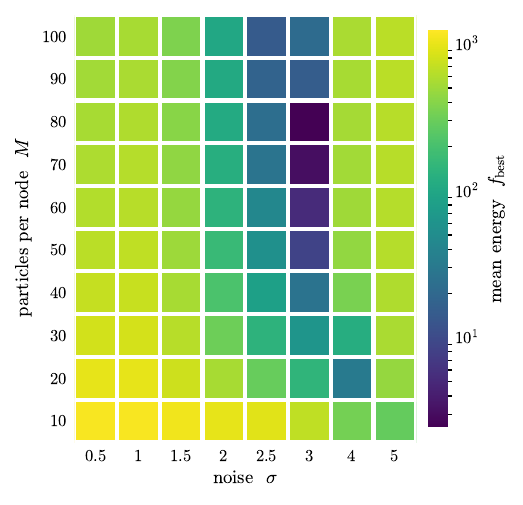}%
\caption{Ring parameter $\miparam=0$.}
\end{subfigure}%
\begin{subfigure}{.33\textwidth}%
\includegraphics[width=\linewidth]{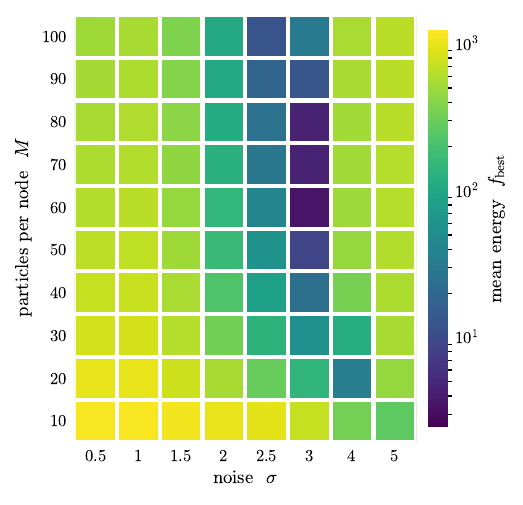}
\caption{Ring parameter $\miparam=\frac{1}{2}$.}
\end{subfigure}%
\begin{subfigure}{.33\textwidth}%
\includegraphics[width=\linewidth]{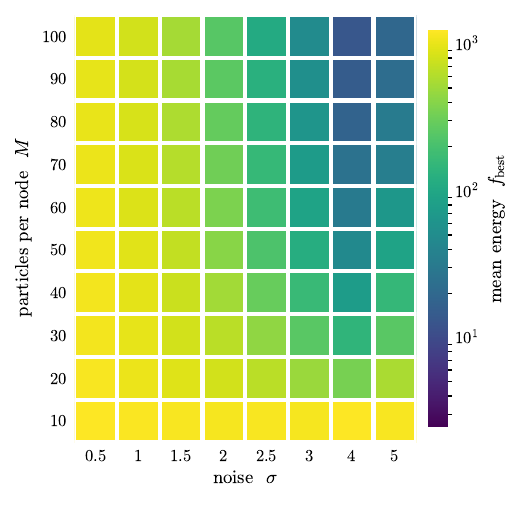}
\caption{Ring parameter $\miparam=1$.}
\end{subfigure}%
\caption{Interplay between the noise strength $\sigma$ and the number of particles $M$ for different parameters $\miparam\in[0,1]$ with the ring matrix $W^{\text{ring}}_\miparam$. The color value of each cell corresponds to the best objective value during a run of \cref{alg:CBO} with the specified hyperparameters, averaged over $20$ different seeds.}
\label{fig:sigmahyper}
\end{figure}
We first study the interplay between the noise parameter $\sigma$ and the number of particles per node $M$. In \cref{fig:sigmahyper} we display the numerical performance with the step size set to $\Delta t=0.02$ and the ring matrix $W_\miparam^{\text{ring}}$ with three different parameters $\miparam$. We observe, that node communication again leads to better results than $\miparam=1$. Beyond that, for all three cases there is a specific band visible in which the algorithms perform best. This band is sharper for $\miparam<1$, and no noticeable difference exists between the plots corresponding to two distinct values of $\miparam < 1$, suggesting that the decentralized variant is slightly more sensitive to the choice of hyperparameters.

\subsubsection{Interplay between step size and number of particles}
\begin{figure}
\centering%
\begin{subfigure}{.33\textwidth}%
\includegraphics[width=\linewidth]{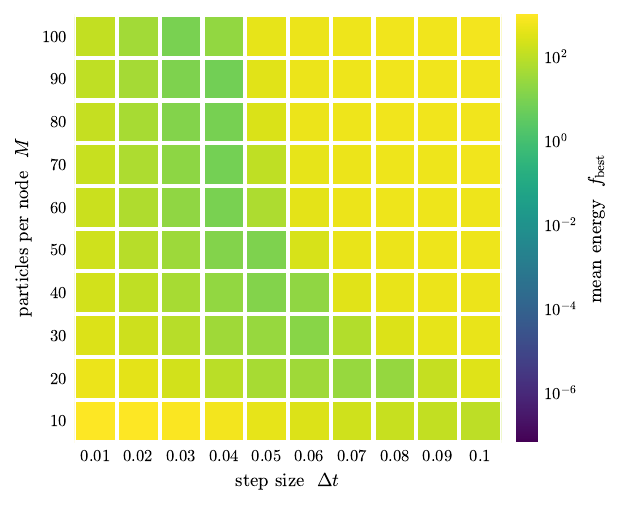}%
\caption{Number of nodes $N=4$}
\end{subfigure}%
\begin{subfigure}{.33\textwidth}%
\includegraphics[width=\linewidth]{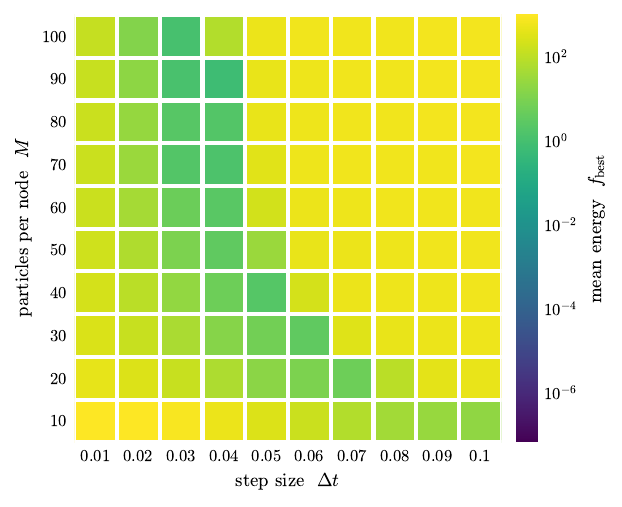}
\caption{Number of nodes $N=8$}
\end{subfigure}%
\begin{subfigure}{.33\textwidth}%
\includegraphics[width=\linewidth]{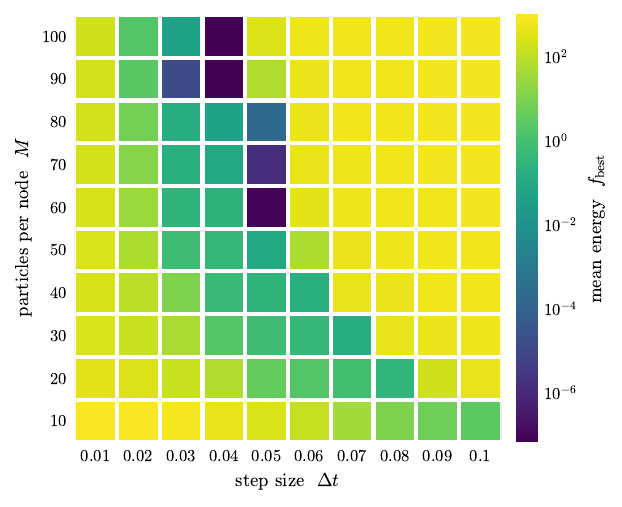}
\caption{Number of nodes $N=16$}
\end{subfigure}%
\caption{Interplay between the step size $\Delta t$ and the particles per node. The test problem is the Rastrigin function in $d=100$ with fixed noise parameter $\sigma=2.5$.}
\label{fig:dthyper}
\end{figure}
In our second study, we examine the interplay between the step size $\Delta t$ and the number of particles $M$ for a fixed noise parameter $\sigma=2.5$ and ring parameter $\miparam=0.5$. In \cref{fig:dthyper} we observe a phase transition that has a similar shape across different numbers of nodes. The effect that is visible in all three panels, is that with an increasing number of particles per node, the step size should decrease for favorable numerical performance. As a side effect of this plot, we also see, that increasing the number of nodes helps the numerical performance overall, a behavior we study separately in \cref{fig:scaledim}.

\subsection{Influence of the mixing matrix and frequency}

In this section we examine how the mixing matrix influences the performance of the decentralized CBO algorithm. Furthermore, we also ask the question, whether the communication cost can be reduced, by synchronizing the nodes only every $S$ steps in \cref{alg:CBO}.
\subsubsection{Comparison of different mixing topologies}
\begin{figure}
\centering
\includegraphics[width=\linewidth]{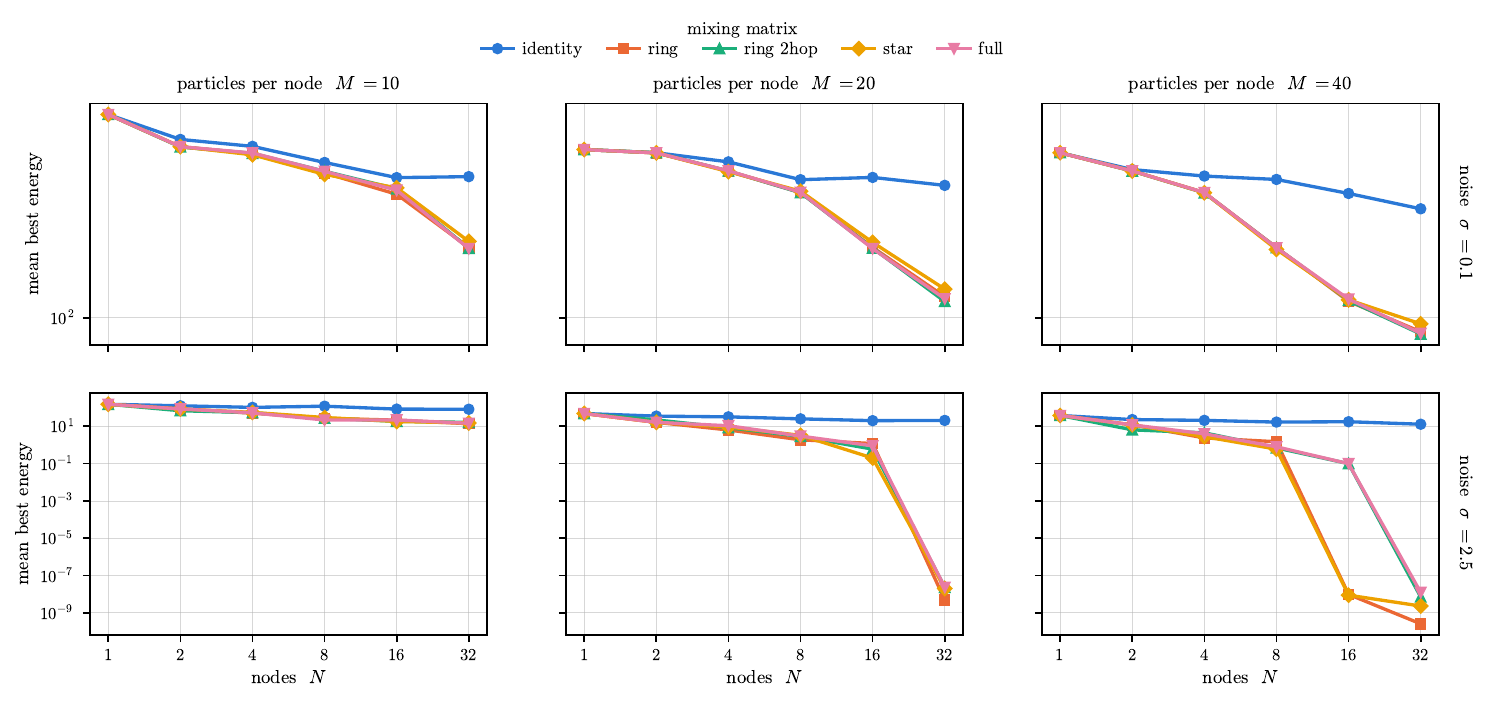}
\caption{Comparison between different mixing matrices. We observe that the connected topologies perform similarly and outperform the non-connected topology, i.e., $W=I$.}
\label{fig:difftopology}
\end{figure}

In our first test in \cref{fig:difftopology} we compare the performance of the described topologies in \cref{sec:choice} when varying the number of nodes, across different noise parameters $\sigma>0$. For all settings, we test on the Rastrigin function in dimension $d=20$, with $\Delta t=0.05$. The first notable result from \cref{fig:difftopology} is that different mixing matrices lead to very similar numerical performance. In particular there is no clear trend observable, which would make one of the topologies the preferable one. From a practical point of view, this suggests, that, whenever possible, working with the ring topology, which is favorable from a communication efficiency point of view, is indeed reasonable and seems not to lower the performance compared to other more costly choices.

\subsubsection{Variable Synchronization Frequency}\label{sec:mixfreq}
On a practical level, the synchronization step in \cref{alg:CBO} can be costly due to communications between the nodes. A natural question is, whether one can reduce the number of communications, by not synchronizing in every step of the CBO algorithm, but only every $S\in\mathbb{N}$ steps. In this section, we consider the behavior of the algorithm, when replacing lines 7--8 in \cref{alg:CBO} with the following modification:
\begin{algorithmic}
\For{$i=1:N$}
\If{$\operatorname{mod}(k, S) = 0$}
\State $Z_{k+1}^i = \sum_{j\in\mathcal{N}_i} W_{ij} Z_{k}^j + (\widehat{C}_{k+1}^i  - \widehat{C}_{k}^i)$
\Else
\State $Z_{k+1}^i = Z_{k}^i + (\widehat{C}_{k+1}^i  - \widehat{C}_{k}^i)$
\EndIf
\EndFor
\end{algorithmic}
\begin{figure}
\begin{subfigure}{.5\textwidth}
\includegraphics[width=\linewidth]{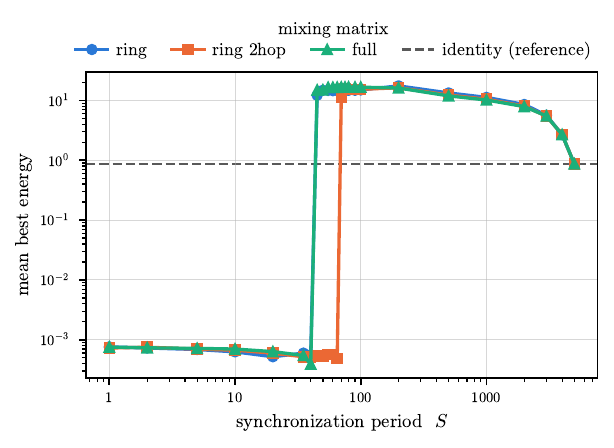}
\caption{Ackley function in $d=50$.}
\end{subfigure}%
\begin{subfigure}{.5\textwidth}%
\includegraphics[width=\linewidth]{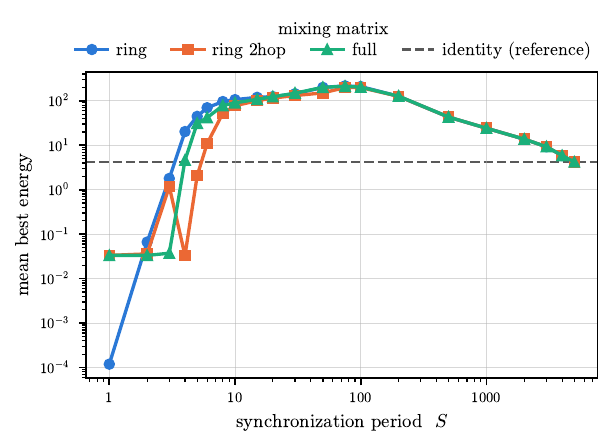}
\caption{Rastrigin function in $d=30$.}
\end{subfigure}
\caption{Numerical performance across longer synchronization periods described in \cref{sec:mixfreq}. Each point on the graph is averaged over $20$ independent runs.}
\label{fig:mixfreq}
\end{figure}
In \cref{fig:mixfreq} we plot the numerical performance, when increasing the length of the period $S$, thus decreasing the communication cost of the algorithm. We use the Rastrigin and Ackley functions as the underlying test problems and refer to \cite{hansen2009real,jamil2013literature} for their definitions. Note, that if we let the CBO algorithm run for $K\in\mathbb{N}$ iterations and per synchronization we have $C\in\mathbb{N}$ communications, the total number of communications is $\lfloor \tfrac{K}{S} \rfloor\cdot C$. Notably, in \cref{fig:mixfreq} we see that on Ackley for all topologies we can increase the period length to around $S=30$, without losing numerical performance. There is then a sharp interface, where increasing the period length even more, drastically decreases the performance. A possible explanation for this phenomenon is, that running the nodes independently for too long lets them diverge far apart, to a point, where even the averaging step cannot recover favorable dynamics. For Rastrigin, we can see a performance loss already for lower synchronization periods, also in this case we observe distinct behaviors across mixing topologies.

For both methods, the mean energies increase up to some worst configuration, after which the mean energies then drop again towards the reference value where the total number of iterations corresponds to the period length, i.e., $K=S$, which means that effectively there was no interaction at all. 
This result is particularly interesting from a practical point of view, because it tells us that in certain scenarios we can reduce the communication cost, by lowering the number of synchronizations without loosing much numerical performance. 
However, this has to be done with care, because at a certain point the performance will degrade, depending on the underlying problem.
\subsection{Behavior across problem dimensions}
\begin{figure}
\centering
\begin{subfigure}{0.5\linewidth}
\includegraphics[width=\linewidth]{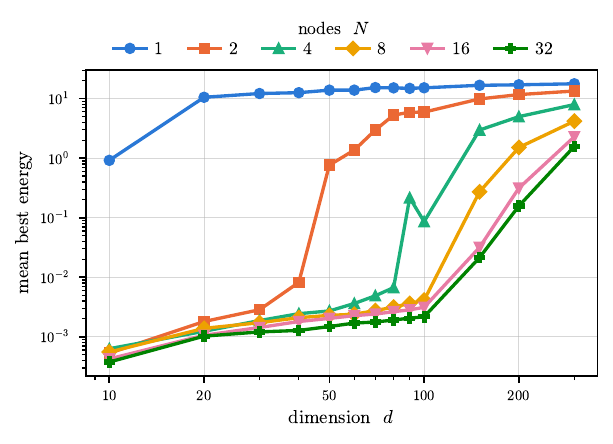}%
\caption{Ackley function.}
\end{subfigure}%
\begin{subfigure}{0.5\linewidth}
\includegraphics[width=\linewidth]{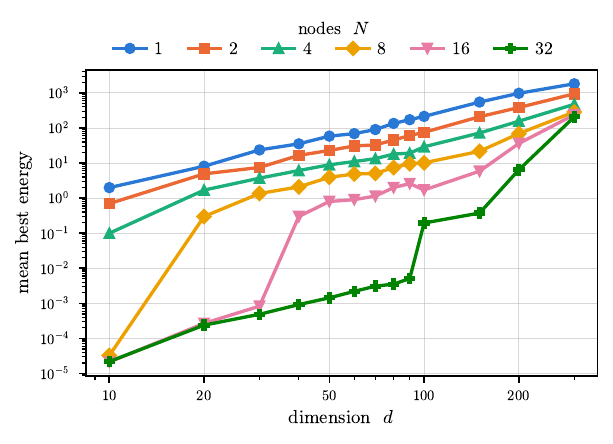}%
\caption{Rastrigin function.}
\end{subfigure}
\caption{Numerical performance with increasing problem dimension, across different number of nodes $N$.}
\label{fig:scaledim}
\end{figure}

In this section we test how the number of nodes can influence the numerical performance of the decentralized CBO algorithm across varying problem dimensions. To this end, we fix $\sigma =2.5, \Delta t = 0.05$, employ $M=20$ particles per node and use the Ackley and Rastrigin functions as the underlying test problems. We again choose the ring matrix $W^{\text{ring}}_\miparam$ with parameter $\miparam=\tfrac{1}{2}$ and the dynamics run for $5,000$ iterations. In \cref{fig:scaledim} we plot the results with the specified settings, averaged over $10$ different runs. In both cases, we observe that the difficulty of the problem increases with increasing dimension, where employing more nodes yields better numerical performance. However, solely increasing the number of nodes can significantly improve the performance only up to a certain number of nodes.

\subsection{Benchmarks}

In this section, we now perform a quantitative benchmark on different functions. We consider the Rastrigin, Ackley, Alpine, Styblinski--Tang, Lunacek Bi-Rastrigin and Katsuura functions and refer to \cite{hansen2009real,jamil2013literature} for their respective definitions. %
In order to prevent rewarding unwanted biases in the optimizers, we apply a random shift of the global optimum in every run. Denoting by $g$ the respective function, we randomly sample a new vector $a\in\R^d$ for every run, keep it fixed for this run, and then use $f(x)=g(x+a)$ as the objective function. Furthermore, we use heuristic 2, the particle clipping, with the domain boundaries detailed in \cref{tab:problem-thresholds}. Furthermore, we now also employ heuristic 3, the resampling scheme.

\begin{table}[t]
\centering
\resizebox{\textwidth}{!}{%
\begin{tabular}{lcccccc}
\toprule
 & Rastrigin & Ackley & Alpine & Styblinski--Tang & Lunacek & Katsuura \\
\midrule
Domain & $[-5.12, 5.12]^d$ & $[-32.768, 32.768]^d$ & $[-10, 10]^d$ & $[-5, 5]^d$ & $[-5, 5]^d$ & $[-5, 5]^d$\\
\bottomrule
\end{tabular}
}
\caption{Test problems and the domains they are typically optimized on. For our experiments, we initialize the particles uniformly within the displayed domain and set the clipping domain to a fixed factor of $4$ times the initialization domain.}
\label{tab:problem-thresholds}
\end{table}%
\subsubsection{Dimension $d=100$}
In \cref{tab:decben} we report a quantitative benchmark in dimension $d=100$. Therein, we count a run as successful if the gap to the known optimal value $f_{\min}$, is below a certain threshold, namely
\begin{align*}
f_{\text{best}} - f_{\min} < 0.0001,
\end{align*}
where in this case $f_{\text{best}}$ denotes the best objective value that the respective optimizers recorded during their runs. This interpretation changes in the case of noisy objective functions detailed further below.
In \cref{tab:decben} we compare the numerical performance of different CBO algorithms, where for each algorithm we allow a total of $20,000$ iterations. Therein, ``Multi-run CBO'' refers to running $N$ different instances of CBO in parallel without communication. Equivalently, this is decentralized CBO with the mixing matrix $W=I$. To compare against a different zero-order optimizer, we include results obtained with the covariance matrix adaptation evolution strategy (CMA-ES) \cite{hansen2006cma}, employing the \texttt{pycma} package \cite{nikolaus_hansen_2026_22738161}.

For the CBO runs with $M=20$ particles we choose a step size of $\Delta t = 0.02$ and, using the insights from \cref{fig:dthyper}, a slightly smaller $\Delta t = 0.01$ for $M=160$. Using the trend seen in \cref{fig:sigmahyper},  for non-decentralized CBO we choose $\sigma=3$ for $M=20$ and $\sigma=5$ for $M=160$. From the same plot we also observe that the decentralized version is more sensitive to choosing $\sigma$ too big and thus we keep $\sigma=3$ for both versions. Furthermore, for both decentralized versions we use the ring mixing matrix with parameter $\miparam=0.5$.

We make comparisons based on the following questions which \cref{tab:decben} gives answers to:
\begin{enumerate}[label=C\arabic*)]\itemsep=0pt%
\item \emph{Multi-run CBO vs Decentralized CBO}: Given $8$ nodes each with $M=20$ agents, does our proposed node communication scheme have a benefit, over running independent runs in parallel?
\item[$\to$] Up to the Katsuura function (which we comment on below) this can always be answered positively, which is in line with our previous observation that communication provides a benefit in most cases.
\item \emph{Serial CBO vs Decentralized CBO a)}: Decentralized CBO with $M=20$ particles on $N=8$ nodes implicitly evolves an ensemble of $N\times M=160$, which we assume distributed across different devices. Although this is not the main interest of the present paper one might ask how a serial CBO algorithm would perform given $160$ particles, such that the number of objective function evaluations is matched.
\item[$\to$] The serial variant with $M=160$ particles performs better than the decentralized version with only $M=20$. This result is expected since in this case the particles can directly communicate.
\item \emph{Serial CBO vs Decentralized CBO b)}: Assume now each computational node allows for a maximum number of $M=160$ particles to be handled. How much can the usage of $N=8$ nodes improve upon a single one.
\item[$\to$] Using more nodes improves the performance on the Rastrigin function significantly. This is the interesting case for this table since this is the test where the serial variant does not perform well.
\end{enumerate}
\begin{remark}
The Katsuura function defined in \cite{hansen2009real} has an extremely rugged landscape with $10^d$ local minimizers within $[-5,5]^d$. We first note, that the good performance of the serial CBO variants is surprising and we were not able to locate a similar performance in the literature. In particular, it appears that in this case the averaging across tracking variables has a negative impact on the numerical performance.
\end{remark}

\begin{table}
\centering
\resizebox{\textwidth}{!}{%
\begin{tabular}{lcccccc}
\toprule
Method & Rastrigin & Ackley & Alpine & Styblinski--Tang & Lunacek & Katsuura \\
\midrule
\multicolumn{7}{l}{\textbf{Success rate}} \\
Multi-run CBO ($N=8$, $M=20$) & 0\% & 0\% & 0\% & 0\% & 0\% & 0\% \\
Serial CBO    ($N=1$, $M=160$) & 25\% & \textbf{100\%} & \textbf{90\%} & \textbf{100\%} & 55\% & \textbf{85\%} \\
\midrule%
Decentralized CBO ($N=8$, $M=20$) & 0\% & 0\% & 0\% & 0\% & 0\% & 0\% \\
Decentralized CBO ($N=8$, $M=160$) & \textbf{95\%} & 95\% & 50\% & 90\% & \textbf{95\%} & 0\% \\
\midrule%
CMA-ES (popsize 160, best of 8 restarts) & 0\% & \textbf{100\%} & \textbf{100\%} & 0\% & 0\% & 0\% \\
\midrule\midrule%
\multicolumn{7}{l}{\textbf{Median gap $f_{\mathrm{best}} - f_{\min}$}} \\
Multi-run CBO ($N=8$, $M=20$) & 1.4e+02 & 14 & 4 & 7.6e+02 & 1.4e+02 & 0.0011 \\
Serial CBO ($N=1$, $M=160$) & 0.99 & 3.5e-08 & 3.2e-08 & 1.9e-05 & 2.2e-06 & \textbf{5.1e-05} \\
\midrule%
Decentralized CBO, ring ($N=8$, $M=20$) & 1.4 & 0.026 & 0.12 & 0.067 & 1.7 & 3.8 \\
Decentralized CBO, ring ($N=8$, $M=160$) & \textbf{7.9e-08} & 1.9e-07 & 0.00065 & \textbf{8.4e-06} & \textbf{7.6e-08} & 3.4 \\
\midrule%
CMA-ES (popsize 160, best of 8 restarts) & 17 & \textbf{1.9e-11} & \textbf{3.6e-11} & 14 & 20 & 3.5 \\
\bottomrule
\end{tabular}%
}
\caption{Performance of different algorithms over 20 runs in dimension $d=100$. We report both the success rates and the median gap of the objective values to the optimal one.}
\label{tab:decben}
\end{table}

\subsubsection{Dimension $d=1000$}
Beyond the previous scenario we test the optimizers for the high-dimensional case $d=1000$. In this case, we do not expect to find the global minimizer up to high precision, but we can still test how good the optimization performance is. In \cref{tab:benhighd} we consider three benchmark functions and report only the median gap of the best found objective value to the optimal one. Note, that in this dimension the standard CMA-ES is not feasible because of the computational effort of building the full covariance matrix. Thus in this case, we instead employ the separable variant introduced in \cite{ros2008simple} as sep-CMA-ES, which uses a diagonal covariance matrix instead. Also note, that we allow $8$ restarts to have a fair comparison to the decentralized version with $8$ nodes.

In \cref{tab:benhighd} we can clearly see, that the decentralized CBO variant in this high-dimensional scenario outperforms the serial CBO variant and also the CMA-ES on the selected test problems.

\begin{table}[t]
\centering
\resizebox{\textwidth}{!}{%
\begin{tabular}{lccc}
\toprule
Method & Rastrigin & Styblinski--Tang & Lunacek \\
\midrule
Serial CBO ($N=1$, $M=160$) & 1.2e+03 & 3.7e+03 & 9.8e+02 \\
Decentralized CBO, ring ($N=8$, $M=160$) & \textbf{93} & \textbf{0.0024} & \textbf{97} \\
\midrule%
sep-CMA-ES (popsize 160, best of 8 restarts) & 9.6e+02 & 1.2e+03 & 9.7e+02 \\
\bottomrule
\end{tabular}%
}
\caption{Benchmark in dimension $d=1000$, we report the median gap $f_{\text{best}} - f_{\min}$, where $f_{\text{best}}$ is the best found objective value in one run of the respective algorithm and $f_{{\min}}$ is the optimal value of the objective. The median is taken over 10 runs.}
\label{tab:benhighd}
\end{table}

\subsection{Noisy objective functions}

We now study the case, in which we only have access to an estimator of the objective function. To simulate this numerically, we add noise to the output, which we model by a stochastic function
\begin{align*}
\widehat{f}(x, \omega) := f(x) + \Xi(\omega), 
\end{align*}
where $\Xi:\Omega\to\R$ is a random variable on the probability space $(\Omega, \mathcal{A}, \mathbb{P})$. Note, that compared to \cref{tab:decben}, we do not formulate a success criterion, since the additive noise makes a precise threshold under which the objective evaluations should fall less meaningful. Instead, we only report the median gap to the optimal value. In this case, also the value of $f_{\text{best}}$ is not the best recorded value, since this is now biased by noise. Instead, we now take the output of each algorithm, denoted by $\bar{x}$ and then report $f_{\text{best}}=f(\bar{x})$. For CBO, the question then is, if this returned point $\bar{x}$ should correspond to the consensus $\widehat{C}_k$ at the final iterate $k$, the particle $\tilde{x}$ which resulted in the lowest value in $\widehat{f}$ for one specific noise realization, or for the decentralized version one of the tracking variables $Z^{i}_k$ for $i=1:N$. In order to decide this, we take the following approach after the algorithm terminates: We evaluate $\widehat{f}(\widehat{C}_k, \omega)$ for CBO, and 
$\widehat{f}(Z^{i}_k, \omega)$ for the decentralized variant, and additionally $\widehat{f}(\tilde{x}, \omega)$ for both, on five different realizations $\omega$ each and then return the point which resulted in the best median noisy objective value. The reported numbers in the tables are then evaluated in the original objective $f$ without noise, a quantity the algorithms do not have access to.

\subsubsection{Gaussian Noise}\label{sec:Gauss}

In \cref{tab:decben_noisy} we report the performance when $\operatorname{Law}(\Xi) = \mathcal{N}(0, \sigma_f^2)$, where in our case, we choose $\sigma_f=0.1$ for all functions. Additionally to the standard CMA-ES, here we also compare to a modified version, UH-CMA-ES \cite{hansen2008method}, which is specifically designed to handle noisy objective functions. Compared to \cref{tab:decben}, we observe that the decentralized variants are more robust to noise in the objective function than the serial CBO variants. 

\begin{table}
\centering
\resizebox{\textwidth}{!}{%
\begin{tabular}{lcccccc}
\toprule
Method & Rastrigin & Ackley & Alpine & Styblinski--Tang & Lunacek & Katsuura \\
\midrule
Serial CBO ($N=1$, $M=160$)   & 4.6 & 5.3 & 2.6 & 3.8 & 3.6 & 1.3 \\
Multi-run CBO ($N=8$, $M=20$) & 56 & 4.2 & 3.4 & 1.1e+02 & 67 & \textbf{0.59} \\
\midrule
Decentralized CBO, ring ($N=8$, $M=20$) & \textbf{0.13} & 0.31 & 0.22 & 0.13 & 0.14 & 3.7 \\
Decentralized CBO, ring ($N=8$, $M=160$) & 0.09 & 0.23 & 0.18 & \textbf{0.094} & \textbf{0.053} & 1.6 \\
\midrule
CMA-ES (popsize 160) & 32 & 0.071 & \textbf{0.1} & 92 & 32 & 4.3 \\
UH-CMA-ES (popsize 160) & 28 & \textbf{0.055} & 0.22 & 85 & 31 & 4.4 \\
\bottomrule
\end{tabular}%
}
\caption{Performance of different algorithms over $20$ runs with Gaussian noise added to the objective function. The reported value is the median objective value gap $f_{\text{best}} - f_{\min}$ where $f_{\text{best}}$ is calculated as described in \cref{sec:Gauss}.}
\label{tab:decben_noisy}
\end{table}

\subsubsection{Heavy-tailed noise}
We also perform the previous experiment in a setting where the distribution of $\Xi$ has heavier tails. Namely, we choose $\Xi = s_f \cdot \tilde{\Xi}$, where $\operatorname{Law}(\tilde{\Xi})$ is the Student's $t$ distribution with $2$ degrees of freedom and we choose the scale $s_f=0.1$ for all problems. The results are reported in \cref{tab:decben_noisystudent}, where we can again observe the noise robustness of the decentralized variants.
\begin{table}
\centering
\resizebox{\textwidth}{!}{%
\begin{tabular}{lcccccc}
\toprule
Method & Rastrigin & Ackley & Alpine & Styblinski--Tang & Lunacek & Katsuura \\
\midrule
Serial CBO ($N=1$, $M=160$) & 9 & 22 & 15 & 6.7 & 11 & 7.4 \\
Multi-run CBO ($N=8$, $M=20$) & 49 & 22 & 5.9 & 96 & 49 & \textbf{2.1} \\
\midrule%
Decentralized CBO, ring ($N=8$, $M=20$) & \textbf{0.26} & 1.9 & 0.34 & \textbf{0.25} & \textbf{0.28} & 8.6 \\
Decentralized CBO, ring ($N=8$, $M=160$) & 0.35 & 3.3 & 0.47 & 0.41 & 0.39 & 8.6 \\
\midrule%
CMA-ES (popsize 160) & 30 & \textbf{0.12} & \textbf{0.14} & 85 & 29 & 10 \\
UH-CMA-ES (popsize 160) & 31 & 21 & 0.22 & 85 & 31 & 9.2 \\
\bottomrule
\end{tabular}%
}
\caption{The same experiment as in \cref{tab:decben_noisy} with heavy-tailed instead of Gaussian noise added to $f$.}
\label{tab:decben_noisystudent}
\end{table}
\subsubsection{Machine learning example}
\def\params{x}
We now consider an application from machine learning, where we are given a finite labeled dataset $\mathcal{T}$ and our task is to find a parametrized classifier $h^\params{}$ with parameters $\params{}\in\R^d$ that solves the associated classification problem. We choose for $h^\params{}$ a simple neural network as in \cite{bungert2025mirrorcbo}, which is defined with weight and bias $\tilde{W}\in\R^{784\times 10}, b\in\R^{10}$ as follows
\begin{align*}
h^\params{}(u):= 
\operatorname{SoftMax}(\operatorname{BN}(\operatorname{ReLU}(\tilde{W}u + b))).
\end{align*}
For a vector $z\in\R^{\tilde{d}}$ we set $\operatorname{ReLU}(z)_i := \max\{0, z_i\}$,
$\operatorname{SoftMax}(z)_i := \frac{\exp((z)_i)}{\sum_{j=1}^{\tilde{d}} \exp((z)_j)}$ and the batchnorm $\operatorname{BN}(z)^{(b)}_i := 
\frac{(z^{(b)})_i - m(z)_i}{\sqrt{s(z)_i^2 + \epsilon}} \, (\gamma_{\text{BN}})_i + (\beta_{\text{BN}})_i$ acting on a batch of inputs $(z^{(1)}, \ldots, z^{(B)})$ with batch mean and standard deviation $m(z), s(z)$, zero-division prevention $\epsilon>0$ and trainable parameters $\gamma_{\text{BN}}, \beta_{\text{BN}}$. This results in a total of $d=7840 + 3 \cdot 10 = 7870$ learnable parameters, which we collect in the parameter variable $\params{}=(\tilde{W}, b, \gamma_{\text{BN}}, \beta_{\text{BN}})$. The supervised learning task is then the following optimization problem
\begin{align*}
\min_{\params{}\in\R^d} \frac{1}{\#\mathcal{T}}
\underbrace{\sum_{(u,v)\in \mathcal{T}} \ell(h^\params{}(u), v)}_{=:f(\params{})},
\end{align*}
\begin{table}[t]
\centering
\resizebox{\textwidth}{!}{%
\begin{tabular}{lcccc}
\toprule
 & \multicolumn{2}{c}{Cross-entropy} & \multicolumn{1}{c}{0-1 loss (not differentiable)} & \multicolumn{1}{c}{Disjoint data per node} \\
\cmidrule(lr){2-3}\cmidrule(lr){4-4}\cmidrule(lr){5-5}
Method & B=64 & B=512 & B=512 & B=64, federated \\
\midrule
Serial CBO ($N=1$, $M=50$) & 85.4\% & 89.1\% & 83.4\% & --- \\
Decentralized CBO, ring ($N=20$, $M=50$) & 90.4\% & 91.8\% & \textbf{91.5\%} & \textbf{90.2\%} \\
SGD (1000 restarts, best) & \textbf{92.9\%} & \textbf{92.9\%} & --- & ---\\
\bottomrule
\end{tabular}%
}
\caption{Results for the machine learning application on the MNIST dataset \cite{lecun2010mnist}. The reported results are the best accuracy on the test data set, taken over $5$ different repeats of each experiment. In this case we modify the hyperparameters for CBO to $\Delta t= 0.4,\sigma=0.1,\alpha_0=1$ and run for $37,500$ iterations. Furthermore, we compare the results to a stochastic gradient descent scheme (SGD) using the same mini-batch loss and run for $10,000$ iterations.}
\label{tab:mnist}
\end{table}%%
where $\ell$ is a loss function acting on the label space. In our experiments, we consider the differentiable cross-entropy $\ell^{\text{ce}}(v,\tilde{v}) = -\sum_{j=1}^C \tilde{v}_j \log(v_j)$ and the non-differentiable $0\text{-}1$ loss given as
\begin{align*}
\ell^{0\text{-}1}(v,\tilde{v}) = \begin{cases}
0 &\text{if } \argmax_{j=1,\ldots,C} \tilde{v}_j = \argmax_{j=1,\ldots,C} v_j,\\
1 &\text{else}.
\end{cases}
\end{align*}
With every query to the function $f$, we do not consider the sum over all training pairs, but only randomly subsampled mini-batches. This means that the stochastic function we evaluate is given as
\begin{align*}
\widehat{f}(\params{}, \omega) = \frac{1}{B}\sum_{(u,v)\in S(\omega)} \ell(h^\params{}(u), v),
\end{align*}
where $S(\omega)$ is a randomly sampled subset $S(\omega)\subset\mathcal{T}$ with maximum size $\#S(\omega) \leq B\in\N$. 

In the first column of \cref{tab:mnist} we compare the performance of the decentralized CBO algorithm to the single-node version and SGD, where we consider classifiers trained on the MNIST dataset \cite{lecun2010mnist} and report the accuracy on the test set, which was not part of the training process. For $B=64$ we observe a significant performance improvement when using multiple nodes in CBO. A possible explanation for this behavior is not only the increased total number of particles of decentralized CBO, but also that each individual node sees a different batch of the data, meaning that per iteration the multi-node version sees $20$ (not necessarily distinct) mini-batches of data, whereas the single node version only $1$. In order to study this effect, we increase the batch size to $B=512$ with both methods and see that indeed the serial version improves a lot. However, the performance of the decentralized version is still better. Moreover, it should be noted that for many problems one cannot arbitrarily increase $B$ since the associated batch of data would not fit into memory. Thus this is a real advantage of the decentralized version, since it preserves the favorable performance, even with a smaller batch size, split across different nodes.

As expected, in \cref{tab:mnist} we observe that stochastic gradient descent (SGD) with many restarts still beats the CBO variants, which use only zero-order information. We thus furthermore consider the $0\text{-}1$ loss, which is not differentiable and cannot be handled directly with SGD. Here, the decentralized version is more robust to this more difficult scenario than the single-node version. 

Beyond that the decentralized approach opens up another interesting possibility, namely federating the data as in \cite{carrillogarciatrillos2024fedcbo}. For this we now split up the dataset randomly (i.i.d.) into disjoint subsets $\mathcal{T} = \mathcal{T}_1 \cup \ldots \cup \mathcal{T}_N$, where node $i$ only has access to its own version of the data $\mathcal{T}_i$. In \cref{tab:mnist} we can see that there is a slight drop in performance, however the decentralized CBO algorithm still performs surprisingly well even in this federated setting.
\section{Concluding Remarks}
We considered the global minimization of functions that can only be evaluated through a noisy, possibly biased, oracle, and we proposed a decentralized parallel variant of CBO in which each computational node evolves its own swarm towards a local tracking variable, updated through a single round of communication, without the need for a central coordinator. Under suitable assumptions on the method's parameters and the underlying communication network, we proved that the agents achieve almost sure consensus. Moreover, under a hitting-probability assumption verified for a simplified dynamics, we proved convergence in expectation to the global minimizer within an explicitly bounded number of iterations. 

The numerical studies confirmed the theoretical findings on the consensus formation. Moreover, different hyperparameter and ablation studies established that the usage of multiple computational agents together with our proposed communication scheme has a beneficial impact on the numerical performance in most cases. However, we observed that the algorithm is sensitive to the correct choice of hyperparameters. If chosen correctly, we established a strong optimization performance on various benchmarks up to dimension $d=1000$ and stochastic objective functions. We also obtained good performance on a machine learning task, where in particular the non-differentiable $0\text{-}1$ loss, which is not possible to optimize with SGD, seems promising. The positive results in the case where data is split across different nodes call for a more in depth numerical study of possible applications, beyond the toy examples provided here.

\section{AI Methodology}
Following the template available at \url{https://simons.berkeley.edu/media/31229/download}, Table \ref{AI_tab} and the paragraph below describe how AI was involved in the research that led to this work and the preparation of this manuscript.
\begin{table}[H]
\begin{center}
\renewcommand{\arraystretch}{1.3}
\begin{tabular}{@{}l@{\qquad}c@{}}
\toprule
\textbf{Activity} & \textbf{AI used} \\
\midrule
Asking the research question & NO \\
Coming up with the approach & NO \\
Proof development & NO \\
Writing and exposition & NO \\
Checking for bugs & YES \\
Other supporting tasks & YES \\
\bottomrule
\end{tabular}
\end{center}
\caption{AI Methodology}\label{AI_tab}
\end{table}
% 1. Asking the research question: proposing a model, a definition, an assumption or a conjecture.
% 2. Coming up with the approach: the organizing plan of the proof, such as a reduction, a relaxation or a sequence of intermediate claims.
% 3. Proof development: the important technical ideas, those one would highlight in a 50-minute talk on the work.
% 4. Writing and exposition: substantive explanation and presentation, such as a proof overview or an example.
% 5. Checking for bugs: finding or correcting consequential errors, such as a circular argument or a missing hypothesis.
% 6. Other supporting tasks: copy-editing, figures, numerical experiments, minor lemmas, literature search.
% For the first four categories, answer Yes if a human making a comparable contribution would have been offered co-authorship; for the last two, any substantial contribution counts.

\begin{sloppypar}
\paragraph{Details}
The research direction, the method and the approach to the theoretical analysis, including proofs, as well as the core numerical implementation of the proposed method were developed by the authors without the employment of any AI tool. Claude Opus 5.5 and Sonnet 5.5 were prompted to review the draft to check the correctness of the mathematical arguments, aid uniforming the notation and copy-editing the text. Claude Opus 5.5 was used to debug the core implementation, prepare the benchmarks, executables for hyperparameter studies and with preparing the plots for the paper. It was furthermore used to prepare the illustrations for the oscillatory behavior in \cref{sec:ballet}.

% [Replace this bracketed guidance with a free-form story of how the work developed and how AI was involved. One might describe what was asked or tried, what was useful or unsuccessful, what the authors contributed or changed, and how the resulting arguments, text, or code were checked. Name the models or tools used, including versions where known. Focus on the parts relevant to the work; there is no need to discuss every checklist item separately.]

\end{sloppypar}

\bigskip
\begin{center}
  \FundingLogos
  
  \vspace{0.5em}
  \begin{tcolorbox}\centering\small
    % CHOOSE the sentence that applies to your grant/programme and replace placeholders.
    % 1) For Horizon Europe / ERC (recommended wording):
   
    Funded by the European Union. Views and opinions expressed are however those of the author(s) only and do not necessarily reflect those of the European Union or the European Research Council Executive Agency. Neither the European Union nor the granting authority can be held responsible for them. This project has received funding from the European Research Council (ERC) under the European Union’s Horizon Europe research and innovation programme (grant agreement No. 101198055, project acronym NEITALG).
    
    % 2) (Alternatively, for Horizon 2020 ERC grants use the H2020 wording:)
    % This project has received funding from the European Research Council (ERC) under the European Union's Horizon 2020 research and innovation programme (grant agreement No. <GRANT NUMBER>, project acronym <ACRONYM>).
  \end{tcolorbox}
\end{center}

\section*{Acknowledgments}
M.F. and T.R. acknowledge the support of the Munich Center for Machine Learning and the ERC Advanced Grant NEITALG, grant agreement No. 101198055. 

G.M. acknowledges the INdAM - GNCS project ``Scalable Algorithms for Machine Learning: Dimensionality Reduction, Stochastic and Distributed Optimization'', CUP E53C25002010001.

The authors gratefully acknowledge the computational and data resources provided by the Leibniz Supercomputing Centre (\url{www.lrz.de}).

\bibliographystyle{plain}
\bibliography{references}

@book{matan, place={Cambridge}, title={Matrix Analysis}, publisher={Cambridge University Press}, author={Horn, Roger A. and Johnson, Charles R.}, year={1985}}

@article{gaussianbound,
author = {Massart, P.},
year = {2003},
pages = {},
title = {Concentration Inequalities and Model Selection},
volume = {1896},
isbn = {978-3-540-48497-4},
journal = {Lecture Notes in Mathematics, Springer},
doi = {10.1007/978-3-540-48503-2}
}

@article{pinnau2017,
  author  = {Pinnau, R. and Totzeck, C. and Tse, O. and Martin, S.},
  title   = {A consensus-based model for global optimization and its mean-field limit},
  journal = {Mathematical Models and Methods in Applied Sciences},
  volume  = {27},
  number  = {1},
  pages   = {183--204},
  doi     = {10.1142/S0218202517400061},
  year = {2017}
}

@article{carrillo2018,
  author  = {Carrillo, J. A. and Choi, Y.-P. and Totzeck, C. and Tse, O.},
  title   = {An analytical framework for a consensus-based global optimization method},
  journal = {Mathematical Models and Methods in Applied Sciences},
  volume  = {28},
  number  = {6},
  pages   = {1037--1066},
  year    = {2018}
}

@article{fornasier2024,
  author  = {Fornasier, M. and Klock, T. and Riedl, K.},
  title   = {Consensus-based optimization methods converge globally},
  journal = {SIAM Journal on Optimization},
  volume  = {34},
  number  = {3},
  pages   = {2973--3004},
  year    = {2024},
  doi     = {10.1137/22M1527805}
}

@article{hajinkim2020,
  author  = {Ha, S.-Y. and Jin, S. and Kim, D.},
  title   = {Convergence of a first-order consensus-based global optimization algorithm},
  journal = {Mathematical Models and Methods in Applied Sciences},
  volume  = {30},
  pages   = {2417--2444},
  year    = {2020}
}

@article{carrillojinlizhu2021,
  author  = {Carrillo, J. A. and Jin, S. and Li, L. and Zhu, Y.},
  title   = {A consensus-based global optimization method for high dimensional machine learning problems},
  journal = {ESAIM: Control, Optimisation and Calculus of Variations},
  volume  = {27},
  pages   = {S5},
  year    = {2021}
}

@article{fornasierrichtarikriedlsun2025,
  author  = {Fornasier, M. and Richt{\'a}rik, P. and Riedl, K. and Sun, L.},
  title   = {Consensus-based optimization with truncated noise},
  journal = {European Journal of Applied Mathematics},
  volume  = {36},
  number  = {2},
  pages   = {292--315},
  year    = {2025}
}

@article{kalisesharmatretyakov2023,
  author  = {Kalise, D. and Sharma, A. and Tretyakov, M. V.},
  title   = {Consensus-based optimization via jump-diffusion stochastic differential equations},
  journal = {Mathematical Models and Methods in Applied Sciences},
  volume  = {33},
  number  = {2},
  pages   = {289--339},
  year    = {2023}
}

@article{riedl2023memory,
  author  = {Riedl, K.},
  title   = {Leveraging memory effects and gradient information in consensus-based optimisation: On global convergence in mean-field law},
  journal = {European Journal of Applied Mathematics},
  pages   = {1--32},
  note    = {First View},
  year    = {2023}
}

@article{ko2022,
  author  = {Ko, D. and Ha, S.-Y. and Jin, S. and Kim, D.},
  title   = {Convergence analysis of the discrete consensus-based optimization algorithm with random batch interactions and heterogeneous noises},
  journal = {Mathematical Models and Methods in Applied Sciences},
  volume  = {32},
  number  = {06},
  pages   = {1071--1107},
  year    = {2022}
}

@article{bellavia2025,
  author  = {Bellavia, S. and Malaspina, G.},
  title   = {A discrete consensus-based global optimization method with noisy objective function},
  journal = {Journal of Optimization Theory and Applications},
  volume  = {206},
  number  = {20},
  year    = {2025}
}

@article{totzeckwolfram2020,
  author  = {Totzeck, C. and Wolfram, M.-T.},
  title   = {Consensus-based global optimization with personal best},
  journal = {Mathematical Biosciences and Engineering},
  volume  = {17},
  number  = {5},
  pages   = {6026--6044},
  year    = {2020}
}

@article{fornasier2020hypersurfaces,
  author  = {Fornasier, M. and Huang, H. and Pareschi, L. and S{\"u}nnen, P.},
  title   = {Consensus-based optimization on hypersurfaces: well-posedness and mean-field limit},
  journal = {Mathematical Models and Methods in Applied Sciences},
  volume  = {30},
  pages   = {2725--2751},
  year    = {2020}
}

@article{fornasier2021spherejmlr,
  author  = {Fornasier, M. and Huang, H. and Pareschi, L. and S{\"u}nnen, P.},
  title   = {Consensus-based optimization on the sphere: Convergence to global minimizers and machine learning},
  journal = {Journal of Machine Learning Research},
  volume  = {22},
  number  = {237},
  pages   = {1--55},
  year    = {2021}
}

@article{carrillohoffmannstuartvaes2022,
  author  = {Carrillo, J. A. and Hoffmann, F. and Stuart, A. M. and Vaes, U.},
  title   = {Consensus-based sampling},
  journal = {Studies in Applied Mathematics},
  volume  = {148},
  number  = {3},
  pages   = {1069--1140},
  year    = {2022}
}

@INPROCEEDINGS{sun2026robotics,
    AUTHOR    = {Xudong Sun AND Armand Jordana AND Massimo Fornasier AND Jalal Etesami AND Majid Khadiv},
    TITLE     = {{Consensus-based optimization (CBO): Towards Global Optimality in Robotics}},
    BOOKTITLE = {Proceedings of Robotics: Science and Systems},
    YEAR      = {2026},
    ADDRESS   = {Sydney, Australia},
    MONTH     = {July},
    DOI       = {10.15607/RSS.2026.XXII.161}
}

@misc{fornasierhuangklemencmalaspina2026,
  author        = {Fornasier, M. and Huang, H. and Klemenc, J. and Malaspina, G.},
  title         = {From Consensus-Based Optimization to Evolution Strategies: Proof of Global Convergence},
  year          = {2026},
  eprint        = {2602.11677},
  archivePrefix = {arXiv},
  note          = {arXiv preprint}
}

@incollection{totzeck2021trends,
  author    = {Totzeck, C.},
  title     = {Trends in Consensus-Based Optimization},
  booktitle = {Active Particles, Volume 3: Advances in Theory, Models, and Applications},
  pages     = {201--226},
  publisher = {Springer},
  year      = {2021}
}

@article{carrillogarciatrillos2024fedcbo,
  author  = {Carrillo, J. A. and Garc{\'i}a Trillos, N. and Li, S. and Zhu, Y.},
  title   = {{FedCBO}: Reaching group consensus in clustered federated learning through consensus-based optimization},
  journal = {Journal of Machine Learning Research},
  volume  = {25},
  number  = {214},
  pages   = {1--51},
  year    = {2024}
}

@misc{garciatrilloscb2o2024,
  author        = {Garc{\'i}a Trillos, N. and Li, S. and Riedl, K. and Zhu, Y.},
  title         = {{CB\textsuperscript{2}O}: Consensus-Based Bi-Level Optimization},
  year          = {2024},
  eprint        = {2411.13394},
  archivePrefix = {arXiv},
  note          = {arXiv preprint}
}

@article{garciatrillosdefending2025,
  author  = {Garc{\'i}a Trillos, N. and Kumar Akash, A. and Li, S. and Riedl, K. and Zhu, Y.},
  title   = {Defending against diverse attacks in federated learning through consensus-based bi-level optimization},
  journal = {Philosophical Transactions of the Royal Society A},
  volume  = {383},
  number  = {2298},
  pages   = {20240235},
  year    = {2025}
}

@inproceedings{geiping2020invertinggradients,
  author    = {Geiping, J. and Bauermeister, H. and Dr{\"o}ge, H. and Moeller, M.},
  title     = {Inverting gradients -- how easy is it to break privacy in federated learning?},
  booktitle = {Advances in Neural Information Processing Systems},
  volume    = {33},
  pages     = {16937--16947},
  year      = {2020}
}

@inproceedings{haim2022reconstructing,
  author    = {Haim, N. and Vardi, G. and Yehudai, G. and Shamir, O. and Irani, M.},
  title     = {Reconstructing training data from trained neural networks},
  booktitle = {Advances in Neural Information Processing Systems},
  volume    = {35},
  pages     = {22911--22924},
  year      = {2022}
}

@article{xiaoboyd2004,
  author  = {Xiao, L. and Boyd, S.},
  title   = {Fast linear iterations for distributed averaging},
  journal = {Systems \& Control Letters},
  volume  = {53},
  number  = {1},
  pages   = {65--78},
  year    = {2004}
}

@article{nedicozdaglar2009,
  author  = {Nedi{\'c}, A. and Ozdaglar, A.},
  title   = {Distributed subgradient methods for multi-agent optimization},
  journal = {IEEE Transactions on Automatic Control},
  volume  = {54},
  number  = {1},
  pages   = {48--61},
  year    = {2009}
}

@article{Bailo_CBX_Python_and_2024,
author = {Bailo, Rafael and Barbaro, Alethea and Gomes, Susana N. and Riedl, Konstantin and Roith, Tim and Totzeck, Claudia and Vaes, Urbain},
doi = {10.21105/joss.06611},
journal = {Journal of Open Source Software},
month = jun,
number = {98},
pages = {6611},
title = {{CBX: Python and Julia Packages for Consensus-Based Interacting Particle Methods}},
url = {https://joss.theoj.org/papers/10.21105/joss.06611},
volume = {9},
year = {2024}
}

@article{bungert2025mirrorcbo,
  title={{MirrorCBO}: A consensus-based optimization method in the spirit of mirror descent},
  author={Bungert, Leon and Hoffmann, Franca and Kim, Doh Yeon and Roith, Tim},
  journal={Mathematical Models and Methods in Applied Sciences},
  volume={35},
  number={14},
  pages={3083--3170},
  year={2025},
  publisher={World Scientific}
}

@inproceedings{ros2008simple,
  title={A simple modification in CMA-ES achieving linear time and space complexity},
  author={Ros, Raymond and Hansen, Nikolaus},
  booktitle={International conference on parallel problem solving from nature},
  pages={296--305},
  year={2008},
  organization={Springer}
}

@article{jamil2013literature,
  title={A literature survey of benchmark functions for global optimisation problems},
  author={Jamil, Momin and Yang, Xin-She},
  journal={International Journal of Mathematical Modelling and Numerical Optimisation},
  volume={4},
  number={2},
  pages={150--194},
  year={2013},
  publisher={Inderscience Publishers Ltd}
}

@misc{nikolaus_hansen_2026_22738161,
  author       = {Nikolaus Hansen and
                  yoshihikoueno and
                  ARF1 and
                  Sait Cakmak and
                  Gabi Kadlecová and
                  Guillermo Abad López and
                  Kento Nozawa and
                  Luca Rolshoven and
                  martinsvardsjo and
                  Paul Stahr and
                  Youhei Akimoto and
                  brieglhostis and
                  Dimo Brockhoff},
  title        = {CMA-ES/pycma: r4.5.0},
  year         = {2026},
  publisher    = {Zenodo},
  version      = {r4.5.0},
  doi          = {10.5281/zenodo.22738161},
}

@article{hansen2006cma,
  title={The CMA evolution strategy: a comparing review},
  author={Hansen, Nikolaus},
  journal={Towards a new evolutionary computation: Advances in the estimation of distribution algorithms},
  pages={75--102},
  year={2006},
  publisher={Springer}
}

@article{fiedler2026consensus,
  title={Consensus-based optimization for linearly separable functions},
  author={Fiedler, Christian and Roith, Tim},
  journal={arXiv preprint arXiv:2609.01317},
  year={2026}
}

@article{beddrich2026constrained,
  title={Constrained consensus-based optimization and numerical heuristics for the few particle regime},
  author={Beddrich, Jonas and Chenchene, Enis and Fornasier, Massimo and Huang, Hui and Wohlmuth, Barbara},
  journal={Journal of Global Optimization},
  pages={1--54},
  year={2026},
  publisher={Springer}
}

@article{hansen2008method,
  title={A method for handling uncertainty in evolutionary optimization with an application to feedback control of combustion},
  author={Hansen, Nikolaus and Niederberger, Andr{\'e} SP and Guzzella, Lino and Koumoutsakos, Petros},
  journal={IEEE Transactions on Evolutionary Computation},
  volume={13},
  number={1},
  pages={180--197},
  year={2008},
  publisher={IEEE}
}

@phdthesis{hansen2009real,
  title={Real-parameter black-box optimization benchmarking 2009: Noiseless functions definitions},
  author={Hansen, Nikolaus and Finck, Steffen and Ros, Raymond and Auger, Anne},
  year={2009},
  school={INRIA}
}

@article{lecun2010mnist,
  title={{MNIST} handwritten digit database},
  author={LeCun, Yann and Cortes, Corinna and Burges, CJ},
  journal={ATT Labs [Online]. Available: http://yann.lecun.com/exdb/mnist},
  volume={2},
  year={2010}
}

@misc{bonandin2026exploiting,
  title         = {Exploiting Structure with Anisotropic Consensus-Based Optimization},
  author        = {Bonandin, Sabrina and Riedl, Konstantin and Veneruso, Sara},
  year          = {2026},
  eprint        = {2607.10205},
  archivePrefix = {arXiv},
  primaryClass  = {math.OC},
  url           = {https://arxiv.org/abs/2607.10205}
}

@article{benfenati2022binary,
  author  = {Benfenati, Alessandro and Borghi, Giacomo and Pareschi, Lorenzo},
  title   = {Binary interaction methods for high dimensional global
             optimization and machine learning},
  journal = {Applied Mathematics and Optimization},
  volume  = {86},
  number  = {1},
  pages   = {Paper No. 9},
  year    = {2022},
  doi     = {10.1007/s00245-022-09836-5}
}

@article{bellavia2026parallelCBO,
  author  = {Bellavia, Stefania and Fornasier, Massimo and Malaspina, Greta and Roith, Tim},
  title   = {A parallelization strategy for consensus based optimization with global convergence guarantees},
  journal = {In preparation},
  year    = {2026},
}
\end{document}